\documentstyle{article}

\newcommand{\N}{{\rm I\kern-.5ex N}}
\newcommand{\Z}{{\sf \vrule height 1.55ex depth-1.2ex width.03em\kern-.11em 
Z \kern-.9ex Z\kern-.11em\vrule height 0.3ex depth0ex width.03em}}
\newcommand{\Q}{{\rm\kern.2ex\vrule height1.55ex depth-.05ex width.03em\kern-.7ex Q}}
\newcommand{\R}{{\rm I\kern-.5ex R}}
\newcommand{\Rvar}{{\rm I\kern-.5ex R}}
\newcommand{\C}{{\rm\kern.3ex\vrule height1.55ex depth-.05ex width.03em\kern-.7ex C}}
\newcommand{\Cvar}{{\, \rm\kern.3ex\vrule height1.1ex depth-.05ex width.03em\kern-.7ex C}}

\newcommand{\nabp}{\nab \hspace{-1.05ex}
\rule[.5ex]{.2ex}{.8ex}   \hspace{1.05ex}}

\newcommand{\spat}{\hspace{4ex}}
\newcommand{\restriction}{\!\mid\!}

\newcommand{\cu}{{\cal U}}

\newcommand{\nab}{\nabla}

\newcommand{\cF}{{\cal F}}
\newcommand{\cG}{{\cal G}}

\newcommand{\cN}{{\cal N}}
\newcommand{\cM}{{\cal M}}

\newcommand{\od}{\odot}
\newcommand{\ot}{\otimes}

\newcommand{\la}{\Lambda}

\newcommand{\om}{\omega}
\newcommand{\io}{\iota}
\newcommand{\vfi}{\varphi}
\newcommand{\vep}{\varepsilon}

\newcommand{\al}{\alpha}

\newcommand{\ga}{\Gamma}
\newcommand{\sde}{\delta}
\newcommand{\de}{\Delta}

\newcommand{\th}{\theta}

\newcommand{\si}{\sigma}
\newcommand{\Mfi}{{\cal M}_{\vfi}}
\newcommand{\Nfi}{{\cal N}_{\vfi}}

\newcommand{\Mpsi}{{\cal M}_{\psi}}
\newcommand{\Npsi}{{\cal N}_{\psi}}
\newcommand{\lafi}{\la_\vfi}

\newcommand{\oal}{\overline{\al}}

\newcommand{\text}[1]{\mbox{#1}}

\newcommand{\cst}{\text{C}$^*$}

\newcommand{\qed}{\ \hfill \rule{2mm}{2mm}}
\newenvironment{demo}{\medskip\noindent\bf Proof :\ \  \rm}{\qed\bigskip\par }

\newtheorem{definition}{Definition}[section]
\newtheorem{proposition}[definition]{Proposition}
\newtheorem{lemma}[definition]{Lemma}
\newtheorem{corollary}[definition]{Corollary}
\newtheorem{remark}[definition]{Remark}
\newtheorem{theorem}[definition]{Theorem}

\newtheorem{notation}[definition]{Notation}
\newtheorem{result}[definition]{Result}

\newtheorem{terminology}[definition]{Terminology}

\begin{document}
\begin{center}
\LARGE\bf
KMS-weights on \cst-algebras
\end{center}

\bigskip

\begin{center}
\rm J. Kustermans  \footnote{Research Assistant of the
National Fund for Scientific Research (Belgium)}

Institut for Matematik og Datalogi

Odense Universitet

Campusvej 55

5230 Odense M

Denmark

e-mail : johank@imada.ou.dk

\bigskip\medskip

\bf April 1997 \rm
\end{center}

\bigskip

\subsection*{Abstract}
In this paper, we build a solid framework for the use of KMS-weights on \cst-algebras. We
will use another definition than the one introduced by Combes in \cite{Comb2}, but we will
prove that they are equivalent. However, the subject of KMS-weights is approached from a
somewhat different angle. We introduce a construction procedure for KMS-weights, prove the
most important properties of them, construct the tensor product of KMS-weights and
construct weights which are absolutely continuous to a given weight.

\bigskip

\section*{Introduction}

In \cite{Comb}, Combes studied lower semi-continuous weights on \cst-algebras and proved
the major result that a lower semi-continuous weight can be approximated by the positive
linear functionals which are majorated by it.

At the moment, these lower semi-continuous weights seem not to behave well enough to be
used in an efficient way. So it is natural to search for  extra conditions conditions on
lower semi-continuous weights which makes them manageable. An attempt in this direction
was undertaken by Jan Verding in  \cite{Verd} where he introduced the so-called regular
weights. These are lower semi-continuous weights which have a well-behaved truncating net.

It turns out that this regularity condition is also very useful in the theory of
\cst-valued weights (see \cite{JK2}).

\medskip

There is however a class of lower semi-continuous weights which behave very well and were
introduced by Combes in definition 4.1 of \cite{Comb2}, the so-called KMS-weights:

\begin{list}{}{\setlength{\leftmargin}{.4 cm}}

\item Consider a \cst-algebra $A$ and a densely defined lower semi-continuous weight
$\vfi$ on $A$ such that there exists a norm continuous one-parameter group $\si$ on $A$
satisfying:
\begin{enumerate}
\item We have for every $t \in \R$ that $\vfi \si_t = \vfi$.
\item For every $x,y \in \Nfi \cap \Nfi^*$, there exists a bounded continuous function $f$
from $S(i)$ into $\,\C$ which is analytic on $S(i)^0$ and such that:
\begin{itemize}
\item We have for every $t \in \R$ that $f(t) = \vfi(\si_t(x) y)$.
\item We have for every $t \in \R$ that $f(t+i) = \vfi(y  \si_t(x))$.
\end{itemize}
\end{enumerate}
Then $\vfi$ is called a KMS-weight on $A$, $\si$ is called a modular group for $\vfi$.

\end{list}

Here, $S(i)$ denotes the horizontal strip in the complex plane between 0 and i.

\medskip

By a proof of Combes in \cite{E-V}, we know that such a KMS-weight can be extended to
normal KMS-weight in the GNS-representation of $\vfi$. In this way, we can make use of the
very rich theory of normal weights on von Neumann algebras and use for instance the Radon
Nikodym theorem for normal weights (see \cite{Pe-Tak}).

\medskip

In short, we get a class of weights which which posses the same rich structure as the
Borel measures on locally compact spaces.

\bigskip

The last years, there is a lot of interest for lower-semicontinuous weights in general
(and KMS-weights in particular) from people working in the field of \cst-algebraic quantum
groups. This because of the role of the left Haar weight in \cst-algebraic quantum group
theory.

At the moment, it is not yet clear whether the left Haar weight should be a KMS-weight.
There are however some indications that this will be the case :
\begin{itemize}
\item Masuda, Nakagami \& Woronowicz are momentarily working on a definition for
reduced \cst-algebraic quantum group and their left Haar weight will be a KMS-weight.
\item In all the known examples, the left Haar weight turns out to satisfy the KMS-condition.
\item A. Van Daele proves in \cite{VD1} that the left Haar functional on an algebraic quantum group
satisfies automatically some sort of weak KMS-condition.
\item It is proven in \cite{Kus} that this weak KMS-condition implies that the left Haar weight
on the reduced \cst-algebraic group arising from the algebraic one, is a KMS-weight.
\item In \cite{Kus1}, we prove that this KMS-condition on this reduced \cst-algebraic quantum group
implies that the left Haar weight on the universal \cst-algebraic quantum group is also
KMS. It appears that the same principle can be used to prove such a result in the general
cas.
\end{itemize}

Especially the universal case necessitates to investigate KMS-weights in a \cst-algebraic
framework.

\bigskip

In this paper, we will work with another definition of a KMS-weight than the one mentioned
above but we will show that both are equivalent. One of the reasons to use the other
definition is the fact that Masuda, Nakagami \& Woronowicz seemed to be going to use this
other definition in their approach to quantum groups.

\medskip

The main aim of this paper is to build a solid framework for KMS-weights on \cst-algebras.
We believe that it gives a technically useful overview of KMS-weights in the \cst-algebra
picture.

Most of the results will look familiar to the ones known in normal weight theory but will
be proven  within the \cst-algebra framework.

\bigskip

In the first section, we fix notations and give an overview of the results about
one-parameter representations and their analytic continuations.

We give the definition of a KMS-weight in the second section. At the same time, we prove
some useful properties concerning lower semi-continuous weights.

In section 3, we have gathered some results from \cite{Verd} about the construction of a
weight starting from a GNS-construction.

The  fourth section concerns the behaviour of one-parameter representations which are
relatively invariant with respect to some closed mapping.

In the fifth section, we prove that a closed linear mapping  from within a
\cst-algebra into a Hilbert space and  satisfying certain KMS-characteristics
gives rise to a KMS-weight on this \cst-algebra.

The most important technical properties concerning KMS-weights are collected in the sixth
section.

In section 7, we construct the tensor product of KMS-weights.

Given a KMS-weight $\vfi$ and a strictly positive element $\sde$ affiliated to the
\cst-algebra which is relatively invariant under a modular group of $\vfi$, we construct
the weight $\vfi(\sde^\frac{1}{2} \, . \, \sde^\frac{1}{2})$ in section 8. We have to use
another method than the one used by Pedersen \& Takesaki because of the relative
invariance.

\bigskip

Let us now fix some terminology and conventions.

\medskip

Consider a complex number $z$, then $S(z)$ will denote the horizontal strip  $\{
\, c \in \C \mid \text{Im }c \in [0,\text{Im }z] \, \}$.  The interior of $S(z)$ will be
denoted by $S(z)^0$.

\medskip

The domain of a mapping $f$ will be denoted by $D(f)$, its image by $\text{Ran}\,f$.

\medskip

For any Banach space $E$, we denote the set of bounded operators on $E$ by $B(E)$ and the
topological dual by $E^*$.

\medskip

Let $A$ be a \cst-algebra. We denote the multiplier algebra of $A$ by $M(A)$.

Any element $\om \in A^*$ has a unique strictly continuous linear extension to $M(A)$ and
we denote this extension by $\overline{\om}$. We put $\om(x) = \overline{\om}(x)$ for
every $x \in M(A)$.

\medskip

For a good introduction into elements affiliated with $A$, we refer to \cite{Lan} and
\cite{Wor6}.

Whenever we speak of a positive element affiliated element with $A$, we mean a (possibly
unbounded) positive and selfadjoint element affiliated with $A$. The same remark applies
to positive operators in Hilbert spaces. A positive element affiliated to $A$ is called
strictly positive if it has dense range. As usual, we can raise a strictly positive
affiliated element to any complex power.

\medskip

Consider a Hilbert space $H$. Let $\pi$ a non-degenerate $^*$-homomorphism  from $A$ into
$B(H)$. Then $\pi$ has a unique extension to a $^*$-representation $\overline{\pi}$. We
put $\pi(x) = \overline{\pi}(x)$ for every $x \in M(A)$.

If $T$ is an element affiliated with $A$, then $\pi(T)$ denotes the closed densely defined
operator in $H$ such that $\pi(D(T))$ is a core for $\pi(T)$ and $\pi(T)(\pi(a)v) =
\pi(T(a))v$ for every $a \in D(T)$ and $v \in H$.

\medskip

In the following, $\circ$ will denote the composition of mappings and elements of $M(A)$
will be considered as linear mappings from $A$ into $A$.

\bigskip

Consider an element $T$ affiliated with $A$ and $a \in M(A)$. Then
\begin{itemize}
\item We call $a$ a left multiplier of $T$ if $a \! \circ \! T$ is bounded
as a mapping from $A$ into $A$ and $\overline{a \, \circ \, T}$ belongs to $M(A)$.  In
this case, we put  $a \, T = \overline{a \! \circ \! T}$
\item We call $a$ a right multiplier of $T$ if $D(T \! \circ \! a) = A$. You can prove
in this case that $T \! \circ \! a$ belongs to $M(A)$ and we put $T \, a = T \! \circ \!
a$.
\end{itemize}

It is not very difficult to prove that $a$ is a left multiplier of $T$ if and only if $a$
is a right multiplier of $T^*$. If this is the case, we have that $(a T)^* = T^* a^*$.

We will freely use calculation rules involving this left and right multipliers, e.g. Let
$a,x \in M(A)$ and $T$ an element affiliated with $A$. If $a$ is a left multiplier of $T$,
then $x a$ is a left multiplier of $T$ and $(x a) \, T = x \, (a \, T)$.

\medskip

Consider elements $S$, $T$ affiliated with $A$ and $a \in M(A)$. Then we call $a$ a middle
multiplier of $S,T$ if  $D(S \!\circ\! a \!\circ\! T) = D(T)$, $S\!\circ\!a\!\circ\! T$ is
bounded and $\overline{S\!\circ\! a\!\circ\! T}$ belongs to $M(A)$. In this case, we put
$S a T = \overline{S \!\circ\! a \!\circ\!  T}$.

\bigskip

\section{One-parameter representations}

In this first section, we  recall the definition of one-parameter representations and an
overview of the most important results. Most of the proofs can be found in \cite{JK3}. A
standard reference for one-parameter representationd is \cite{Zsido}.

\bigskip

Whenever we use the notion of integrability of a function with values in a Banach space,
we mean the strong form of integrability (e.g. Analysis II, S. Lang) :

\medskip

Consider a measure space $(X,{\cal M},\mu)$ and a Banach space $E$.

\begin{itemize}
\item It is obvious how to define integrability for step functions from $X$ into $E$.
\item Let $f$ be a function from $X$ into $E$. Then $f$ is $\mu$-integrable if and only if there exists a sequence of integrable step functions $(f_n)_{n=1}^\infty$ from $X$ into $E$ such that :
\begin{enumerate}
\item We have for almost every $x \in X$ that $(f_n(x))_{n=1}^\infty$ converges to $f(x)$ in the norm topology.
\item The sequence $(f_n)_{n=1}^\infty$ is convergent in the $L_1$-norm.
\end{enumerate}

In this case, the sequence $(\int f_n \, d\mu)_{n=1}^\infty$ is convergent and the
integral of $\int f \, d\mu$ is defined to be the limit of this sequence (Of course, one
has to prove that this limit is independent of the choice of the sequence
$(f_n)_{n=1}^\infty$).
\item It is possible to define a form of $\mu$-measurability (see Lang) in such a way
that a function $f$ is integrable if and only if the function $f$ is measurable and
$\|f\|$ is integrable.
\end{itemize}

\bigskip

We start this section with a basic but very useful lemma.

\begin{lemma} \label{o.lem1}
Consider Banach spaces $E$, $F$, \  $\la$ a closed linear mapping from within $E$ into
$F$. Let $f$ be a function from $\R$ into $D(\la)$ such that
\begin{itemize}
\item $f$ is integrable.
\item The function $\R \rightarrow F : t \mapsto \la(f(t))$ is integrable.
\end{itemize}
Then $\int f(t) \, dt$ belongs to $D(\la)$ and $\la\bigl(\,\int f(t) \, dt\,\bigr)
= \int \la(f(t)) \, dt$.
\end{lemma}
\begin{demo}
Define $G$ as the graph of the mapping $\la$. By assumption, we have that $G$ is a closed
subspace of $E \oplus F$. Next, we define the mapping $g$ from $\R$ into $G$ such that
$g(t) = \bigl(f(t),\la(f(t))\bigr)$ for every $t \in \R$. It follows that $g$ is
integrable and
$$\int g(t) \, dt = \bigl( \, \int f(t) \, dt \, , \int \la(f(t)) \, dt \, \bigr) \ .$$
Because $G$ is a closed subspace of $E \oplus F$, we have that $\int g(t) \, dt$ belongs
to $G$. This implies that $\bigl(\, \int f(t) \, dt \, , \int \la(f(t)) \, dt\,\bigr)$
belongs to $G$.
\end{demo}

\bigskip

When we speak of analyticity, we will always mean norm analyticity. But we have the
following well known result (which follows from the Uniform Boundedness principle).

\begin{result} \label{anal}
Consider a Banach space $E$ and a subspace $F$ of $E^*$ such that $\|x\| = \sup\{ \,
\|\om(x)\| \mid \om \in F \text{ with } \|\om\| \leq 1 \}$ \ \ for every $x \in E$. Let
$O$ be an open subset of $\,\C$ and $f$ a function from $O$ into $E$.

Then $f$ is analytic $\Leftrightarrow$ We have for every $\om \in F$ that $\om \! \circ \!
f$ is analytic.
\end{result}

Notice that the result is true for $F=E^*$.

\subsection{One-parameter representations on Banach spaces}

We introduce the notion of strongly continuous one-parameter groups and discuss the
analytic continuations of them. The standard reference for the material in this section is
\cite{Zsido}.

\begin{terminology}
Consider a Banach space $E$.
By a one-parameter representation on $E$, we will always mean a mapping
$\al$ from $\R$ into $B(E)$ such that :
\begin{itemize}
\item We have for every $s,t \in \R$ that $\al_{s+t} = \al_s \al_t$.
\item $\al_0= \io$
\item We have for every $t \in \R$ that $\| \al_t \| \leq 1$.
\end{itemize}
We call $\al$  strongly continuous $\Leftrightarrow$ We have for every $a \in E$ that the
mapping $\R \rightarrow E : t \mapsto \al_t(a)$ is continuous.
\end{terminology}

This definition implies for every $t \in \R$ that $\al_t$ is invertible in $B(E)$, that
$(\al_t)^{-1} = \al_{-t}$ and that $\|\al_t\|=1$.

\begin{remark}\rm We would like to mention the following special cases :
\begin{itemize}
\item If $H$ is a Hilbert space and $u$ is a strongly continuous one-parameter
representation on $H$ such that $u_t$ is unitary for every $t \in \R$, we call $u$ a
strongly continuous unitary one-parameter group on $H$.
\item If $A$ is a \cst-algebra and $\al$ is a strongly continuous one-parameter representation
on $A$ such that $\al_t$ is a \cst-automorphism on $A$ for every $t \in \R$, we call $\al$
a norm continuous one-parameter group on $A$.
\item If $A$ is a \cst-algebra and $u$ is a strongly continous one-parameter representation
on $A$ such that $u_t$ is a unitary element in $M(A)$ for every $t \in \R$, we call $u$ a
strictly continuous unitary one-parameter group on $A$.
\end{itemize}
\end{remark}

Now we will summarize the theory of analytic continuations of such norm continuous
one-parameter representations. For the  most part of this subsection, we will fix a Banach
space $E$ and a strongly continuousone-parameter representation $\al$ on $E$.

\begin{definition} \label{pa1.def1}
Consider $z \in \C$. We define the mapping $\al_z$ from within $E$ into $E$
such that :
\begin{itemize}
\item The domain of $\al_z$ is by definition the set \
$\{\, a \in E \mid$ There exists a function $f$ from $S(z)$ into $E$ such that
\begin{enumerate}
\item $f$ is continuous on $S(z)$
\item $f$ is analytic on $S(z)^0$
\item We have  that $\al_t(a) = f(t)$ for every $t \in \R$ \ \ \ \ \ \ $\}$.
\end{enumerate}
\item Choose $a$ in the domain of $\al_z$ and $f$ the function from $S(z)$ into $E$ such that
\begin{enumerate}
\item $f$ is continuous on $S(z)$
\item $f$ is analytic on $S(z)^0$
\item We have  that $\al_t(a) = f(t)$ for every $t \in \R$
\end{enumerate}
Then we have by definition that $\al_z(a) = f(z)$.
\end{itemize}
\end{definition}

In the case that $z$ belongs to $\R$, this definition of $\al_z$ corresponds with $\al_z$
which we started from. Therefore, the notation is justified.

\begin{remark}\rm
\begin{itemize}
\item Consider $z \in \C$ and $y \in S(z)$. Then it is clear that $D(\al_z)$ is a
subset of $D(\al_y)$.
\item It is also clear that $D(\al_y)=D(\al_z)$ for $y,z \in \C$ with $\text{\rm Im }
 y = \text{\rm Im } z$.
\item Let $z \in \C$ and $a \in D(\al_z)$. Then the  function $S(z) \rightarrow E :
u \mapsto \al_u(a)$ is continuous on $S(z)$ and analytic on $S(z)^0$ (because this
function must be equal to the function $f$ from the definition).
\item Consider an element $a$ in $E$. We say that $a$ is analytic with respect to
$\al$ if $a$ belongs to $D(\al_z)$ for every $z \in \C$. If $a$ is analytic with respect
to $\al$, then the function $\C \rightarrow  E : u \mapsto \al_u(a)$ is analytic.
\end{itemize}
\end{remark}

Now we give a list of the most important properties of these analytic continuations. The
most important tool for proving the results is the Phragmen-Lindelof theorem.

\begin{proposition}
The mapping $\al_z$ is a closed linear operator in $E$ which is densely defined and has
dense range.
\end{proposition}

We would like to mention that the proof of the closedness is substantially easier than in
the case of strongly continuous one-parameter groups on von Neumann algebras (see
\cite{Zsido}).

\medskip

\begin{result}
Consider $z \in \C$ and $a \in D(\al_z)$. Then the function $S(z) \rightarrow E : u
\mapsto \al_u(a)$ is bounded.
\end{result}

\medskip

\begin{proposition}
Consider $z \in \C$ and $t \in \R$. Then $\al_z \al_t = \al_t \al_z = \al_{z+t}$.

Furthermore, the following equalities hold :
\begin{itemize}
\item We have that $\al_t(D(\al_z)) = D(\al_z)$ and $\al_t(\text{\rm Ran}\, \al_z) = \text{\rm
Ran } \al_z$.
\item We have that $D(\al_{z+t}) = D(\al_z)$ and $\text{\rm Ran}\, \al_{z+t} =$
$\text{\rm Ran}\, \al_z$.
\end{itemize}
\end{proposition}

\medskip

\begin{proposition}
Consider $z \in \C$. Then $\al_z$ is injective and $(\al_z)^{-1} = \al_{-z}$. Therefore
$\text{\rm Ran}\, \al_z = D(\al_{-z})$ and $D(\al_z) = \text{\rm Ran}\, \al_{-z}$.
\end{proposition}

\medskip

\begin{proposition}
Consider $y,z \in \C$. Then
\begin{enumerate}
\item $\al_y \, \al_z \subseteq \al_{y+z}$.
\item If $y$ and $z$ lie at the same side of the real axis, we have that
$\al_y \, \al_z = \al_{y+z}$.
\end{enumerate}
\end{proposition}

\medskip

As usual, we smear elements to construct elements which behave well with respect to $\al$.

\begin{notation}
Consider $a \in E$, $n > 0$ and $z \in \C$. Then we define the element $a(n,z)$ in $E$
such that
$$a(n,z) = \frac{n}{\sqrt{\pi}} \int \exp(-n^2 (t-z)^2) \, \al_t(a) \, dt
\ .$$
Then we have that $\|a(n,z)\| \leq \|a\| \, \exp(n^2 (\text{Im }z)^2)$.
\end{notation}

\medskip

For $n > 0$ and $z \in \C$, we will use the notation
$$ a(n) = a(n,0) = \frac{n}{\sqrt{\pi}} \int \exp(-n^2 t^2) \, \al_t(a) \, dt \ .$$

\medskip

\begin{proposition} \label{pa1.prop1}
Let $a \in E$, $n > 0$ and $z \in \C$. Then $a(n,z)$ is analytic  with respect to $\al$
and $\al_y(a(n,z)) = a(n,z+y)$ for every $y \in \C$.
\end{proposition}

\medskip

\begin{proposition}
Consider   $n > 0$, $y,z \in \C$ and $a \in D(\al_y)$. Then $\al_y(a(n,z))$  $=
\al_y(a)(n,z)$.
\end{proposition}

\medskip

A proof of the next result can be found in \cite{Wor5}.

\begin{proposition} \label{pa.prop1}
Let $z$  be a complex number. Consider a dense subset $K$ of $E$. Then  the set \newline
$\langle \, a(n) \mid a \in K , n \in \N \, \rangle $ is a core for $\al_z$.
\end{proposition}

\bigskip\bigskip

We end this subsection with some special cases :

\medskip

The first one is a familiar one and the proof can be found in corollary 9.21 of
\cite{Stra}.

\begin{proposition}
Consider a Hilbert space $H$ and an injective positive operator $M$ in $H$. Define the
strongly continuous unitary one-parameter group $u$ on $H$ such that $u_t = M^{it}$ for
every $t \in \R$. Then $u_z = M^{i z}$ for every $z \in \C$.
\end{proposition}

\medskip

The next result is a \cst-version of the previous one.

\begin{proposition}
Consider a \cst-algebra $A$ and a strictly positive element $\sde$ affiliated with  $A$.
Define the strictly continuous unitary one-parameter group $u$ on $A$ such that we have
for every $t \in \R$ that $u_t =  \sde^{i t}$ considered  as a left multiplier. Then $u_z
= \sde^{i z}$ for every $z \in \C$.
\end{proposition}

\medskip

The previous result implies immediately the following result.

\begin{proposition}
Consider a \cst-algebra $A$ and a strictly positive element $\sde$ affiliated with  $A$.
Define the strictly continuous unitary one-parameter group $u$ on $A$ such that we have
for every $t \in \R$ that $u_t =  \sde^{i t}$ considered  as a right multiplier. Let $z$
be a complex number and $a$ an element in $A$, then :
\begin{itemize}
\item The element $a$ belongs to $D(u_z)$ $\Leftrightarrow$ $a$ is a left multiplier of
$\,\sde^{i z}$ and $a \, \sde^z$ belongs to $A$.
\item If $a$ belongs to $D(u_z)$, then $u_z(a) = a \, \sde^{i z}$.
\end{itemize}
\end{proposition}

\medskip

The following proposition is also a \cst-version of a known result in Hilbert space theory
(see proposition 9.24 of \cite{Stra}).

\begin{proposition}
Consider a \cst-algebra $A$ and a strictly positive element $\sde$ affiliated with  $A$.
Define the  norm continuous one-parameter group $\al$ on $A$ such that $\al_t(a) =
\sde^{-it} a \, \sde^{i t}$ for every $t \in \R$ and $a \in A$.

Let $z$ be a complex number and $a$ an element in $A$, then :
\begin{itemize}
\item The element $a$ belongs to $D(\al_z)$ $\Leftrightarrow$ $a$ is a middle multiplier of
$\sde^{-i z},\, \sde^{i z}$ and $\sde^{-i z} a \, \sde^{i z}$ belongs to $A$.
\item If $a$ belongs to $D( \al_z)$, then $\al_z(a) = \sde^{-i z} a \, \sde^{i z}$.
\end{itemize}
\end{proposition}

\medskip

A proof of the last 3 propositions can be found in \cite{JK3}.

\bigskip

\subsection{One-parameter groups on \cst-algebras}

For this section, we will fix a norm-continuous one-parameter group $\al$ on a
\cst-algebra $A$. We will investigate the properties of the analytic continuations in this
case a little bit further. Proofs of the results in this section can be found in
\cite{JK3}.

\medskip

By the remarks of the previous section , we have for any $z \in \C$ a closed linear
operator $\al_z$ from within $A$ into $A$ which is generally unbounded.

\medskip

Because $\al_t$ is a $^*$-homomorphism for every $t \in \R$, we have the  two following
algebraic properties.

\begin{proposition}
The mapping $\al_z$ is a multiplicative linear operator in $A$.
\end{proposition}

\begin{proposition}
Consider $z \in \C$. Then $D(\al_{\,\overline{z}}) = D(\al_z)^*$ and
$\text{\rm Ran}\, \al_{\,\overline{z}} = (\text{\rm Ran}\, \al_z)^*$.
We have for every $a \in D(\al_z)$ that $\al_z(a)^* = \al_{\,\overline{z}}(a^*)$.
\end{proposition}

This implies for all $z \in \R i$ and $a \in D(\al_z)$ that $\al_z(a)^*$ belongs to
$D(\al_z)$ and $\al_z(\al_z(a)^*)^* = a$.

\bigskip

For the rest of this section, we want to concentrate on strictly analytic continuations of
$\al$.

\medskip

Concerning analyticity, we have the following result. The proof of this fact uses result
\ref{anal} and the fact that every continuous linear functional on $A$ is of the form $a
\, \om$ with $\om \in A^*$ and $a \in A$.

\begin{result}
Consider an open subset $O$ of the complex plane, $f$ a function from $O$ into $M(A)$.

Then $f$ is analytic $\Leftrightarrow$ We have for every $\om \in A^*$ that the function
$\overline{\om} \! \circ \! f$ is analytic \newline $\Leftrightarrow$ We have for every $a
\in A$ that the function $O \rightarrow A : z \mapsto f(z) \, a$ is analytic.
\end{result}

\bigskip

It is possible to give the following definition.

\begin{definition}
Consider a complex number $z$. Then the mapping $\al_z$ is closable for the strict
topology on $M(A)$ and we denote the strict closure by $\overline{\al}_z$.
\end{definition}

In the case where $z$ belongs to $\R$, the previous  definition implies that  $\oal_z$ is
the unique $^*$-homomorphism from $M(A)$ into $M(A)$ which extends $\al_z$. It is also
clear that $\al_z \subseteq \oal_z$ for every $z \in \C$.

The two previous remarks justify the following notation.
For every $z \in \C$ and $a \in D(\oal_z)$, we put $\al_z(a) = \oal_z(a)$.

\bigskip

\begin{result} \label{pa2.prop1}
Let $a$ be an element in $M(A)$. Then the function $\R \rightarrow
M(A) :  t \mapsto \al_t(a)$ is strictly continuous.
\end{result}

\medskip

The following theorem is one of the most important results of \cite{JK3}.

\begin{theorem}
Consider $z \in \C$. Then we have the following properties
\begin{itemize}
\item Let $a$ be an element in $M(A)$.

\vspace{2mm}

Then $a$ belongs to $D(\,\oal_z)$ $\Leftrightarrow$ there exists a function $f$ from
$S(z)$ into $M(A)$ such that
\begin{enumerate}
\item $f$ is strictly continuous on $S(z)$
\item $f$ is analytic on $S(z)^0$
\item We have  that $\al_t(a) = f(t)$ for every $t \in \R$
\end{enumerate}
\item Let $a$ be an element in  $D(\,\oal_z)$ and $f$ the function from $S(z)$ into $M(A)$ such that
\begin{enumerate}
\item $f$ is strictly continuous on $S(z)$
\item $f$ is analytic on $S(z)^0$
\item We have  that $\al_t(a) = f(t)$ for every $t \in \R$
\end{enumerate}
Then we have that $\al_z(a) = f(z)$.
\end{itemize}
\end{theorem}

\bigskip

It is also possible to give another characterization of $\oal_z$. First look at the
following result.

\begin{proposition}
Consider $z \in \C$ and $a \in D(\oal_z)$ , $b \in D(\al_z)$. Then $a b$ and $b a$ belong
to $D(\al_z)$ and $\al_z(a b) = \al_z(a) \, \al_z(b)$ and $\al_z(b a) = \al_z(b) \,
\al_z(a)$.
\end{proposition}

There is also a converse of this :

\begin{proposition}
Consider elements $a,b \in M(A)$ and $z \in \C$. Then
\begin{enumerate}
\item If we have for every $c \in D(\al_z)$ that $a c$ belongs to $D(\al_z)$ and $\al_z(a c) = b \, \al_z(c)$,
then $a$ belongs to $D(\oal_z)$ and $\al_z(a) = b$.
\item If we have for every $c \in D(\al_z)$ that $c a$ belongs to $D(\al_z)$ and $\al_z(c a) = \al_z(c) \, b$,
then $a$ belongs to $D(\oal_z)$
and $\al_z(a) = b$.
\end{enumerate}
\end{proposition}

\bigskip

These strictly analytic extensions satisfy the same kind of properties as the norm
analytic extensions.

Now we give a list of these properties.

\medskip

\begin{remark}\rm
\begin{itemize}
\item Consider $z \in \C$ and $y \in S(z)$. Then it is clear that $D(\oal_z)$ is a
subset of $D(\oal_y)$.
\item It is also clear that $D(\oal_y)=D(\oal_z)$ for $y,z \in \C$ with $\text{Im }
y = \text{Im } z$.
\item Let $z \in \C$ and $a \in D(\oal_z)$. Then the  function $S(z) \rightarrow M(A)
: u \mapsto \al_u(a)$ is strictly continuous on $S(z)$ and analytic on $S(z)^0$.
\item Consider an element $a$ in $M(A)$. We say that $a$ is strictly analytic with
respect to $\al$ if $a$ belongs to $D(\oal_z)$ for every $z \in \C$. If $a$ is strictly
analytic with respect to $\al$, then the function $\C \rightarrow  M(A) : u \mapsto
\al_u(a)$ is analytic.
\end{itemize}
\end{remark}

\medskip

\begin{proposition}
The mapping $\oal_z$ is a strictly closed multiplicative linear operator in $M(A)$.
\end{proposition}

\medskip

\begin{proposition}
Consider $z \in \C$. Then $D(\oal_{\,\overline{z}}) = D(\oal_z)^*$ and $\text{\rm Ran}\,
\oal_{\,\overline{z}} = (\text{\rm Ran}\, \oal_z)^*$. We have for every $a \in D(\al_z)$
that $\al_z(a)^* = \al_{\overline{z}}(a^*)$.
\end{proposition}

This implies that  for all $z \in \R i$ and $a \in D(\oal_z)$ that $\al_z(a)^*$ belongs to
$D(\oal_z)$ and $\al_z(\al_z(a)^*)^* = a$.

\medskip

\begin{result}
Consider $z \in \C$ and $a \in D(\oal_z)$. Then the function $S(z) \rightarrow M(A) : u
\mapsto \al_u(a)$ is bounded.
\end{result}

\medskip

\begin{proposition}
Consider $z \in \C$ and $t \in \R$. Then $\oal_z \oal_t = \oal_t \oal_z = \oal_{z+t}$.

Furthermore the following equalities hold :
\begin{itemize}
\item We have that $\oal_t(D(\oal_z)) = D(\oal_z)$ and $\oal_t(\text{\rm Ran}\, \oal_z)
= \text{\rm Ran}\, \oal_z$.
\item We have that $D(\oal_{z+t}) = D(\oal_z)$ and $\text{\rm Ran}\, \oal_{z+t} =$
$\text{\rm Ran}\, \oal_z$.
\end{itemize}
\end{proposition}

\medskip

\begin{proposition}
Consider $z \in \C$. Then $\oal_z$ is injective and $(\oal_z)^{-1} = \oal_{-z}$. Therefore
$\text{\rm Ran}\, \oal_z = D(\oal_{-z})$ and $D(\oal_z) = \text{\rm Ran}\, \oal_{-z}$.
\end{proposition}

\medskip

\begin{proposition}
Consider $y,z \in \C$. Then
\begin{enumerate}
\item $\oal_y \, \oal_z \subseteq \oal_{y+z}$.
\item If $y$ and $z$ lie at the same side of the real axis, we have that
$\oal_y \, \oal_z = \oal_{y+z}$.
\end{enumerate}
\end{proposition}

\medskip

As before, we smear elements to construct elements which behave well with respect to
$\al$.

\begin{notation}
Consider $a \in M(A)$, $n > 0$ and $z \in \C$. Then we define the element $a(n,z)$ in
$M(A)$ such that
$$a(n,z) \, b  = \frac{n}{\sqrt{\pi}} \int \exp(-n^2 (t-z)^2) \, \al_t(a) \,b \, dt
$$
for every $b \in A$. Then we have that $\|a(n,z)\| \leq \|a\| \, \exp(n^2
(\text{Im}z)^2)$.
\end{notation}

\medskip

For $n > 0$ and $z \in \C$, we will use the notation $a(n) = a(n,0)$, so
$$a(n) \, b = \frac{n}{\sqrt{\pi}} \int \exp(-n^2 t^2) \, \al_t(a) \, b \, dt $$
for every $b \in A$.

\medskip

\begin{proposition}
Let $a \in M(A)$, $n > 0$ and $z \in \C$. Then $a(n,z)$ is strictly analytic  with respect
to $\al$ and $\al_y(a(n,z)) = a(n,z+y)$ for every $y \in \C$.
\end{proposition}

\medskip

\begin{proposition}
Consider   $n > 0$, $y,z \in \C$ and $a \in D(\oal_y)$. Then $\al_y(a(n,z))$  $=
\al_y(a)(n,z)$.
\end{proposition}

\medskip

Another useful property of smearing is contained in the following result.

\begin{proposition} \label{pa2.prop6}
Let $z \in \C$ and $n > 0$. Consider $a \in M(A)$ and a bounded net $(a_i)_{i \in I}$ in
$M(A)$ such that $(a_i)_{i \in I}$ converges strictly to a. Then $(a_i(n,z))_{i \in I}$ is
also bounded and converges strictly to $a(n,z)$.
\end{proposition}

\bigskip

\section{The definition of KMS-weights}

In this section, we introduce the definition of a KMS-weight. We will also introduce the
definition of a so-called regular weight. In a later section, it is shown that a
KMS-weight is automatically regular. We will also prove a result about such a regular
weight which will be very useful in a later section. Also some results about lower
semi-continuous weights will be included. The standard reference for lower semi-continuous
weights is \cite{Comb}.

\medskip

Let us first introduce notations and conventions concerning weights.

Consider a  \cst-algebra $A$  and a densely defined weight $\vfi$
on $A$. We use the following notations :
\begin{itemize}
\item ${\cal M}^+_\vfi = \{\, a \in A^+ \mid \vfi(a) < \infty  \,\} $
\item ${\cal N}_\vfi = \{\, a \in A \mid \vfi(a^*a) < \infty \,\} $
\item ${\cal M}_\vfi = \text{span\ } {\cal M}^+_\vfi
= {\cal N}_\vfi^* {\cal N}_\vfi$ .
\end{itemize}

\medskip

A GNS-construction of $\vfi$ is by definition a triple
$(H,\pi,\la)$ such that
\begin{itemize}
\item $H$ is a Hilbert space
\item $\la$ is a linear map from ${\cal N}_\vfi$ into $H$ such that
\begin{enumerate}
\item  $\la({\cal N}_\vfi)$ is dense in $H$
\item  We have that $\langle \la(a),\la(b)
\rangle = \vfi(b^*a)$ for every $a,b \in {\cal N}_\vfi$.
\end{enumerate}
\item $\pi$ is a representation of $A$ on $H$ such that $\pi(a)\,\la(b)
= \la(a b)$ for every $a \in A$ and $b \in {\cal N}_\vfi$.
\end{itemize}

\bigskip

The following concepts play a central role in the theory of lower
semi-continuous weights :
\begin{itemize}
\item We define the set ${\cal F}_\varphi = \{\, \omega \in A^*_+ \mid
\omega \leq \varphi \,\}$.
\item Put ${\cal G}_\varphi = \{\, \alpha\,\omega \mid \omega \in {\cal F}_\varphi \, , \, \alpha \in  \,\, ]0,1[ \,\,\}$.
Then ${\cal G}_\varphi$ is a directed subset of ${\cal F}_\varphi$.
\end{itemize}

A proof of the last result can be found in proposition \ref{extr1.prop1}.

\bigskip

The major result about lower semi-continuous weights was proven by Combes (proposition 1.7
of \cite{Comb}) :

\begin{theorem} \label{we1.thm1}
Suppose that $\vfi$ is lower semi-continuous, then we have the following results :
\begin{itemize}
\item We have for every $x \in A^+$ that
$\varphi(x) = \sup \, \{\, \omega(x) \mid \omega \in {\cal F}_\varphi \, \}  $
\item Consider $x \in A^+$, then the net
$\bigl(\omega(x)\bigr)_{\omega \in {\cal G}_\varphi}$ converges to $\vfi(x)$
\item Consider $x \in {\cal M}_\varphi$, then  the net
$\bigl(\omega(x)\bigr)_{\omega \in {\cal G}_\varphi}$ converges to $\vfi(x)$
\end{itemize}
\end{theorem}

\medskip

By proposition 2.4 of \cite{Comb}, we have the following result.

\begin{result} \label{we1.res2}
Consider $\om \in \cF_\vfi$. Then there exists a unique element $T_\om \in B(H)$ with $0
\leq T_\om \leq 1$ such that $\langle T_\om \, \la(a) , \la(b) \rangle$ for every
$a,b \in \Nfi$.
\end{result}

The theorem above implies that $(T_\om)_{\om \in \cG_\vfi}$ converges strongly to $1$ if
$\vfi$ is lower semi-continuous.

\medskip

\begin{result}  \label{we1.res3}
Suppose that $\vfi$ is lower semi-continuous. Then we have the following properties.
\begin{enumerate}
\item The mapping $\la$ is closed.
\item The $^*$-representation $\pi$ is non-degenerate and $\pi(a) \, \la(b) =
\la(a b)$ for every $a \in M(A)$ and $b \in {\cal N}_\vfi$.
\end{enumerate}
\end{result}
\begin{demo}
\begin{enumerate}
\item Choose a sequence $(a_n)_{n=1}^\infty \in \Nfi$, $a \in A$ and $v \in H$ such
that $(a_n)_{n=1}^\infty$ converges to $a$ and $(\la(a_n))_{n=1}^\infty$ converges to $v$.

Choose $\om \in \cF_\vfi$. We have for every $n \in \N$ that $\langle T_\om \, \la(a_n) ,
\la(a_n) \rangle = \om(a_n^* a_n)$, which implies that
$\langle T_\om \, v , v \rangle = \om(a^* a)$. So we see that $\om(a^* a) \leq \|v\|^2$.

Hence, theorem \ref{we1.thm1} implies that $a$ belongs to $\Nfi$.

\medskip

Choose $\om \in \cF_\vfi$.

Take $b \in \Nfi$. We have for every $n \in \Nfi$ that $ \langle T_\om \, \la(a_n) ,
\la(b)\rangle = \om(b^* a_n)$, which implies that
$$\langle T_\om \, v , \la(b) \rangle = \om(b^* a ) = \langle T_\om \, \la(a) , \la(b)
\rangle  \ .$$
This implies that $T_\om \, v = T_\om \, \la(a)$.

The lower semi-continuity of $\vfi$ implies that $(T_\om)_{\om \in \cG_\vfi}$ converges to
1. Hence $v = \la(a)$.

\item From \cite{Comb}, we know already that $\pi$ is non-degenerate. Choose $a \in M(A)$
and $b \in \Nfi$.

Take  $e \in A$. Then $e \, a$ belongs to $A$ which implies that $e \, a
\, b$ belongs to $\Nfi$ and $\la(e\,a\,b) = \pi(e\,a) \la(b) = \pi(e)
\pi(a) \la(b)$.

It is clear that $\Nfi$ is a left ideal in $M(A)$ so $a \, b$ belongs to $\Nfi$. This
implies that $e \, a \, b$ belongs to $\Nfi$ and $\la(e\,a\,b) = \pi(e) \la(a b)$.

So we see that $\pi(e) \pi(a) \la(b) = \pi(e) \la(a b)$.

Therefore, the non-degeneracy of $\pi$ implies that $\la(a b) = \pi(a) \la(b)$.
\end{enumerate}
\end{demo}

\medskip

\bigskip\bigskip

If $\vfi$ is lower semi-continuous, the weight $\vfi$ has a natural extension to a weight
$\overline{\vfi}$ on $M(A)$ by putting $$\overline{\vfi}(x)= \sup\,\{\, \om(x) \mid \om
\in {\cal F}_\vfi \,\}$$ for every $x \in M(A)^+$.

\medskip

Then $\overline{\vfi}$ is the unique strictly lower semi-continuous weight on $M(A)$ which
extends $\vfi$. The unicity follows from the fact that any element in $M(A)^+$ can be
strictly approximated by an increasing net in $A^+$.

\medskip

We define ${\overline{\cal M}}_\vfi = {\cal M}_{\overline{\vfi}}$ and ${\overline{\cal
N}}_\vfi = {\cal N}_{\overline{\vfi}}$. For any $x \in {\overline{\cal M}}_\vfi$, we put
$\vfi(x) = \overline{\vfi}(x)$.

\medskip

Theorem \ref{we1.thm1} implies that $\overline{\vfi}$ is an extension of $\vfi$. So we
have that $\Mfi^+ = \overline{\cM}_\vfi^+ \cap A$, $\Nfi = \overline{\cN}_\vfi \cap A$
and $\Mfi \subseteq \overline{\cM}_\vfi \cap A$

\medskip

Then we have that
\begin{itemize}
\item Consider $x \in M(A)^+$, then the net
$\bigl(\omega(x)\bigr)_{\omega \in {\cal G}_\varphi}$ converges to $\vfi(x)$
\item Consider $x \in \overline{\cal M}_\varphi$, then  the net
$\bigl(\omega(x)\bigr)_{\omega \in {\cal G}_\varphi}$ converges to $\vfi(x)$
\end{itemize}

\bigskip

\begin{lemma}
Suppose that $\vfi$ is lower semi-continuous. Then the mapping $\la : \Nfi \mapsto H$ is
strictly closable.
\end{lemma}
\begin{demo}
Choose a net $(a_i)_{i \in I}$ in $\Nfi$ and $v \in H$ such that $(a_i)_{i \in I}$
converges strictly to 0 and such that $(\la(a_i))_{i \in I}$ converges to $v$.

Take $e \in A$. Then we have that $(e a_i)_{i \in I}$ converges to 0. We have also for
every $i \in I$ that $e a_i$ belongs to $\Nfi$ and $\la(e a_i) = \pi(e) \la(a_i)$. Hence,
$(\la(e a_i))_{i \in I}$ converges to $\pi(e) v$. So the closedness of $\la$ implies that
$\pi(e) v = 0$.

Therefore, the non-degeneracy of $\pi$ implies that $v=0$.
\end{demo}

\medskip

\begin{definition}
Suppose that $\vfi$ is lower semi-continuous. Then we define the mapping $\overline{\la}$
from within $M(A)$ into $H$ such that $\Nfi$ is a strict core for $\overline{\la}$ and
$\overline{\la}$ extends $\la$. We put $\la(a) = \overline{\la}(a)$ for every $a \in
D(\overline{\la})$
\end{definition}

\medskip

Using result \ref{we1.res3}.2, it is not very difficult to prove for every $x \in M(A)$
and every $a \in D(\overline{\la})$ that $x a$ belongs to $D(\overline{\la})$ and $\la(x
a) = \pi(a) \la(x)$.

\medskip

\begin{proposition}
Suppose that $\vfi$ is lower semi-continuous. Then
$\bigl(H,\overline{\la},\overline{\pi}\,\bigr)$ is a GNS-construction for
$\overline{\vfi}$.
\end{proposition}
\begin{demo}
\begin{itemize}
\item Choose $a \in D(\overline{\la})$. By definition, there exists a net $(a_i)_{i \in I}$
in $\Nfi$ such that $(a_i)_{i \in I}$ converges strictly to $a$ an $(\la(a_i))_{i \in I}$
converges to $\la(a)$.

Take $\om \in \cF_\vfi$. Then we have for every $i,j \in I$ that $\langle T_\om \la(a_i) ,
\la(a_j) \rangle = \om (a_j^* a_i)$.

This implies that $\langle T_\om \la(a) , \la(a) \rangle = \om(a^* a)$.

By the definition of $\overline{\vfi}$ and because $(T_\om)_{\om \in \cG_\vfi}$ converges
strongly to 1, we get that $a$ belongs to ${\overline{\cN}}_\vfi$ and that $\vfi(a^* a) =
\langle \la(a) , \la(a) \rangle$.
\item Choose $a \in {\overline{\cN}}_\vfi$. Take an approximate unit $(e_i)_{i \in I}$ of $A$.
Then $(e_i \, a)_{i \in I}$ converges strictly to $a$.

We have for every $i \in I$ that $e_i \, a \in \overline{\cN}_\vfi \cap A = \Nfi$.

Because $(a^* e_i \, a)_{i \in I}$ is an increasing net in $\overline{\cM}_\vfi^+$, which
converges strictly to $a^* a$, the strict lower semi-continuity of $\overline{\vfi}$
implies that $\bigl(\vfi(a^* e_i \, a)\bigl)_{i \in I}$ converges to $\vfi(a^* a)$.

We have moreover for every $j,k \in I$ with $j \leq k$ that $0 \leq e_k - e_j \leq 1$, so
$$ \| \la(e_j \, a) - \la(e_j \, a))  \|
= \vfi( a^* \, (e_k -e_j)^2\, a)
\leq \vfi( a^*\,  (e_k - e_j)\, a)  = \vfi(a^* e_k \, a) - \vfi(a^* e_j \, a)$$

So we get that $(\la(e_i \, a))_{i \in I}$ is cauchy and hence convergent in $H$.

Hence, we get by the definition of $\overline{\la}$ that $a$ belongs to
$D(\overline{\la})$.
\end{itemize}
\end{demo}

\medskip

\begin{remark} \rm
By definition, we have that $\Nfi$ is a strict core for $\overline{\la}$. But we have even
more:

Consider $a \in \overline{\cN}_\vfi$. Then there exists a net $(a_k)_{k \in K}$ in $\Nfi$
such that
\begin{itemize}
\item We have for every $k \in K$ that $\|a_k\| \leq \|a\|$ and
$\|\la(a_k)\| \leq \|\la(a)\|$.
\item The net $(a_k)_{k \in K}$ converges strictly to $a$ and the net
$(\la(a_k))_{k \in K}$ converges strongly$^*$ to $\la(a)$.
\end{itemize}
This follows immediately by multiplying $a$ to the left by an approximate unit of $A$.
\end{remark}

\medskip

We will not use this extension to the multiplier in this paper. However, most of the
results in this paper concerning $\vfi$, $\Mfi$ , $\Mfi^+$, $\Nfi$, $\la$  have obviuous
variants concerning $\vfi$, $\overline{\cM}_\vfi^+$, $\overline{\cM}_\vfi$,
$\overline{\cN}_\vfi$, $\overline{\la}$.

They can be proven using the results in this paper (using strict approximation arguments)
or by similar proofs.

\bigskip \bigskip

Now we give the definition of a KMS-weight.

\begin{definition}
Consider a C$^*$-algebra $A$ and a densely defined lower semi-continuous weight $\vfi$ on
$A$ such that there exist a norm-continuous one parameter group $\si$ on $A$ such that
\begin{enumerate}
\item We have that $\vfi \, \si_t = \vfi$ for every $t \in \R$.
\item We have that $\vfi(a^* a) = \vfi(\si_\frac{i}{2}(a)\,\si_\frac{i}{2}(a)^*)$ for every
$a \in D(\si_\frac{i}{2})$.
\end{enumerate}
Then $\vfi$ is called a KMS-weight and $\si$ is called a modular group
for $\vfi$.
\end{definition}

Later we will prove that the second condition on $\si$ can be weakened.

\medskip

This definition of a KMS-weight is different from the usual one (see definition 4.1 of
\cite{Comb2}) but we will prove in  theorem \ref{we4.thm1} that this definition is
equivalent with the usual one.

\bigskip

In \cite{Ver}, Jan Verding introduced the notion of so-called regular weights. We will
formulate this definition and prove a useful property concerning regular weights.

\begin{definition}
Consider a \cst-algebra $A$. Let $\vfi$ be a weight on $A$  with GNS-construction
$(H,\la,\pi)$. We call $\vfi$ regular $\Leftrightarrow$ $\vfi$ is densely defined, lower
semi-continuous and there exists a net $(u_i)_{i \in I}$ in $\Nfi$ such that
\begin{itemize}
\item We have for every $i \in I$ that :
 \begin{enumerate}
 \item There exists $S_i \in B(H)$ such that $S_i \la(a) = \la(a u_i)$ for
  every $a \in \Nfi$.
 \item $\|u_i\| \leq 1$ and $\|S_i\| \leq 1$.
 \end{enumerate}
\item Furthermore, we assume that
\begin{enumerate}
\item $(u_i)_{i \in I}$ converges strictly to 1.
\item $(S_i)_{i \in I}$ converges strongly to 1.
\end{enumerate}
\end{itemize}
We call the net $(u_i)_{i \in I}$ a truncating net for $\vfi$.
\end{definition}

\medskip

It is not difficult to see that the definition of regularity is independent of the choice
of  GNS-construction.

\bigskip

\begin{proposition} \label{we1.prop2}
Consider a regular weight $\vfi$ on a \cst-algebra $A$ with GNS-construction
$(H,\la,\pi)$.  Let $B$ be a sub-$^*$-algebra of $\Nfi \cap \Nfi^*$ such that $B$ is a
core for $\la$ and  let $a$ be an element in $\Nfi \cap \Nfi^*$. Then there exists a
sequence $(a_n)_{n=1}^\infty$ in $B$ such that
\begin{enumerate}
\item $(a_n)_{n=1}^\infty \rightarrow a$
\item $\bigl(\la(a_n)\bigr)_{n=1}^\infty \rightarrow \la(a)$
\item $\bigl(\la(a_n^*)\bigr)_{n=1}^\infty \rightarrow \la(a^*)$.
\end{enumerate}
\end{proposition}
\begin{demo}
Choose a truncating net $(u_i)_{i \in I}$ for $\vfi$. Define for every $i \in I$ the
operator $S_i \in B(H)$ such that $S_i \la(a) = \la(a u_i)$ for every $a \in \Nfi$.

\medskip

Choose $n \in \N$. Because $(u_i)_{i \in I}$ is a bounded net which converges strictly to
1 and $(S_i)_{i \in I}$ is a bounded net which converges strongly to 1, there exist an
element $j \in I$ such that
\begin{itemize}
\item
$\| \la(a) - \pi(u_j^*) S_j \la(a) \| \leq \frac{1}{n}$
\item
$\| \la(a^*) - \pi(u_j^*) S_j \la(a^*) \| \leq \frac{1}{n}$
\item
$\| a - u_j^*  a u_j \| \leq \frac{1}{n} .$
\end{itemize}
By the definition of $S_j$, we get that
\begin{itemize}
\item
$\| \la(a) - \pi(u_j^*) \pi(a) \la(u_j) \| \leq \frac{1}{n}$
\item
$\| \la(a^*) - \pi(u_j^*) \pi(a^*) \la(u_j) \| \leq \frac{1}{n} .$
\end{itemize}

Choose $n \in \N$. Because $B$ is a core for $\la$, there exists an element $b_n \in B$
such that
\begin{itemize}
\item
$\|u_j -  b_n \| \leq  \min\{ \frac{1}{2n} \frac{1}{\|a\|+1} \, , \, \frac{1}{n}
\frac{1}{(\|a\|+1) \, (\|\la(u_j)\|+1)} \} $
\item
$ \| \la(u_j) - \la(b_n) \| \leq \frac{1}{2 n} \frac{1}{\|a\|+1} $
\item
$\| b_n \| \leq 2$
\end{itemize}
Furthermore, there exists an element $d_n \in B$ such that
$$\|a - d_n\| \leq  \min \{ \frac{1}{n} \frac{1}{\|b_n\|^2 + 1} \, , \,
\frac{1}{n} \frac{1}{(\|b_n\|+1)(\|\la(b_n)\|+1)} \}$$

We define $a_n = b_n^* d_n b_n$. Then $a_n$ belongs to $B$
and $a_n^* = b_n^* d_n^* b_n$.

Moreover,

\begin{enumerate}

\item Using the above estimations, we have that
\begin{eqnarray*}
& & \| \la(a_n) - \la(a) \| = \| \pi(b_n^*) \pi(d_n) \la(b_n) - \la(a) \| \\
& & \spat \leq \| \pi(b_n^*) \pi(d_n) \la(b_n) - \pi(b_n^*) \pi(a) \la(b_n) \|  +  \|
\pi(b_n^*) \pi(a) \la(b_n) - \pi(b_n^*) \pi(a) \la(u_j) \| \\
& & \spat\ \ \ + \| \pi(b_n^*) \pi(a) \la(u_j) - \pi(u_j^*) \pi(a) \la(u_j) \| +
 \| \pi(u_j^*) \pi(a) \la(u_j) - \la(a) \| \\
& & \spat \leq   \| b_n \| \, \| \la(b_n) \| \, \| d_n - a \| +
 \| b_n \| \, \| a \| \, \| \la(b_n) - \la(u_j) \|
+ \| b_n - u_j \| \, \|a\| \, \| \la(u_j) \| + \frac{1}{n} \\
& & \spat \leq  \frac{1}{n} + 2 \,  \| a \| \, \| \la(b_n) - \la(u_j) \|
 + \frac{1}{n} + \frac{1}{n}
\leq  \frac{1}{n} + \frac{1}{n} + \frac{1}{n} + \frac{1}{n} = \frac{4}{n}
\end{eqnarray*}

\item Similarly, one can prove
$\|\la(a_n^*) - \la(a^*) \| \leq \frac{4}{n}$.

\item  Using the above estimations once again, we get that
\begin{eqnarray*}
& & \|a_n - a\|  =  \|b_n^* d_n b_n - a \| \\
& & \spat \leq  \| b_n^* d_n b_n - b_n^* a b_n \| + \| b_n^* a b_n - b_n^* a u_j \|
+ \| b_n^* a u_j - u_j^* a u_j \| + \| u_j^*  a u_j - a \| \\
& & \spat \leq  \|b_n\|^2 \|d_n - a \| + \|b_n\| \, \|a\| \, \| b_n - u_j \|
+ \| b_n - u_j\| \, \| a \| \, \|u_j\| + \| u_j^* a u_j - a \| \\
& & \spat \leq  \frac{1}{n} + 2 \, \|a\| \, \| b_n -  u_j \|
+ \| b_n - u_j \| \, \| a \| + \frac{1}{n}
\leq \frac{1}{n} + \frac{1}{n} + \frac{1}{n} + \frac{1}{n} =
\frac{4}{n} .
\end{eqnarray*}

\end{enumerate}

Looking at 1, 2 and 3, we see that $(a_n)_{n=1}^\infty \rightarrow a$,
$\bigl(\la(a_n)\bigr)_{n=1}^\infty \rightarrow \la(a)$ and
$\bigl(\la(a_n^*)\bigr)_{n=1}^\infty
\rightarrow \la(a^*)$.
\end{demo}

\medskip

\begin{corollary} \label{we1.cor1}
Consider a regular weight $\vfi$ on a \cst-algebra $A$ with GNS-construction
$(H,\la,\pi)$. Let $a$ be an element in $\Nfi \cap \Nfi^*$. Then there exists a sequence
$(a_n)_{n=1}^\infty$ in $\Mfi$ such that
\begin{enumerate}
\item $(a_n)_{n=1}^\infty \rightarrow a$
\item $\bigl(\la(a_n)\bigr)_{n=1}^\infty \rightarrow \la(a)$
\item $\bigl(\la(a_n^*)\bigr)_{n=1}^\infty \rightarrow \la(a^*)$.
\end{enumerate}
\end{corollary}

\bigskip

We need the following result from \cite{Pe-Tak} (lemma 3.1). The proof remains valid for
lower semi-continuous weights on \cst-algebras.

\begin{lemma}
Consider a \cst-algebra $A$ and a densely defined lower semi-continuous weight $\vfi$ on
$A$. Let $f$ be a function from $\R$ into $\Mfi^+$ which is continuous and integrable.

Then the function $\R \rightarrow \R^+  : t \mapsto \vfi(f(t))$ is lower semi-continuous
and
$$\vfi\bigl(\,\int f(t) \, dt\,\bigr) = \int \vfi(f(t)) \, dt \ .$$
\end{lemma}

\bigskip

\begin{lemma}
Consider a \cst-algebra $A$ and a densely defined lower semi-continuous weight $\vfi$
on $A$, $\al$ a norm continuous one-parameter group on $A$ such that there exists
a strictly positive number $\lambda$ such that $\vfi \, \al_t = \lambda^t \, \vfi$ for every $t
\in \R$. Let $z \in \C$, $n \in \N$ and $a \in \Mfi$.

Then we have that the element $\int \exp(-(t-z)^2) \, \al_t(a) \, dt$ belongs to
 $\Mfi$ and $$\vfi\bigl( \,\int \exp(-(t-z)^2) \, \al_t(a) \, dt \, \bigr)
= \lambda^z \, \sqrt{\pi} \, \exp(\frac{(\ln \lambda)^2}{4}) \, \vfi(a) \ .$$
\end{lemma}
\begin{demo}
Choose $b \in \Mfi^+$. Define the function $f:\R \rightarrow \C:
t \mapsto \exp(-(t-z)^2)$.

There exist continuous positive functions
$f_1,\ldots\!,f_4$ such that $f_1,\ldots\!,f_4 \leq |f|$ and
$f = f_1 - f_2 + i f_3 - i f_4$.

\medskip

Fix $j \in \{1,\ldots\!,4\}$. It is clear that the function $\R \rightarrow \Mfi^+: t
\mapsto f_j(t) \al_t(b)$ is continuous. Because we also have that $\|f_j(t) \al_t(b)\| =
f_j(t) \|b\| \leq |f(t)| \, \|b\|$ for every $t \in \R$, we see that the
function $\R \rightarrow \Mfi^+: t
\mapsto f_j(t) \al_t(b)$ is integrable. So, the previous lemma implies that the function
$\R \rightarrow \R^+: t \mapsto  \vfi(f_j(t) \al_t(b))$ is lower semi-continuous and
$$\vfi\bigl(\, \int f_j(t) \al_t(b) \, dt\, \bigr) = \int \vfi( f_j(t) \al_t(b)) \, dt
= \bigl(\,\int f_j(t) \lambda^t \, dt\,\bigr) \, \vfi(b) \ .$$
Because $f_j \leq |f|$, the function $\R \rightarrow \R^+ : t \mapsto f_j(t) \lambda^t$ is
integrable. Therefore we see that $\int f_j(t) \al_t(b) \, dt$ belongs to $\Mfi^+$.

\medskip

Using the fact that $f = f_1 - f_2 + i f_3 - i f_4$, we get that
$\int f(t) \al_t(b) \, dt $ belongs to $\Mfi$ and
$$ \vfi\bigl( \, \int f(t) \al_t(b) \, dt \, \bigr) =
\bigl(\,\int f(t) \lambda^t \, dt\,\bigr) \, \vfi(b)  \ .$$ Therefore,
\begin{eqnarray*}
& & \vfi\bigl( \, \int f(t) \al_t(b) \, dt \, \bigr)
= \bigl(\,\int \exp(-(t-z)^2) \, \lambda^t \, dt\,\bigr) \,
\vfi(b) \\ & & \spat =  \bigl(\,\int \exp(-t^2) \, \lambda^{t+z} \, dt\,\bigr) \, \vfi(b)
= \lambda^z \, \sqrt{\pi} \, \exp(\frac{(\ln \lambda)^2}{4}) \, \vfi(b) \ .
\end{eqnarray*}
The lemma follows by linearity.
\end{demo}

\begin{proposition} \label{we1.prop3}
Consider a \cst-algebra $A$ and a densely defined lower semi-continuous weight $\vfi$ on
$A$, $\al$ a norm continuous one-parameter group on $A$ such that there  exists a
strictly positive number $\lambda$ such that $\vfi \, \al_t = \lambda^t \, \vfi$ for every $t
\in \R$. Let $z \in \C$ and $a \in \Mfi \cap D(\al_z)$ such that $\al_z(a)$
belongs to $\Mfi$. Then we have that $\vfi(\al_z(a)) = \lambda^z \, \vfi(a)$.
\end{proposition}
\begin{demo}
Define
$$ x = \int \exp(-t^2) \, \al_t(\al_z(a)) \, dt \ .$$
Because $\al_z(a)$ belongs to $\Mfi$, the previous lemma implies that
$x$ belongs to $\Mfi$  and
$$\vfi(x) = \sqrt{\pi} \, \exp(\frac{(\ln \lambda)^2}{4}) \, \vfi(\al_z(a)) \ .$$
We have also that $$x = \int \exp(-(t-z)^2) \, \al_t(a) \, dt \ ,$$ so the previous lemma
implies in this case that
$$\vfi(x) = \lambda^z \, \sqrt{\pi} \, \exp(\frac{(\ln \lambda)^2}{4}) \, \vfi(a) \ .$$
Comparing these two different expressions for $\vfi(x)$ gives us that
$\vfi(\al_z(a)) = \lambda^z \, \vfi(a)$.
\end{demo}

\bigskip

\begin{result}
Consider a densely defined lower semi-continuous non-zero weight $\vfi$ on a \cst-algebra
$A$. Let $\al$ be a norm continuous one parameter group on $A$ such that there exists for
every $t$ in $\R$ a number $\lambda_t$ such that $\vfi \, \al_t = \lambda_t \, \vfi$ for
every $t \in \R$. Then there exists a unique strictly positive number $\lambda$ such that
$\vfi \, \al_t = \lambda^t \, \vfi$ for every $t \in \R$.
\end{result}
\begin{demo} The technique of this proof comes from proposition 5.7 of \cite{Pe-Tak} but can be more
easily applied in this case. We have certainly that $\lambda_t > 0$ for every $t \in \R$.
Put $\lambda = \lambda_1 > 0$.

Because $\al_{s+t} = \al_s \al_t$ for every $s,t \in \R$, we get easily that
$\lambda_{s+t} = \lambda_s \lambda_t$ for every $s,t \in \R$. As usual, this implies that
$\lambda_q = \lambda^q$ for every $q \in Q$.

The lower semi-continuity of $\vfi$ implies that the mapping $\R \rightarrow \R^+_0 : t
\mapsto \lambda_t$ is lower semi continuous. So the mapping $\R \rightarrow \R : t
\mapsto \lambda_t - \lambda^t$ is lower semi-continuous. This implies that the set $\{\, t
\, \in \R \mid \lambda_t \leq \lambda^t \,\}$ is closed in $\R$. We know already that this
set contains $\Q$, so it must be equal to $\R$.

Therefore $\lambda_t \leq \lambda^t$ for every $t \in \R$. This implies that $\lambda_t =
\lambda^t$ for every $t \in \R$ (if there would exist an element in $\R$ for which this
equality does not hold, this element or its opposite would violate the previous
inequality).
\end{demo}

\medskip

The proof of the following result is due to J. Verding (see \cite{Verd}).

\begin{result} \label{we1.res1}
Consider a \cst-algebra $A$ and a densely defined lower semi-continuous weight $\vfi$ on
$A$. Take a GNS-construction $(H,\la,\pi)$ for $\vfi$. Let $\al$ be a norm continuous one
parameter group on $A$ such that there exists a strictly positive number $\lambda$ such
that $\vfi \, \al_t = \lambda^t \, \vfi$ for every $t \in \R$. Then there exists a unique
injective positive operator $T$ in $H$  such that $T^{it} \la(a) = \lambda^{-\frac{t}{2}}
\, \la(\al_t(a))$ for every $a \in \Nfi$ and $t \in \R$.
\end{result}
\begin{demo}
For every $t \in \R$, there exist a unitary operator $u_t$ on $H$ such that $u_t \la(a) =
\lambda^{-\frac{t}{2}} \, \la(\al_t(a))$ for every $a \in \Nfi$. The mapping $\R
\rightarrow B(H) : t \mapsto u_t$ is clearly a unitary representation of $\R$ on $H$

\medskip

Fix $\om \in \cG_\vfi$ for the moment. Because $\om \leq \vfi$, there exists a unique
element $T_\om \in B(H)^+$ such that $\langle T_\om \la(a) , \la(b) \rangle = \om(b^* a)$
for every $a,b \in \Nfi$ (see result \ref{we1.res2}).

Because $\vfi$ is lower semi-continuous,  the remarks after result \ref{we1.res2} imply that
$(T_\om)_{\om \in \cG_\vfi}$ converges strongly to 1. This implies that the set $\{\,
T_\om \la(a) \mid \om \in \cG_\vfi, a \in \Nfi \,\}$ is dense in $H$. \ \ \ \ \ \ (*)

\medskip

Choose $\th \in \cG_\vfi$ and $x,y \in \Nfi$. We have for every $t \in \R$ that
$$\langle u_t \la(x) , T_\th \la(y) \rangle
= \lambda^{-\frac{t}{2}} \, \langle \la(\al_t(x)) , T_\th \la(y) \rangle
= \lambda^{-\frac{t}{2}} \, \th(y^* \al_t(x)) \ . $$
This implies that the mapping $\R \rightarrow \C : t \mapsto
\langle u_t \la(x) , T_\th \la(y) \rangle$ is continuous.

Therefore (and  because $\|u_t\| \leq 1$ for every $t \in \R$), result (*) implies that
the mapping $\R \rightarrow \C : t \mapsto \langle u_t v , w \rangle$ is continuous for
all $v,w \in H$.

\medskip

Consequently, the mapping $\R \rightarrow B(H) : t  \mapsto u_t$ is a strongly continuous
unitary group representation of $\R$ on $H$. The result follows by the Stone theorem.
\end{demo}

\bigskip

\begin{terminology}
Consider a \cst-algebra $A$ and an element $T$ affiliated with $A$. A truncating sequence
for $T$ is by definition a sequence $(e_n)_{n=1}^\infty$ in $M(A)$ such that :
\begin{enumerate}
\item $(e_n)_{n=1}^\infty$ is bounded and converges strictly to 1.
\item We have for every $n \in N$ that $e_n$ is a left and right multiplier of $|T|$ and
$|T| \, e_n = e_n \, |T|$ .
\end{enumerate}
\end{terminology}

It is not difficult to check that the functional calculus for $|T|$ guarantees the
existence of such a truncating sequence for $T$

\bigskip

Let us fix for the moment an element $T$ affiliated with a \cst-algebra $A$ and a
truncating sequence $(e_n)_{n=1}^\infty$ for $T$. Then we have the following properties.
\begin{itemize}
\item  Take $m \in \N$. We know for every $x \in A$ that $e_m x$ belongs
to $D(|T|)$, so $e_m x$ belongs to $D(T)$. Furthermore
$$ \|T (e_m x)\| = \|\,|T| (e_m x) \| \ .$$
This implies that $T e_m$ is a bounded linear operator on $A$ with $\| T e_m \| = \| \,
|T| e_m \|$. So $e_m$ is a right multiplier of $T$.

\item Let $x \in D(T)$. Then we have for every $n \in N$ that
$$\|T(e_n x) - T(x)\| = \| \,|T|(e_n x) - |T|(x) \|
= \|\, e_n \,|T|(x) - |T|(x)\| \ .$$
This implies that $\bigl(T(e_n x)\bigr)_{n=1}^\infty$ converges to $T(x)$.

Combining this with the closedness of $T$, this implies for every $x \in A$ that $x$ belongs to $D(T)$
if and only if the sequence $\bigl(T(e_n x)\bigr)_{n=1}^\infty$ is convergent in $A$.
\end{itemize}

\medskip

\begin{proposition} \label{we1.prop1}
Consider a \cst-algebra $A$ and an element $T$ affiliated with $A$. Let $\vfi$ be a densely
defined lower semi-continuous weight on $A$  with GNS-construction $(H,\la,\pi)$.
Let $a$ be an element in $D(T) \cap \Nfi$,
then :
\begin{enumerate}
\item We have that $T(a)$ belongs to $\Nfi$ $\Leftrightarrow$
      $\la(a)$ belongs to $D(\pi(T))$.
\item If $T(a)$ belongs to $\Nfi$, then $\la(T(a)) = \pi(T) \la(a)$.
\end{enumerate}
\end{proposition}
\begin{demo}
Take a truncating sequence $(e_n)_{n=1}^\infty$ for $T$. It is not too difficult to check
that $\pi(T) \pi(e_n) = \pi(T e_n) \in B(H)$ and
$\pi(e_n) |\pi(T)| \subseteq |\pi(T)| \pi(e_n) = \pi(|T| e_n) \in B(H)$.
\begin{itemize}
\item
Suppose that $T(a)$ belongs to $\Nfi$. Because $a$ belongs to $D(|T|)$ and $(|T|\,a)^*
(|T|\,a) = (Ta)^*(Ta)$, we see that $|T|\,a$ also belongs to $\Nfi$.

Fix $n \in \N$.  We know  that $|T| e_n$ belongs to $M(A)$ and $e_n |T| \subseteq |T|
e_n$. We know also that $\pi(|T|) \pi(e_n)$ belongs to $B(H)$ and $|\pi(T)| \pi(e_n) =
\pi(|T| e_n)$.

Hence, $\pi(e_n) \la(a)$ belongs to $D(|\pi(T)|)$ and
\begin{eqnarray*}
|\pi(T)|(\pi(e_n) \la(a)) & = & \pi(|T| e_n) \la(a)
= \la((|T| e_n) a) \\
& = & \la(e_n (|T|a)) = \pi(e_n) \la(|T| a) \ .
\end{eqnarray*}

This last equality implies that $\bigl(\,|\pi(T)|(\pi(e_n) \la(a))\,\bigr)_{n=1}^\infty$
converges to    $\la(|T| a)$. Therefore $\la(a)$ must belong to $D(|\pi(T)|)$, so it must
belong to $D(\pi(T))$.

\item
Suppose that $\la(a)$ belongs to $D(\pi(T))$. Then $\la(a)$ will also belong to $D(|\pi(T)|)$.

For every $n \in N$, we have that $T e_n$ belongs to $M(A)$, which implies that $(T e_n)a$
belongs to $\Nfi$ and
$$\la((T e_n) a)= \pi(T e_n) \la(a)
= \pi(T) (\pi(e_n) \la(a)) \ .$$

Because $\la(a)$ belongs to $D(|\pi(T)|)$,
this equality implies for every $n \in \N$ that
\begin{eqnarray*}
& & \| \la( (T e_n) a ) - \pi(T) \la(a) \|
= \| \pi(T)( \pi(e_n) \la(a) - \la(a) ) \|  \\
& & \!\!\spat =  \| \, |\pi(T)|( \pi(e_n) \la(a) - \la(a))\|
=  \| \, \pi(e_n) \, |\pi(T)| \la(a) - |\pi(T)| \la(a) \|
\end{eqnarray*}
So $\bigl(\la((T e_n)a)\bigr)_{n=1}^\infty$ converges to $\pi(T) \la(a)$.

We also know that $\bigl(T (e_n a) \bigr)_{n=1}^\infty$ converges to $T(a)$, so the
closedness of $\la$ implies that $T(a)$ belongs to $\Nfi$ and $\la(T(a))= \pi(T) \la(a)$.
\end{itemize}
\end{demo}

\bigskip

\section{Weights arising from GNS-constructions}

In this section, we gather some results from \cite{Verd} and \cite{Q-V}. Except lemma
\ref{extr1.lem1}, all the results and proofs in this section are due to Jan Verding and
can be found in \cite{Verd} and \cite{Q-V}. However, due to the lack of availability of
these two works, we feel it necessary to include the proofs here.

\bigskip

Consider a \cst-algebra $A$,  a Hilbert space $H$, a dense left ideal $N$ of $A$ and a
linear map $\la$ from $N$ into $H$ with dense range. Assume furthermore the existence of a
$^*$-homomorphism $\pi$ from $A$ into $B(H)$  such that $\pi(x) \la(a) = \la(x a )$ for
every $x \in A$ and $a \in N$.

\medskip

\begin{notation}  \label{extr1.not1}
We define the sets $$\cF = \{\, \om \in A^+ \mid \text{ We have that \  } \om(a^* a)
\leq \|\la(a)\|^2 \text{ for every } a \in N \, \} $$ and
$$\cG = \{ \, \lambda \, \om \mid \om \in \cF , \lambda \in ]0,1[ \,\, \} \ .$$
\end{notation}

\begin{remark} \rm
Consider $\om \in \cF$. Then it is not difficult to see that there exists a unique element
$T \in B(H)$ with $0 \leq T \leq 1$ such that $\langle T \la(a) , \la(b) \rangle = \om(b^*
a)$ for every $a,b \in N$. We get also that $T$ belongs to $\pi(A)'$. The techniques used
to prove these results can also be found in the proof of the next lemma.
\end{remark}

\bigskip

\begin{lemma} \label{extr1.lem1}
Consider a positive sesquilinear form on $N$ such that $s(a b_1,b_2)
= s(b_1,a^* b_2)$ for all $a\in A$ and all $b_1,b_2\in N$.  Suppose moreover that
there exists a positive linear functional $\th$ on $A$ such $s(b,b)\leq \th(b^*b)$ for
every $b \in N$. Then there exists a unique positive linear functional $\omega$ on $A$
with $\om \leq \th$ such that $\omega(b_2^* \, b_1) = s(b_1,b_2)$ for every $b_1,b_2 \in
N$.
\end{lemma}

\begin{demo}
Let $(\pi,H,v)$ be GNS-object for $\th$ ($v$ is a cyclic vector).

Because $s$ is a positive sesquilinear form on N, we can use the Cauchy-Schwarz inequality
for $s$. So we have for every $b_1,b_2\in N$ that

$$ |s(b_1,b_2)|^2 \leq  s(b_1,b_1)\,\,s(b_2,b_2)
\leq  \th(b_1^*\,b_1)\,\,\th(b_2^*\,b_2)
=  \|\pi(b_1) v\|^2\,\, \|\pi(b_2) v\|^2  .$$
Therefore we can define a continuous positive sesquilinear form $t$
on $H$ such that
$t(\pi(b_1) v,\pi(b_2)v) = s(b_1,b_2)$ for every $b_1,b_2 \in N$.
It is clear that $t$ is positive and $\|t\| \leq 1$.

So there exist an element $T \in B(H)$ with $0\leq T \leq 1$ such that
$t(x,y) = \langle T x ,y \rangle$ for every $x,y \in H$.

This implies that $\langle T \pi(b_1) v , \pi(b_2) v \rangle = s(b_1,b_2)$ for every
$b_1,b_2 \in N$.

\medskip

Next we show that $T$ belongs to $\pi(A)'$. Therefore, choose $a \in N$.

We have for every $b_1,b_2 \in N$ that
\begin{eqnarray*}
\langle T \pi(a)\, \pi(b_1)v ,\pi(b_2)v \rangle & = &
\langle T \pi(a b_1)v ,\pi(b_2)v \rangle = s(a b_1,b_2) \\
& = & s(b_1,a^* b_2) =  \langle T \pi(b_1) v,\pi(a^* b_2)v \rangle \\
& = & \langle T \pi(b_1)v , \pi(a^*)\pi(b_2)v \rangle
=  \langle \pi(a) T \pi(b_1) v , \pi(b_2)v \rangle .
\end{eqnarray*}
This implies that $T \pi(a) = \pi(a) T$.

\medskip

Now we define the continuous linear functional $\omega$ on $A$ such
$\omega(x) = \langle T \pi(x)v , v \rangle $ for every $x \in A$.

Using the fact that $T$ belongs to $\pi(A)'$, we have for every $b_1,b_2 \in N$ that
$$\om(b_2^*\,b_1) = \langle T \pi(b_2^* \, b_1) v , v \rangle
= \langle T \pi(b_1) v , \pi(b_2) v \rangle = s(b_1,b_2). $$

Consequently, we have for every  $b \in N$ that $\om(b^* \, b) = s(b,b)$,
implying that $0 \leq \om(b^* b) \leq \th(b^* b)$.
This implies easily that $0 \leq \om \leq \th$.
\end{demo}

\medskip

\begin{lemma}
Consider a unital \cst-algebra $C$. Let $T_1,T_2$ be elements in $C$ with
$0 \leq T_1,T_2 \leq 1$, let $\gamma$ be a number in $]0,1[$.
Then there exist an element $T \in C$ with $0 \leq T \leq 1$ and such
that $\gamma \, T_1 \leq T$,  $\gamma \, T_2 \leq T$
and $T \leq \frac{\gamma}{1-\gamma} (T_1+T_2)$.
\end{lemma}
\begin{demo}
For the moment, fix $i \in \{1,2\}$. We  define
$$S_i = \frac{\gamma \, T_i}{1-\gamma \, T_i} \in C.$$

It is then easy to check that
$$ \gamma \, T_i = \frac{S_i}{1+S_i}   \text{\ \ \ \ \ \ \ (a)} \hspace{1.5cm}
\text{ and } \hspace{1.5cm} S_i \leq \frac{\gamma}{1-\gamma}\,
T_i \text{\ \ \ \ \ \ \ (b)} \ .$$

Next, we define
$$T = \frac{S_1 + S_2}{1+S_1+S_2} \in C .$$
We get immediately that $0 \leq T \leq 1$.
By (b), we have that $T \leq S_1+S_2 \leq \frac{\gamma}{1-\gamma} (T_1 + T_2)$.

We know that the function $\R^+ \rightarrow \R^+ : t \mapsto \frac{t}{1+t}$ is operator
monotone (see \cite{Ped}). This implies for every $i=1,2$ that
$$ T \geq \frac{S_i}{S_i+1} = \gamma \, T_i \ ,$$
where we used (a) in the last equality.
\end{demo}

\medskip

We now use these two lemmas to proof the following useful proposition.

\begin{proposition} \label{extr1.prop1}
The set ${\cal G}$ is upwardly directed for the natural order on $A^*_+$.
\end{proposition}
\begin{demo}
Choose $\om_1, \om_2 \in {\cal F}$ and $\lambda_1,\lambda_2 \in ]0,1[$. Then there exist a
number $\gamma \in ]0,1[$ such that $\lambda_1,\lambda_2 < \gamma$.

By the remark after notation \ref{extr1.not1},  we know that there exist $T_1,T_2$ in
$\pi(A)'$ such that $0  \leq T_1,T_2 \leq 1$ and such that  $\om_1(y^* x) = \langle T_1
\la(x) , \la(y) \rangle$ and $\om_2(y^* x) = \langle T_2 \la(x) ,
\la(y) \rangle$ for every $x,y \in N$.

\medskip

The previous lemma implies the existence of $T
\in \pi(A)'$ with $0 \leq T \leq 1$ and such that $\gamma \, T_1 \leq T$,  $\gamma
\, T_2 \leq T$ and $T \leq \frac{\gamma}{1-\gamma} (T_1+T_2)$.

\medskip

Put $\lambda = \max(\frac{\lambda_1}{\gamma},\frac{\lambda_2}{\gamma}) \in ]0,1[$. \ Then
$\lambda T \geq \lambda \gamma \, T_1 \geq \lambda_1 T_1 $ and analogously, $\lambda T
\geq \lambda_2 T_2$. \ \ \ \ \ (a)

\medskip

Define $\th  = \frac{\gamma}{1-\gamma}(\om_1 + \om_2) \in A^*_+$

Next, we define the mapping $s$ from $N \times N$ into $\C$ such that $s(x,y) =
\langle T \la(x) , \la(y) \rangle$ for every $x,y \in N$. Then $s$ is sesquilinear.
Furthermore:
\begin{itemize}
\item We have for $x \in N$ that
$$ 0 \leq \langle T \la(x) , \la(y) \rangle \leq \frac{\gamma}{1-\gamma} \, \langle
(T_1 + T_2) \la(x) , \la(y) \rangle  = \frac{\gamma}{1-\gamma} \, (\om_1(x^* x) +
\om_2(x^* x)) , $$ which implies that $0 \leq s(x,x) \leq \th(x^* x)$.
\item We have for every $x,y \in N$ and $a \in A$ that
\begin{eqnarray*}
\langle T \la(a x) , \la(y) \rangle  & = & \langle  T \pi(a) \la(x) , \la(y) \rangle
= \langle \pi(a) T \la(x) , \la(y) \rangle \\
& = & \langle T \la(x), \pi(a^*) \la(y) \rangle  = \langle T \la(x) , \la(a^* y) \rangle
\end{eqnarray*}
which implies that $s(a x,y) = s(x, a^* y)$.
\end{itemize}

This allows us to apply the previous lemma. Therefore, we get the existence of $\om
\in  A_+^*$ with $\om \leq \th$ and such that $s(x,y) = \om(y^* x)$ for every $x,y \in
N$.

\medskip

Then we have for every $x \in N$ that $$\om(x^* x) = s(x,x) = \langle T \la(x) , \la(x)
\rangle \leq \|\la(x)\|^2$$ which implies that $\om$ belongs to $\cF$.

By (a), we have moreover for every $x \in N$ that
$$\lambda \, \om(x^* x) = \lambda \, \langle T \la(x) , \la(x) \rangle \geq \lambda_1 \,
\langle T_1 \la(x) , \la(x) \rangle = \lambda_1 \, \om_1(x^* x) \ $$
so we see that $\lambda_1 \, \om_1 \leq \lambda \, \om$. We get in a similar way that
$\lambda_2 \, \om_2 \leq \lambda \, \om$
\end{demo}

\medskip

This allows us to give the following definition :

\begin{definition} \label{extr1.def1}
We define the weight $\vfi$ on $A$ such that $\vfi(x) = \sup \, \{ \, \om(x) \mid \om \in
\cF \, \}$ for every $x \in A^+$.

Then $\vfi$ is a densely defined lower semi-continuous weight on $A$ such that $N
\subseteq \Nfi$ and $\vfi(x^* x) \leq \|\la(x)\|^2$ for every $x \in N$.
\end{definition}

\bigskip

In the last proposition of this section, we will improve the relationship between $\vfi$
and the mapping $\la$ \ (under some extra conditions).

\medskip

\begin{lemma}
Let $E$ be a normed space, $H$ a Hilbert space and $\la$ linear mapping from within $E$
into $H$. Let $(x_i)_{i \in I}$ be a net in $D(\la)$ and $x$ an element in $E$ such that
$(x_i)_{i\in I}$ converges to $x$ and $(\la(x_i))_{i\in I}$ is bounded. Then there exists
a sequence $(y_n)_{n=1}^\infty$ in the convex hull of $\{\, x_i \mid i \in I \,\}$ and an
element $v \in H$ such that $(y_n)_{n=1}^\infty$ converges to $x$ and
$(\la(y_n))_{n=1}^\infty$ converges to $v$.
\end{lemma}
\begin{demo}
By the Banach-Alaoglu theorem, there exists a subnet $(x_{i_j})_{j\in J}$ of $(x_i)_{i\in
I}$ and $v \in H$ such that $(\la(x_{i_j}))_{j \in J}$ converges to $v$ in the weak
topology on $H$. (For this, we need $H$ to be a Hilbert space.)

\medskip

Fix $n \in \N$. Then there exists $j_n \in J$ such that $\|x_{i_j} - x\| \leq \frac{1}{n}$
for all $j \in J$ with $j \geq j_n$.

Now $v$ belongs to the weak-closed convex hull of the set $\{\,\la(x_{i_j}) \mid j \in J
\text{ such that } j \geq j_n \,\}$, which is the same as  the norm-closed convex hull.

Therefore, there exist $\lambda_1,\ldots\!,\lambda_m \in \R^+$ with $\sum_{k=1}^m
\lambda_k = 1$ and elements $\alpha_1,\ldots\!,\alpha_m \in J$  with
$\alpha_1,\ldots\!,\alpha_m \geq j_n$  such that
$$ \| v - \sum_{k=1}^m \lambda_k \la(x_{i_{\alpha_k}}) \| \leq \frac{1}{n} .$$

Put $y_n =  \sum_{k=1}^{m} \lambda_k x_{i_{\alpha_k}}$. Then
$y_n \in D(\la)$, and $\la(y_n) = \sum_{k=1}^m \lambda_k \la(x_{i_{\alpha_k}})$.

Therefore, we have immediately that $\|v - \la(y_n)\| \leq \frac{1}{n}$.

Furthermore,
$$\|x - y_n\| = \left\|\sum_{k=1}^{m}\lambda_k(x - x_{i_{\alpha_k}})\right\| \leq
\sum_{k=1}^{m} \lambda_k  \frac{1}{n} = \frac{1}{n}$$

\medskip

Therefore, we find that $(y_n)_{n=1}^\infty$ converges to $x$ and that
$(\la(y_n))_{n=1}^\infty$ converges to $v$.
\end{demo}

\medskip

Using this lemma, we can prove the result we need later on.

\begin{proposition} \label{extr1.prop2}
Suppose that
\begin{itemize}
\item The mapping $\la$ is closed.
\item There exist a net $(u_i)_{i \in I}$ in $N$ and a net $(S_i)_{i \in I}$ in $B(H)$
such that
\begin{enumerate}
\item We have for every $a \in N$ and $i \in I$ that
that $\la(a u_i) = S_i \la(a)$.
\item We have for every $i \in I$ that $\|S_i\| \leq 1$.
\item $(u_i)_{i \in I}$ converges strictly to 1 and $(S_i)_{i \in I}$ converges
strongly to 1.
\end{enumerate}
\end{itemize}
Then $(H,\la,\pi)$ is a GNS-construction for $\vfi$ and $\cF = \cF_\vfi$.

\medskip

We have moreover the following results :

\medskip

Define for every $i \in I$ the element $\om_i \in A^*_+$ such that $\om_i(x) = \vfi(u_i^*
x u_i)$ for $x \in A$. Then we have for every $i \in I$ that $\om_i$ belongs to $\cF$.
Furthermore,

\begin{enumerate}
\item Consider $x \in \Mfi$, then $\bigl(\om_i(x)\bigr)_{i \in I}$ converges to $\vfi(x)$ .
\item Consider $x \in A^+$, then $\vfi(x) = \sup\,\{ \, \om_i(x) \mid i \in I \,\}$ .
\end{enumerate}
\end{proposition}
\begin{demo}
Choose $i \in I$. Then we have for every $a \in N$ that
$$\om_i(a^* a) = \vfi(u_i^* a^* a u_i) = \langle \la(a u_i) , \la(a u_i) \rangle
= \langle S_i \la(a) , S_i \la(a) \rangle \ . \hspace{1.5cm} \text{(*)}$$

Because $\|S_i\| \leq 1$, this implies that $\om_i(a^* a) \leq \|\la(a)\|^2$ for every $a
\in N$. Therefore $\om_i$ belongs to $\cF$.

\begin{itemize}
\item Choose $a \in N$. Because $(S_i)_{i \in I}$ converges strongly to 1, equality (*)
implies that $\bigl(\om_i(a^* a)\bigr)_{i \in I}$ converges to $\|\la(a)\|^2$ which
implies that $\|\la(a)\|^2 \leq \vfi(a^* a)$, so $\vfi(a^* a) = \| \la(a)\|^2$.

\item Consider $a \in A$ such that the net $\bigl( \om_i(a^* a) )_{i \in I}$ is bounded.

We have for every $i \in I$ that $\om_i \in \cF$ which by the definition of $\vfi$ implies
that $$ \|\la(a u_i)\|^2 = \langle \la(a u_i) , \la(a u_i) \rangle
= \vfi(u_i^* a^* a u_i) = \om_i(a^* a) $$
Therefore the net $\bigl(\la(a u_i)\bigr)_{i \in I}$ is bounded.

Because $(a u_i)_{i \in I}$ converges to $a$, the previous lemma implies the existence of
a sequence $(b_n)_{n=1}^\infty$ in $D(\la)$ and $v \in H$ such that $(b_n)_{n=1}^\infty$
converges to $a$ and $(\la(b_n))_{n=1}^\infty$ converges to $v$. Therefore the closedness
of $\la$ implies that $a$ belongs to $N$.

\medskip

Because of the definition of $\vfi$, this last result implies also that $\Nfi \subseteq N$
which implies that $N = \Nfi$.
\end{itemize}
\end{demo}

\section{Relatively invariant one-parameter groups}

In this paper, we will have to use one-parameter groups which are relatively invariant
with respect to some closed mappings arising from the GNS-construction for weights. In
this section, we will gather some results about such objects.

\bigskip

In  this section, we  consider two Banach spaces $E$,$F$, a subspace $N$ of $E$ and a
closed linear map $\la$ from $N$ into $F$. Let $\al$ be a  strongly continuous
one-parameter representation on $E$ such that there exists a strongly continuous
one-parameter representation $u$ on $F$ together with a strictly  positive number
$\lambda$ such that we have for every $t \in \R$ and $a \in N$ that $\al_t(a)$ belongs to
$N$ and $\la(\al_t(a)) = \lambda^\frac{t}{2} \, u_t \la(a)$.

\bigskip

The techniques used in the proof of the following lemma will be frequently used.

\begin{lemma} \label{we2.lem1}
Consider $y,z \in \C$, $n \in \N$ and $a \in N$. Then
$$\int \exp(- n^2(t-z)^2) \, \lambda^{yt} \,\al_t(a) \, dt$$
belongs to $N$ and
$$\la\bigl(\,\int \exp(-n^2(t-z)^2) \, \lambda^{yt} \,  \al_t(a) \, dt\,\bigr)
= \int \exp(-n^2(t-z)^2) \, \lambda^{(y+\frac{1}{2})t} \, u_t \la(a) \, dt \ .$$
\end{lemma}
\begin{demo}
By assumption, we have for every $t \in \R$ that
$\exp(-n^2(t-z)^2) \, \lambda^{yt} \, \al_t(a)$ belongs to $N$
and
$$\la\bigl( \exp(-n^2(t-z)^2) \, \lambda^{yt} \, \al_t(a) \bigr)
= \exp(-n^2(t-z)^2) \, \lambda^{(y+\frac{1}{2})t} \, u_t \la(a) \ .$$
This implies that the function
$$ \R \rightarrow F : t \mapsto
\la\bigl(\exp(-n^2(t-z)^2) \, \lambda^{yt} \, \al_t(a)\bigr)$$ is integrable.
Using lemma \ref{o.lem1}, we get that
$$\int \exp(-n^2(t-z)^2) \, \lambda^{yt} \, \al_t(a) \, dt$$
belongs to $N$ and
\begin{eqnarray*}
\la\bigl(\,\int \exp(-n^2(t-z)^2) \, \lambda^{yt} \, \al_t(a) \, dt\,\bigr)
& = & \int \la\bigl( \exp(-n^2(t-z)^2) \, \lambda^{yt} \, \al_t(a)\bigr) \, dt \\ & = &
\int \exp(-n^2(t-z)^2) \, \lambda^{(y+\frac{1}{2})t} \, u_t \la(a) \, dt .
\end{eqnarray*}
\end{demo}

The following lemma follows easily from the previous one.

\begin{lemma} \label{we2.lem2}
Consider $a \in N$, $n \in \N$ and define
$$b = \frac{n}{\sqrt{\pi}} \int \exp(-n^2 t^2) \, \al_t(a) \, dt \ .$$
Then $b$ is analytic with respect to $\al$ and we have for every $z \in \C$ that
$\al_z(b)$ belongs to $N$
\end{lemma}

\medskip

A first application can be found in the following result :

\begin{proposition} \label{we2.prop5}
Define the set
$$ C = \{\, a \in N \mid a \text{ is analytic with respect to } \al
\text{ and }  \al_z(a) \text{ belongs to } N \text{ for every } z \in \C \,\} \ .
$$
Then \begin{minipage}[t]{14cm} \begin{enumerate}
\item $C$ is a core for $\la$.
\item Let $z$ be a complex number. If $N$ is dense in $A$, then $C$ is a core for $\al_z$.
\end{enumerate} \end{minipage}
\end{proposition}
\begin{demo}
Define for every $n \in \N$ and $a \in N$ the element
$$a(n) = \frac{n}{\sqrt{\pi}} \int \exp(-n^2 t^2) \, \al_t(a) \, dt \ , $$ which belongs
to $C$ by the previous lemma.
\begin{enumerate}
\item Choose $x \in N$. By lemma \ref{we2.lem1}, it follows that
$$\la(x(n)) = \frac{n}{\sqrt{\pi}} \int \exp(-n^2 t^2) \, \lambda^\frac{t}{2} \,
u_t \la(x) \, dt   $$ for every $n \in N$. This implies that
$\bigl(\la(x(n))\bigr)_{n=1}^\infty$ converges to $\la(x)$. It is also clear that
$(x(n))_{n=1}^\infty$ converges to $x$.
\item If $N$ is dense in $A$, proposition \ref{pa.prop1} learns that the
set $\langle\, a(n) \mid a \in N , n \in \N \,\rangle$ is a core for $D(\al_z)$. This
implies immediately that $C$ is a core for $D(\al_z)$ in this case.
\end{enumerate}
\end{demo}

\bigskip

\begin{proposition}  \label{we2.prop4}
Consider $z \in \C$ and $a \in N \cap D(\al_z)$ such that $\al_z(a)$ belongs to $N$. Then
$\la(a)$ belongs to $D(u_z)$ and $u_z \la(a) = \lambda^{-\frac{z}{2}} \, \la(\al_z(a))$.
\end{proposition}
\begin{demo}
Choose $n \in \N$ and define
$$v_n = \frac{n}{\sqrt{\pi}} \int \exp(-n^2 t^2) \, u_t \la(a) \, dt \in F \ .$$
It is clear that $v_n$ belongs to $D(u_z)$ and
$$u_z(v_n) = \frac{n}{\sqrt{\pi}} \int \exp(-n^2 (t-z)^2) \, u_t \la(a) \ .$$
By lemma \ref{we2.lem1}, we know that
$$\frac{n}{\sqrt{\pi}} \int \exp(-n^2 (t-z)^2) \, \lambda^{-\frac{t}{2}}
\, \al_t(a) \, dt $$
belongs to $N$ and
$$\la\bigl(\, \frac{n}{\sqrt{\pi}} \int \exp(-n^2 (t-z)^2) \, \lambda^{-\frac{t}{2}}
\, \al_t(a) \, dt\,\bigr)
= u_z(v_n) \text{\ \ \ \ \ \ (*)}$$

\medskip

Because $a$ belongs to $D(\al_z)$, we have that
$$
\frac{n}{\sqrt{\pi}} \int \exp(-n^2 (t-z)^2) \, \lambda^{-\frac{t}{2}}
\, \al_t(a) \, dt =
\lambda^{-\frac{z}{2}} \, \frac{n}{\sqrt{\pi}} \int \exp(-n^2 t^2) \,
\lambda^{-\frac{t}{2}} \, \al_t(\al_z(a)) \, dt \ .
$$
Using lemma \ref{we2.lem1} once more, this implies that
$$\la\bigl(\,\frac{n}{\sqrt{\pi}} \int \exp(-n^2 (t-z)^2) \, \lambda^{-\frac{t}{2}} \,
\al_t(a) \, dt \, \bigr)
=  \lambda^{-\frac{z}{2}} \, \frac{n}{\sqrt{\pi}} \int \exp(-n^2 t^2) \,  u_t
\la(\al_z(a)) \, dt \ .$$
Comparing this equality with equality (*), we get that
$$u_z(v_n) =   \lambda^{-\frac{z}{2}} \, \frac{n}{\sqrt{\pi}} \int \exp(-n^2 t^2)
\,  u_t \la(\al_z(a)) \, dt \ .$$
This implies that $\bigl(u_z(v_n)\bigr)_{n=1}^\infty$ converges to $\lambda^{-\frac{z}{2}}
\, \la(\al_z(a))$. It is also clear that $(v_n)_{n=1}^\infty$ converges to $\la(a)$.
Therefore, the closedness of $u_z$ implies that $\la(a)$ belongs to $D(u_z)$ and $u_z
\la(a) = \lambda^{-\frac{z}{2}} \, \la(\al_z(a))$.
\end{demo}

\begin{proposition} \label{we2.prop3}
Consider $z \in \C$ and $a \in N \cap D(\al_z)$ such that $\la(a)$ belongs to $D(u_z)$.
Then $\al_z(a)$ belongs to $N$ and $\la(\al_z(a)) = \lambda^\frac{z}{2} \, u_z \la(a)$.
\end{proposition}
\begin{demo}
Choose $n \in \N$ and define
$$b_n = \frac{n}{\sqrt{\pi}} \int \exp(- n^2 t^2) \, \al_t(\al_z(a)) \, dt  \ . $$
We also have that
$$b_n = \frac{n}{\sqrt{\pi}} \int \exp(-n^2 (t-z)^2) \, \al_t(a) \, dt \ .$$
Therefore, lemma \ref{we2.lem1} implies that
$b_n$ belongs to $N$ and
$$\la(b_n) = \frac{n}{\sqrt{\pi}} \int \exp(-n^2 (t-z)^2) \, \lambda^\frac{t}{2}
\,  u_t \la(a) \, dt \ .  $$
Because $\la(a)$ belongs to $D(u_z)$, we get from the previous equation
that
$$ \la(b_n)
= \lambda^\frac{z}{2} \, \frac{n}{\sqrt{\pi}} \int
\exp(- n^2 t^2) \, \lambda^\frac{t}{2} \, u_t(u_z \la(a)) \, dt
\ .$$

It is clear from this last equation that $\bigl(\la(b_n)\bigr)_{n=1}^\infty$ converges to
$\lambda^\frac{z}{2} \, u_z \la(a)$. It is also clear from the definition that
$(b_n)_{n=1}^\infty$ converges to $\al_z(a)$. Therefore, the closedness of $\la$ implies
that $\al_z(a)$ belongs to $N$ and $\la(\al_z(a)) = \lambda^\frac{z}{2} \, u_z \la(a)$.
\end{demo}

\begin{corollary}
Consider $z \in \C$ and $a \in N \cap D(\al_z)$. Then $\la(a)$ belongs
to $D(u_z)$ $\Leftrightarrow$ $\al_z(a)$ belongs to $N$.
\end{corollary}

\bigskip

\section{KMS-weights arising from certain GNS-constructions} \label{we3}

In \cite{Verd} and \cite{Q-V}, Jan Verding introduced a construction procedure for weights
starting from a GNS-construction. We want to repeat this procedure, but we want to end up
with a KMS-weight. In order to do so, we have to impose stronger conditions on the
ingredients of the construction procedure. In a last part, we describe the natural case
where a KMS-weight is obtained from a left Hilbert-algebra.

\bigskip

First we prove some technical results.

\medskip

Consider a \cst-algebra $A$, a Hilbert space $H$, a dense left ideal $N$ of $A$ and a
closed linear mapping $\la$ from $N$ into $H$. Let $\al$ be a norm continuous
one-parameter group on $A$ and suppose there exists a positive injective operator $\nab$
in $H$ and a strictly posive number $\lambda$ such that $\al_t(a) \in  N$ and
$\la(\al_t(a)) = \lambda^\frac{t}{2} \, \nab^{it} \la(a)$ for every $t \in \R$ and $a \in
N$.

\begin{lemma}   \label{we2.lem3}
Consider $a \in N \cap N^*$ and $n \in \N$ and define
$$b = \frac{n}{\sqrt{\pi}} \int \exp(-n^2 t^2) \, \al_t(a) \, dt \ .$$
Then $b$ is analytic with respect to $\al$ and $\al_z(b)$ belongs to $N \cap N^*$ for
every $z \in \C$.
\end{lemma}
\begin{demo}
Choose $z \in \C$. By lemma \ref{we2.lem2}, we know already that $\al_z(b)$
belongs to $N$.

We also have that $a^*$ belongs to $N$ and that
$$b^* = \frac{n}{\sqrt{\pi}} \int \exp(-n^2 t^2) \, \al_t(a^*) \, dt \ ,$$
so lemma \ref{we2.lem2} implies in this case that $\al_{\,\overline{z}}(b^*)$ belongs to
$N$. Because $\al_z(b)^* = \al_{\,\overline{z}}(b^*)$, we see that $\al_z(b)$ belongs to
$N^*$.
\end{demo}

\medskip

Define the set
$$
C = \{\, a \in N \cap N^* \mid a \text{ is analytic with respect to } \al \text{ and }
\al_z(a) \text{ belongs to } N \cap N^* \text{ for every } z \in \C
\,\} \ .
$$

\medskip

\begin{proposition} \label{we2.prop2}
We have the following properties :
\begin{enumerate}
\item  For every $x \in N \cap N^*$, there exists a sequence $(x_n)_{n=1}^\infty$ in $C$
such that $(x_n)_{n=1}^\infty$ converges to $x$, $(\la(x_n))_{n=1}^\infty$ converges to
$\la(x)$ and $(\la(x_n^*))_{n=1}^\infty$ converges to $\la(x^*)$.
\item If $N \cap N^*$ is a core for $\la$, then $C$ is a core for $\la$.
\item Let $z$ be a complex number. Then $C$ is a core for $\al_z$.
\end{enumerate}
\end{proposition}
\begin{demo}
Define for every $n \in \N$ and $a \in N \cap N^*$ the element
$$a(n) = \frac{n}{\sqrt{\pi}} \int \exp(-n^2 t^2) \, \al_t(a) \, dt \ , $$ which belongs
to $C$ by the previous lemma.
\begin{enumerate}
\item Choose $x \in N \cap N^*$. Again, we see immediately that
$(x(n))_{n=1}^\infty$ converges to $x$. Because $x$ belongs to $N$, the proof of
proposition \ref{we2.prop5} learns us that $\bigl(\la(x(n))\bigr)_{n=1}^\infty$ converges
to $\la(x)$.

We have for every $n \in \N$ that $x(n)^* = x^*(n)$, so the same proof tells us that
$\bigl(\la(x(n)^*)\bigr)_{n=1}^\infty$ converges to $\la(x^*)$.
\item This follows immediately from the previous statement.
\item This is proven in the same way as the second statement of proposition \ref{we2.prop5}.
\end{enumerate}
\end{demo}

\bigskip

The proof of the following result is due to A. Van Daele and J. Verding.

\begin{proposition} \label{pa2.prop2}
Consider $z \in \C$. Then there exists a net $(u_k)_{k \in K}$ in $C \cap A^+$ consisting
of analytic elements for $\al$ and such that
\begin{itemize}
\item We have for every $k \in K$ that $\|u_k\| \leq 1$ and $\|\al_z(u_k)\| \leq 1$.
\item The nets $(u_k)_{k \in K}$ and $(\al_z(u_k))_{k \in K}$ converge strictly to 1.
\end{itemize}
\end{proposition}
\begin{demo}
Because $N$ is a dense left ideal in $A$, we know that there exists a net $(e_q)_{q
\in Q}$ in $N \cap A^+$ such that $\|e_q\| \leq 1$ for every $q \in Q$ and such that
$(e_q)_{q \in Q}$ converges strictly to 1.

\medskip

For every $q \in Q$ and $\sde > 0$, we  define the element
$$e_q(\sde) = \frac{\sde}{\sqrt{\pi}} \int \exp(-\sde^2 t^2) \, \al_t(e_q) \, dt
\ \ \  \in \ \ A^+$$
which is clearly analytic with respect to $\al$ and satisfies $\|e_q(\sde)\|
\leq 1$ and $\|\al_z(e_q(\sde))\| \leq \exp(\sde^2 \, (\text{Im }z)^2)$.

Lemma \ref{we2.lem3} implies that $e_q(\sde) \in C$ for every $q \in Q$ and $\sde > 0$.

\medskip

By proposition \ref{pa2.prop6}, we have for every $\sde > 0$ that
$\bigl(e_q(\sde)\bigr)_{q \in Q}$ and $\bigl(\al_z(e_q(\sde))\bigr)_{q \in Q}$ are bounded
nets which converge strictly to 1.

\medskip

Let us now define the set $$K = \{ \, (F,n) \mid F \text{\ is a finite subset of  \ } A
\text{\ \ and \ \ } n \in \N  \, \} .$$ On $K$ we put an order such that
$$(F_1,n_1) \leq (F_2,n_2) \ \
\Leftrightarrow \ \ F_1 \subseteq F_2 \text{\ and \ } n_1 \leq n_2 $$
for every $(F_1,n_1),(F_2,n_2) \in K$.
In this way, $K$ becomes a directed set.

\medskip

Let us fix $k=(F,n) \in K$.

Firstly, there exist an element $\sde_k > 0$ such that $|\exp(-\sde_k^2 \, (\text{Im
}z)^2) - 1 | \leq \frac{1}{2 n} \frac{1}{\|x\|+1}$ for every $x \in F$.

Secondly, there exist $q_k \in Q$ such that
\begin{itemize}
\item We have that $\| e_{q_k}(\sde_k) \, x - x \| \leq \frac{1}{2 n}$
      for evey $x \in F$.
\item We have that $\| \al_z(e_{q_k}(\sde_k)) \, x - x \|
      \leq \frac{1}{2 n}$ for every $x \in F$.
\item We have that $\| x \, \al_z(e_{q_k}(\sde_k)) - x \|
      \leq \frac{1}{2 n}$ for every $x \in F$.
\end{itemize}
Now we define $u_k = \exp(-\sde_k^2\,(\text{Im }z)^2) \,\, e_{q_k}(\sde_k) \in C \cap
A^+$.

\medskip

It follows that $u_k$ is analytic with respect to $\al$ and  $\|u_k\| \leq 1$. We have
also that
\begin{eqnarray*}
\|\al_z(u_k)\| & = & \exp(-\sde_k^2\,(\text{Im }z)^2) \,\, \| \al_z(e_{q_k}(\sde_k))\| \\
& \leq & \exp(-\sde_k^2\,(\text{Im }z)^2) \, \exp(-\sde_k^2 \,(\text{Im }z)) = 1
\end{eqnarray*}

\medskip

Now we prove that $(u_k)_{k \in K}$ converges strictly to 1.

Choose $x \in A$ and $\vep > 0$. Then there exists $n_0$ in $\N$ such that $\frac{1}{n_0}
\leq \vep$. Put $k_0 = (\{x\},n_0) \in K$.

Take $k=(F,n) \in K$ such that $k \geq k_0$, so $x$ belongs to $F$
and $n_0 \leq n$. Therefore,
\begin{eqnarray*}
& & \| u_k \, x - x \| =
\| \exp(-\sde_k^2 \,(\text{Im }z)^2) \,e_{q_k}(\sde_k) \, x - x \| \\
& & \spat \leq  \| \exp(-\sde_k^2 \, (\text{Im }z)^2) \, e_{q_k}(\sde_k) \, x
- \exp(-\sde_k^2\, (\text{Im }z)^2) \, x \|
+ \| \exp(-\sde_k^2 \, (\text{Im }z)^2) \, x - x \| \\
& & \spat = \exp(-\sde_k^2 \, (\text{Im }z)^2) \,\,
\| e_{q_k}(\sde_k) \, x - x \| + | \exp(-\sde_k^2 \, (\text{Im }z)^2) - 1 | \,\, \|x \| \\
& & \spat \leq \| e_{q_k}(\sde_k) \, x - x \| + \frac{1}{2 n (\|x\|+1)} \, \|x\|
\leq  \frac{1}{2 n} + \frac{1}{2 n} = \frac{1}{n} \leq \vep \  .
\end{eqnarray*}
Because $(u_k)_{k \in K}$ consists of selfadjoint elements, this implies that $(u_k)_{k \in K}$ converges stricly to 1.

Completely analogously, one proves that $(\al_z(u_k)\, x)_{k \in K}$ converges to $x$ for every $x \in A$ (just replace $u_k$ by $\al_z(u_k)$ in the proof above).
Similarly one proves that $(x \, \al_z(u_k))_{k \in K}$ converges to $x$ for every $x \in A$. Consequently, we have that
$(\al_z(u_k))_{k \in K}$ converges strictly to 1.
\end{demo}

\bigskip\medskip

\subsection{The first construction procedure.}

In the first construction procedure,
we will use the following ingredients:

Consider a \cst-algebra $A$ and a Hilbert space $H$. Let $N$ be a dense left ideal of $A$,
$\la$ a closed linear mappping from $N$ into $H$ with dense range and $\pi$ a
non-degenerate representation of $A$ on $H$ such that $\pi(x) \la(a) = \la(x a)$ for every
$x \in A$ and $a \in N$.

Furthermore, assume the existence of
\begin{itemize}
\item a norm continuous one-parameter group $\si$ on $A$
\item an injective positive operator $\nab$ in $H$
\end{itemize}
such that:
\begin{itemize}
\item We have for every $a \in N$ and $t \in \R$ that $\si_t(a)$ belongs
      to $N$ and $\la(\si_t(a))=\nab^{it} \la(a)$.
\item There exists a core $K$ for $\la$ such that
      \begin{enumerate}
      \item $K \subseteq D(\si_\frac{i}{2})$, \
      $\si_{\frac{i}{2}}(K)^* \subseteq N$ and $\| \la(x) \| = \|\la(
      \si_{\frac{i}{2}}(x)^*)\|$ for every $x \in K$.
      \item $\si_t(K) \subseteq K$ for every $t \in \R$.
      \end{enumerate}
\end{itemize}

\medskip

Our first objective is to extend the property about $K$ to $N \cap D(\si_\frac{i}{2})$.

\medskip

In this section, we will use the following notation:
For every $a \in A$ and $n \in \N$, we put
$$a(n) = \frac{n}{\sqrt{\pi}} \int \exp(-n^2 t^2)
\, \si_t(a) \, dt \ ,$$
it is clear that $a(n)$ is analytic with respect to $\si$.

\medskip

Furthermore, we know from lemma \ref{we2.lem1} that we have for every $a \in N$ and $n
\in \N$ that $a(n)$ belongs to $N$ and
\begin{equation}
\la(a(n)) = \frac{n}{\sqrt{\pi}} \int \exp(-n^2 t^2) \, \nab^{it} \la(a) \, dt
\label{we3.eq1} .
\end{equation}

\bigskip

\begin{lemma}
We have the following properties:
\begin{enumerate}
\item We have for every $x \in K$ and $n \in \N$ that $x(n) \in N \cap
D(\si_\frac{i}{2})$ and $\si_\frac{i}{2}(x(n))^* \in N$.
\item
We have for every $x,y \in K$ and $m,n \in \N$ that
$$\langle \la(y(n)) , \la(x(m)) \rangle
= \langle \la(\si_\frac{i}{2}(x(m))^*),
\la(\si_\frac{i}{2}(y(n))^*) \rangle \ .$$
\end{enumerate}
\end{lemma}
\begin{demo}
\begin{enumerate}
\item
Choose $x \in K$ and $n \in \N$. By equation \ref{we3.eq1}, we know already that $x(n)$
belongs to $N$. We also have that $x(n)$ belongs to $D(\si_\frac{i}{2})$. The fact that
$x$ belongs to $D(\si_\frac{i}{2})$ implies that
$$ \si_\frac{i}{2}(x(n))^*
= (\si_\frac{i}{2}(x)(n))^*
= \si_\frac{i}{2}(x)^*(n) \ \ \ \ \ \ \ \ \ (a)$$
Therefore, equation \ref{we3.eq1} and  the fact that $\si_\frac{i}{2}(x)^* \in N$ imply that
$\si_\frac{i}{2}(x(n))^*$ belongs to $N$.
\item
Choose $x,y \in K$ and $m,n \in \N$. Using equation (a) from the first part of the proof
and equation \ref{we3.eq1} once more, we get the following equalities
\begin{eqnarray*}
& & \langle \la(\si_\frac{i}{2}(x(m))^*) ,
\la(\si_\frac{i}{2}(y(n))^*) \rangle
= \langle \la(\si_\frac{i}{2}(x)^*(m)) ,
\la(\si_\frac{i}{2}(y)^*(n)) \rangle \\
& & \spat =  \frac{m \, n}{\pi} \int \int
\exp(-(m^2 s^2+n^2 t^2)) \,
\langle \nab^{is} \la(\si_\frac{i}{2}(x)^*) ,
\nab^{it} \la(\si_\frac{i}{2}(y)^*) \rangle \, ds \, dt    \hspace{1.5cm} \text{(b)}
\end{eqnarray*}
By polarisation, we get that
$$\langle \la(a) , \la(b) \rangle  =
\langle \la(\si_\frac{i}{2}(b)^*) , \la(\si_\frac{i}{2}(a)^*)
\rangle $$
for every $a,b \in K$.

Fix $s,t \in \R$. By assumption, we know that $\si_s(x),\si_t(y)$
belong to $K$. So, using the previous equality, we see that
\begin{eqnarray*}
& & \langle \nab^{is} \la(\si_\frac{i}{2}(x)^*) ,
\nab^{it} \la(\si_\frac{i}{2}(y)^*) \rangle
=  \langle \la(\si_s(\si_\frac{i}{2}(x)^*)) ,
\la(\si_t(\si_\frac{i}{2}(y)^*)) \rangle  \\
&  & \spat = \langle \la(\si_\frac{i}{2}(\si_s(x))^*) ,
\la(\si_\frac{i}{2}(\si_t(y))^*) \rangle
=  \langle \la(\si_t(y)),\la(\si_s(x))\rangle = \langle \nab^{it} \la(y) ,
\nab^{is} \la(x) \rangle
\end{eqnarray*}
Substituting this equality into equation (b) gives us that
\begin{eqnarray*}
\langle \la(\si_\frac{i}{2}(x(m))^*) , \la(\si_\frac{i}{2}(y(n))^*) \rangle
& = & \frac{m \, n}{\pi} \int \int \exp(-(m^2 s^2+n^2 t^2)) \, \langle \nab^{it} \la(y) ,
\nab^{is} \la(x) \rangle \, ds \, dt \\
&  = & \langle \la(y(n)) , \la(x(m)) \rangle \ ,
\end{eqnarray*}
where we used equation \ref{we3.eq1} in the last equality.
\end{enumerate}
\end{demo}

\begin{lemma}
We have the following properties:
\begin{enumerate}
\item
We have for every $x \in N$ and $n \in \N$ that $x(n) \in N \cap
D(\si_\frac{i}{2})$ and $\si_\frac{i}{2}(x(n))^* \in N$.
\item
We have for every $x,y \in N$ and $m,n \in \N$ that
$$\langle \la(y(n)) , \la(x(m)) \rangle
= \langle \la(\si_\frac{i}{2}(x(m))^*),
\la(\si_\frac{i}{2}(y(n))^*) \rangle \ .$$
\end{enumerate}
\end{lemma}
\begin{demo}
\begin{enumerate}
\item
Choose $x \in N$ and $n \in \N$. We know by equation \ref{we3.eq1} that $x(n)$ belongs to
$N$. We have also that $x(n)$ belongs to $D(\si_\frac{i}{2})$.

Because $K$ is a core for $\la$, there exists a sequence $(x_k)_{k=1}^\infty$ in $K$ such
that $(x_k)_{k=1}^\infty \rightarrow x$ and $\bigl(\la(x_k)\bigr)_{k=1}^\infty \rightarrow
\la(x)$. From the previous lemma we know that $x_k(n)$ belongs to $N \cap
D(\si_\frac{i}{2})$ and that $\si_\frac{i}{2}(x_k(n))^*$ belongs to $N$ for every $k \in
\N$.

For every $k \in \N$, we have that
$$ \si_\frac{i}{2}(x_k(n))^*
= \frac{n}{\sqrt{\pi}} \int \exp(-n^2 (t+\frac{i}{2})^2) \, \si_t(x_k^*) \, dt $$
and similarly,
$$ \si_\frac{i}{2}(x(n))^*
= \frac{n}{\sqrt{\pi}} \int \exp(-n^2 (t+\frac{i}{2})^2) \, \si_t(x^*) \, dt \ .$$
Comparing these two equations and remembering that $(x_k)_{k=1}^\infty$ converges to $x$,
we get that \newline $\bigl(\,\si_\frac{i}{2}(x_k(n))^*\,\bigr)_{k=1}^\infty$ converges to
$\si_\frac{i}{2}(x(n))^*$. \ \ \ (a)

Because $\bigl(\la(x_k)\bigr)_{k=1}^\infty$ converges to $\la(x)$, it is clear from equation
\ref{we3.eq1} that $\bigl(\,\la(x_k(n))\,\bigr)_{k=1}^\infty$ converges to $\la(x(n))$. By
the previous lemma, we have for every $k,l \in \N$ that
$$ \| \la(\si_\frac{i}{2}(x_l(n))^*)
-\la(\si_\frac{i}{2}(x_k(n))^*) \|
= \| \la(x_l(n) ) - \la(x_k(n)) \|
$$
Therefore, the convergence of $\bigl(\,\la(x_k(n))\,\bigr)_{k=1}^\infty$ implies that
$\bigl(\,\la(\si_\frac{i}{2}(x_k(n))^*)\,\bigr)_{k=1}^\infty$ is a Cauchy sequence and
hence convergent.  \ \ \ (b)

Using the closedness of $\la$ and conclusions (a) and (b), we get that
$\si_\frac{i}{2}(x(n))^*$ belongs to $N$ and that
$\bigl(\,\la(\si_\frac{i}{2}(x_k(n))^*)\,\bigr)_{k=1}^\infty$ converges to
$\la(\si_\frac{i}{2}(x(n))^*)$.

\item
Choose $x,y \in N$ and $m,n \in \N$.
Then there exists sequences $(x_k)_{k=1}^\infty$, \ $(y_k)_{k=1}^\infty$
in $K$ such that $(x_k)_{k=1}^\infty \rightarrow x$, \  $(y_k)_{k=1}^\infty
\rightarrow y$, \ $\bigl(\la(x_k)\bigr)_{k=1}^\infty \rightarrow \la(x)$
and $\bigl(\la(y_k)\bigr)_{k=1}^\infty \la(y) \rightarrow \la(y)$.

By the first part of the proof, we know that $\bigl(\,\la(x_k(m))\,\bigr)_{k=1}^\infty
\rightarrow \la(x(m))$, \ $\bigl(\,\la(y_k(n))\,\bigr)_{k=1}^\infty \rightarrow \la(y(n))$, \
$\bigl(\,\la(\si_\frac{i}{2}(x_k(m))^*)\,\bigr)_{k=1}^\infty
\rightarrow \la(\si_\frac{i}{2}(x(m))^*)$
and $\bigl(\,\la(\si_\frac{i}{2}(y_k(n))^*)\,\bigr)_{k=1}^\infty
\rightarrow \la(\si_\frac{i}{2}(y(n))^*)$.
By the previous lemma, we know that
$$\langle \la(y_k(n)) , \la(x_k(m)) \rangle =
\langle \la(\si_\frac{i}{2}(x_k(m))^*) ,
\la(\si_\frac{i}{2}(y_k(n))^*) \rangle $$
for every $k \in \N$. This implies that
$$\langle \la(y(n)), \la(x(m)) \rangle =
\langle \la(\si_\frac{i}{2}(x(m))^*) ,
\la(\si_\frac{i}{2}(y(n))^*) \rangle \ .$$
\end{enumerate}
\end{demo}

\medskip

At last, we arive at a final form.

\begin{proposition} \label{we3.prop2}
Consider $x \in N \cap D(\si_\frac{i}{2})$. Then
$\si_\frac{i}{2}(x)^*$ belongs to $N$ and
$\|\la(x)\| = \| \la(\si_\frac{i}{2}(x)^*) \|$.
\end{proposition}
\begin{demo}
By the previous lemma, we have for every $n \in \N$ that $x(n)$ belongs to $N \cap
D(\si_\frac{i}{2})$,  $\si_\frac{i}{2}(x(n))^*$ belongs to $N$ and $\|\la(x(n))\| = \|
\la(\si_\frac{i}{2} (x(n))^*)\|$.

By equation \ref{we3.eq1}, we have for every $n \in \N$ that
$$\la(x(n)) = \frac{n}{\sqrt{\pi}} \int
\exp(-n^2 t^2) \, \nab^{it} \la(x) \, dt  $$
which implies that $\bigl(\,\la(x(n))\,\bigr)_{n=1}^\infty
\rightarrow \la(x) $.

By the previous lemma, we have for every $m,n \in \N$ that
$$\| \la(\si_\frac{i}{2}(x(m))^*)
- \la(\si_\frac{i}{2}(x(n))^*) \|
= \| \la(x(m)) - \la(x(n)) \| \ .$$
Hence, the convergence of $\bigl(\,\la(x(n)) \,\bigr)_{n=1}^\infty$ implies that
$\bigl(\,\la(\si_\frac{i}{2}(x(n))^*)\,\bigr)_{n=1}^\infty$ is a Cauchy sequence and hence
convergent.

Because $x$ belongs to $D(\si_\frac{i}{2})$, we have for every $n \in \N$ that $$
\si_\frac{i}{2}(x(n))^* = \frac{n}{\sqrt{\pi}} \int \exp(-n^2 t^2) \, \si_t(\si_\frac{i}{2}
(x)^*) \, dt \, , $$ from which it follows easily  that $\bigl(\,
\si_\frac{i}{2}(x(n))^*\,\bigr)_{n=1}^\infty$ converges to $\si_\frac{i}{2}(x)^*$.

Hence, the closedness of $\la$ implies that $\si_\frac{i}{2}(x)^*$ belongs to $N$ and
$\bigl(\,\la(\si_\frac{i}{2}(x(n))^*)\,\bigr)_{n=1}^\infty \rightarrow
\la(\si_\frac{i}{2}(x)^*)$.

Because $\|\la(\si_\frac{i}{2}(x(n))^*)\| = \|\la(x(n))\|$ for every $n \in \N$, we will
also have that $\| \la(\si_\frac{i}{2}(x)^*)\| = \|\la(x)\|$.
\end{demo}

Using the fact that $\si_\frac{i}{2}(\si_\frac{i}{2}(x^*)^*) = x$
for every $x \in D(\si_\frac{i}{2})$, the following corollary follows readily.

\begin{corollary}
Consider $x \in D(\si_\frac{i}{2})$. Then $x$ belongs to
$N$ if and only if $\si_\frac{i}{2}(x)^*$ belongs to $N$.
\end{corollary}

Because $K \subseteq N \cap D(\si_\frac{i}{2})$, we have that
$N \cap D(\si_\frac{i}{2})$ is a core for $\la$.
Therefore,  $\la(N \cap D(\si_\frac{i}{2}))$ will certainly be dense
in $H$. This allows us to introduce the following definition.

\begin{definition}
We define the anti-unitary operator $J$ on $H$ such that
$J \la(x) = \la(\si_\frac{i}{2}(x)^*)$ for every $x \in N \cap
D(\si_\frac{i}{2})$.
\end{definition}

It follows immediately that $J \la(N \cap D(\si_\frac{i}{2}))
\subseteq \la(N \cap D(\si_\frac{i}{2}))$ and that
$J( J \la(x)) = \la(x)$ for every $x \in N \cap D(\si_\frac{i}{2}))$. This implies that
$J^2=1$.

\medskip

The next proposition is crucial in the construction of a weight out of our ingredients.

\begin{proposition}  \label{we3.prop1}
Consider $a \in D(\si_\frac{i}{2})$ and $x \in N$. Then
$x a$ belongs to $N$ and $\la(x a) = J \pi(\si_\frac{i}{2}(a))^* J
\la(x)$
\end{proposition}
\begin{demo}
Choose $b,y \in N \cap D(\si_\frac{i}{2})$.
Then we have that $y b \in N \cap D(\si_\frac{i}{2})$.
So, by the definition of $J$, we have that
\begin{eqnarray*}
J \la(y b)  & = & \la(\si_\frac{i}{2}(y b)^*)
=  \la(\si_\frac{i}{2}(b)^* \si_\frac{i}{2}(y)^*) \\
& = & \pi(\si_\frac{i}{2}(b)^*) \la(\si_\frac{i}{2}(y)^*)
= \pi(\si_\frac{i}{2}(b))^* J \la(y) \ .
\end{eqnarray*}
Therefore, $\la(y b) = J \pi(\si_\frac{i}{2}(b))^* J \la(y)$.

\medskip

Proposition \ref{we2.prop5} implies the existence of sequences $(a_n)_{n=1}^\infty$,
$(x_n)_{n=1}^\infty$ in $N \cap D(\si_\frac{i}{2})$ such that $(a_n)_{n=1}^\infty$
converges to $a$, $(x_n)_{n=1}^\infty$ converges to $x$,
$(\si_\frac{i}{2}(a_n))_{n=1}^\infty$ converges to $\si_\frac{i}{2}(a)$ and
$(\la(x_n))_{n=1}^\infty$ converges to $\la(x)$.

This implies immediately that $(x_n a_n)_{n=1}^\infty$ converges to $x a$.

By the first part of this proof, we know that $x_n a_n$ belongs to $N$ and $$\la(x_n a_n)
= J \pi(\si_\frac{i}{2}(a_n))^* J \la( x_n)$$ for every $n \in \N$. This implies that
$$\bigl(\,\la(x_n a_n)\,\bigr)_{n=1}^\infty \rightarrow J \pi(\si_\frac{i}{2}(a))^*
J \la(x) \ .$$
The closedness of  $\la$ implies that $x a$ belongs to $N$ and $\la(x a) = J
\pi(\si_\frac{i}{2}(a))^* J \la(x)$.
\end{demo}

\medskip

Define the set

The following lemma is due to J. Verding and can be found in \cite{Verd}.

\begin{lemma} \label{we3.lem1}
There exists a net $(u_k)_{k \in K} \in  N \cap A^+$ and a net $(T_k)_{k \in K}$ in
$\pi(A)'$ such that
\begin{enumerate}
\item We have for every $k \in K$ that $\|u_k\| \leq 1$ and $\|T_k\| \leq 1$.
\item We have for every $k \in k$ that $u_k$ is analytic with respect to
      $\si$.
\item $(u_k)_{k \in K}$ converges strictly to 1 and
$(T_k)_{k \in K}$ converges strongly$^*$ to 1.
\item For every $a \in N$ and $k \in K$, we have that $\la(a u_k) = T_k \la(a)$.
\end{enumerate}
\end{lemma}
\begin{demo} By proposition \ref{pa2.prop2} we get the existence of a net
$(u_k)_{k \in K}$ in $N \cap A^+$ consisting of elements which are analytic with respect
to $\si$ and such that
\begin{itemize}
\item We have for every $k \in K$ that $\|u_k\| \leq 1$ and
$\|\si_\frac{i}{2}(u_k)\| \leq 1$.
\item The nets $(u_k)_{k \in K}$ and $(\si_\frac{i}{2}(u_k))_{k \in K}$ converge
strictly to 1.
\end{itemize}

Choose $k \in K$ and define the element $T_k = J \pi(\si_\frac{i}{2}(u_k))^* J \in
B(H)$. We have immediately that $\|T_k\| \leq 1$. Proposition \ref{we3.prop1} implies that
$\la(a u_k) = T_k \la(a)$ for every $a \in N$.

We have also that $T_k$ belongs to $\pi(A)'$ :
Choose $a \in A$. Then we have for every $b \in N$ that
$$\pi(a) T_k \la(b) = \pi(a) \la(b u_k) = \pi(a b ) \la(u_k)
= T_k \la(ab) = T_k \pi(a) \la(b) \ .$$
Hence, $\pi(a) T_k = T_k \pi(a)$.

\medskip

It is clear that $(T_k)_{k \in K}$ converges strongly$^*$ to 1.
\end{demo}

\medskip

This proposition allows us to use the construction procedure of J.Verding to get a weight
(see definition \ref{extr1.def1} and proposition \ref{extr1.prop2}).

\begin{definition} \label{we3.def1}
Define the set
$${\cal F} = \{\, \om \in A_+^* \mid \om(a^* a) \leq  \|\la(a)\|^2 \text{\ for every \ }
a \in N \,\} \ .$$
We define the mapping $\vfi$ from $A^+$ to $[0,\infty]$ such that
$$ \vfi(x) = \sup \{\, \om(x) \mid \om \in {\cal F} \,\} $$
for every $x \in A^+$. Then $\vfi$ is a densely defined lower semi-continuous weight on
$A$. We have in this case also that $\Nfi = N$ and $\vfi(b^* a) = \langle \la(a) , \la(b)
\rangle$ for every $a,b \in N$. Therefore, $(H,\pi,\la)$ is a GNS-construction for $\vfi$.
\end{definition}

Proposition \ref{we3.prop2} implies immediately that $\vfi$ is a KMS-weight with respect
to $\si$.

\bigskip\bigskip

\subsection{The second construction procedure.}

In the next part of this section, we want to give a second construction procedure for
KMS-weights, starting with slightly different conditions on the ingredients than in the
previous case.

\medskip

In this section we will use the following ingredients:

Consider a \cst-algebra $A$ and a Hilbert space $H$. Let $N$ be a dense left ideal of $A$,
$\la$ a closed linear mappping from $N$ into $H$ with dense range and $\pi$ a
non-degenerate representation of $A$ on $H$ such that $\pi(x) \la(a) = \la(x a)$ for every
$x \in A$ and $a \in N$.

Furthermore, assume the existence of
\begin{itemize}
\item a norm continuous one-parameter group $\si$ on $A$
\item an injective positive operator $\nab$ in $H$
\item an anti-unitary operator $I$ on $H$
\end{itemize}
such that:
\begin{itemize}
\item We have for every $a \in N$ and $t \in \R$ that $\si_t(a)$ belongs
      to $N$ and $\la(\si_t(a))=\nab^{it} \la(a)$.
\item We have for every $x \in N$ and $a \in D(\si_\frac{i}{2})$ that
      $x a$ belongs to $N$ and  $\la(x a)  = I \pi(\si_\frac{i}{2}(a))^* I^* \la(x)$.
\end{itemize}

If we examine the proof of lemma \ref{we3.lem1} we see that we only use proposition
\ref{we3.prop1} and the fact that $\la$ is invariant under $\si$, so the conclusion of
lemma \ref{we3.lem1} still holds in this case. This allows us to define a weight $\vfi$ in
the same way as before (see definition \ref{extr1.def1} and proposition
\ref{extr1.prop2}).

\medskip

\begin{definition} \label{we3.def2}
Define the set
$${\cal F} = \{\, \om \in A_+^* \mid \om(a^* a) \leq  \|\la(a)\|^2 \text{\ for every \ }
 a \in N \,\} \ .$$
We define the mapping $\vfi$ from $A^+$ to $[0,\infty]$ such that
$$ \vfi(x) = \sup \{\, \om(x) \mid \om \in {\cal F} \,\} $$
for every $x \in A^+$. Then  $\vfi$ is a densely defined lower semi-continuous weight on
$A$. We have that $\Nfi = N$ and $\vfi(b^* a) = \langle \la(a) , \la(b) \rangle$ for every
$a,b \in N$. Therefore, $(H,\pi,\la)$ is a GNS-construction for $\vfi$.
\end{definition}
\medskip

We still have to prove that $\vfi$ is a KMS-weight with respect to $\si$.

\medskip

Define the set
$$
 C = \{\, a \in N \cap N^* \mid a \text{ is analytic with respect to } \si  \text{ and }
 \si_z(a) \text{ belongs to } N \cap N^* \text{ for every } z \in \C \,\} \ .
$$
Then $C$ is a sub-$^*$-algebra of $A$ such that $C \subset N \cap D(\si_\frac{i}{2})$ and
$\si_\frac{i}{2}(C)^* \subseteq N$. It is also clear that $\si_z(C) \subseteq C$ for every
$z \in \C$.

\begin{lemma} \label{we3.lem2}
We have that $C$ is a core for $\la$.
\end{lemma}
\begin{demo}
We prove first that $N \cap N^*$ is a core for $\la$.

\begin{list}{}{\setlength{\leftmargin}{.4 cm}}

\item Choose $x \in N$. Because $N^*$ is dense in $A$, there exists a bounded net
$(e_k)_{k \in K}$ in $N^*$ such that $(e_k)_{k \in K}$ converges strictly to 1. For every
$k \in K$, we have that $e_k a$ belongs to $N^* N$, which is a subset of $N \cap N^*$.
Because $\la(e_k a) = \pi(e_k) \la(a)$ for every $k \in K$, we see that $\bigl(\la(e_k a)\bigr)_{k \in K}$
converges to $\la(a)$. It is also clear that $(e_k a)_{k \in K}$ converges to $a$.

\end{list}

Therefore, proposition \ref{we2.prop2} implies that $C$ is a core for $\la$.
\end{demo}

\begin{proposition} \label{we3.prop3}
We have that $\vfi$ is a KMS-weight with modular group $\si$.
\end{proposition}
\begin{demo}
Choose $a \in C$.

Take $x \in C$. Because $a$ belongs to $D(\si_\frac{i}{2})$, we have by assumption that
\begin{eqnarray*}
& & \langle \pi(x) \la(a) , \la(a) \rangle
=   \langle I \pi(\si_\frac{i}{2}(a))^* I^*  \la(x) , \la(a) \rangle
=  \langle \la(x) , I \pi(\si_\frac{i}{2}(a)) I^* \la(a) \rangle \\
& & \spat =  \langle \la(x) , I \pi(\si_{-\frac{i}{2}}(a^*))^* I^* \la(a)  \rangle
=  \langle \la(x) , I \pi(\si_\frac{i}{2}(\si_{-i}(a^*)))^* I^* \la(a)
\rangle \\
&  & \spat = \langle \la(x) , \pi(a) \la(\si_{-i}(a^*)) \rangle
=  \langle \la(a^* x) , \la(\si_{-i}(a^*)) \rangle \ ,  \text{\ \ \ \ \ \ (*)}
\end{eqnarray*}
where, in the second last equality, we used one of the assumptions of this second
construction procedure once again. We know that $\si_{-\frac{i}{2}}(a^*)$ belongs to $N
\cap D(\si_{-\frac{i}{2}})$ and $\si_{-\frac{i}{2}}(\si_{-\frac{i}{2}}(a^*)) =
\si_{-i}(a^*) \in N$. Hence, proposition \ref{we2.prop4} implies that
$\la(\si_{-\frac{i}{2}}(a^*))$ belongs to $D(\nab^\frac{1}{2})$ and
$$\nab^\frac{1}{2} \la(\si_{-\frac{i}{2}}(a^*))
= \la(\si_{-i}(a^*))\ .$$

Similarly, we have that $a^* x$ belongs to $N \cap D(\si_{-\frac{i}{2}})$
and $\si_{-\frac{i}{2}}(a^* x)$ belongs to $N$. Therefore, $\la(a^* x)$
belongs to $D(\nab^\frac{1}{2})$ and
$$\nab^\frac{1}{2} \la(a^* x) = \la(\si_{-\frac{i}{2}}(a^* x)) \ .$$

Using these two results, equation (*) implies that
\begin{eqnarray*}
& & \langle \pi(x) \la(a) , \la(a) \rangle
=  \langle \la(a^* x) , \nab^\frac{1}{2} \la(\si_{-\frac{i}{2}}(a^*)) \rangle
 =   \langle \nab^\frac{1}{2} \la(a^* x) , \la(\si_\frac{i}{2}(a)^*) \rangle \\
& &  \spat = \langle \la(\si_{-\frac{i}{2}}(a^* x)) , \la(\si_\frac{i}{2}(a)^*)
\rangle
= \langle \la(\si_{-\frac{i}{2}}(a^*) \, \si_{-\frac{i}{2}}(x)) ,
\la(\si_\frac{i}{2}(a)^*) \rangle \\
& & \spat =\langle I \pi(\si_\frac{i}{2}(\si_{-\frac{i}{2}}(x)))^* I^*
\la(\si_{-\frac{i}{2}}(a^*)) , \la(\si_\frac{i}{2}(a)^* \rangle \\
& & \spat =  \langle I \pi(x)^* I^* \la(\si_\frac{i}{2}(a)^*) ,
\la(\si_\frac{i}{2}(a)^*) \rangle.
\end{eqnarray*}
Because $C$ is dense in $A$ and $\pi$ is non-degenerate, this equality implies easily that
we can replace $x$ in this equality by 1. Hence,
$$ \langle \la(a) , \la(a) \rangle
= \langle \la(\si_\frac{i}{2}(a)^*) , \la(\si_\frac{i}{2}(a)^*) \rangle \ . $$
The conclusion follows from proposition \ref{we3.prop2}.
\end{demo}

\bigskip\bigskip

\subsection{The last one.}

In the last part of this section, we describe a natural way to get a KMS-weight on a
\cst-algebra from a left Hilbert algebra.
\label{pag1}

\medskip

Consider a Hilbert space $H$ and a left Hilbert algebra $\cu$ on $H$ with modular operator
$\nab$, modular conjugation $J$ and associated faithful semi-finite normal weight
$\tilde{\vfi}$ on ${\cal L}(\cu)$. Hence $$\cN_{\tilde{\vfi}} = \{\, L_v \mid v \in H
\text{ such that } v \text{ is left bounded with respect to } \cu \, \}$$ and
$\tilde{\vfi}((L_w)^* (L_v)) = \langle v , w \rangle$ for every $v,w$ in $H$ which are
left bounded with respect to $\cu$. We have also that
$$\cN_{\tilde{\vfi}} \cap \cN_{\tilde{\vfi}}^* = \{\,L_v \mid v \in \cu'' \,\} \ . $$

\medskip

Let $A$ be a \cst-algebra and $\pi$ a nondegenerate $^*$-representation of $A$ on $H$ such that $\pi(A) \subseteq
{\cal L}(\cu)$ and suppose there exists a one-parameter group $\si$ on  $A$ such that
$\pi(\si_t(a)) = \nab^{it} \pi(a) \nab^{-it}$ for every $t \in \R$ and $a \in A$.

By definition, $\tilde{\vfi} \, \pi$ will denote the mapping $\tilde{\vfi}
\circ (\pi\restriction_{A^+})$ from $A^+$ into $[0,\infty]$. We put
$\vfi = \tilde{\vfi} \, \pi$.

Then $\vfi$ is a lower semi-continuous weight on $A$ such that
\begin{itemize}
\item $\Mfi^+ = \{\, a \in A \mid \pi(a) \text{ belongs to }
{\cal M}_{\tilde{\vfi}}^+ \,\}$
\item $\Mfi \subseteq \{\, a \in A \mid \pi(a) \text{ belongs to }
{\cal M}_{\tilde{\vfi}} \,\}$
\item  $\Nfi = \{\, a \in A \mid \pi(a) \text{ belongs to } {\cal N}_{\tilde{\vfi}} \,\}$
\item $\Nfi \cap \Nfi^* = \{\, a \in A \mid \pi(a) \text{ belongs to }
{\cal N}_{\tilde{\vfi}} \cap \cN_{\tilde{\vfi}}^* \,\}$
\end{itemize}
and $\vfi(a) = \tilde{\vfi}(\pi(a))$ for every $a \in \Mfi$.

\medskip

Define the following mapping $\la$ from $\Nfi$ into $H$ : \ For every $a \in \Nfi$, there
exists a unique element $v \in H$ such that $v$ is left bounded with respect to $\cu$ and
$L_v = \pi(a)$ and we define $\la(a) = v$, so we get that $L_{\la(a)} = \pi(a)$.

In the following, we will suppose that $\Nfi$ is dense in $A$ and
that $\la(\Nfi)$ is dense in $H$.

It is not difficult to check in this case that  $(H,\la,\pi)$ is a GNS-construction for
$\vfi$.

It follows easily that $\la(\Nfi \cap \Nfi^*) \subseteq \cu''$ and
$\la(a)^* = \la(a^*)$ for every $a \in \Nfi \cap \Nfi^*$.

\begin{proposition}
We have that $\vfi$ is a KMS-weight with modular group $\si$
and such that $\la(\si_t(a)) = \nab^{it} \la(a)$ for every $a \in \Nfi$.
\end{proposition}
\begin{demo}
\begin{enumerate}
\item Choose $a \in \Nfi$ and $t \in \R$. Then $\la(a)$ is left bounded with
respect to $\cu$ and $L_{\la(a)} = \pi(a)$. So $\pi(\si_t(a)) = \nab^{it} \pi(a)
\nab^{-it} = \nab^{it} L_{\la(a)} \nab^{-it}$. Hence, the Tomita theorem implies that
$\nab^{it} \la(a)$ is left bounded with respect to $\cu$ and $L_{\nab^{it} \la(a)} =
\nab^{it} L_{\la(a)} \nab^{-it} = \pi(\si_t(a))$. This implies that $\si_t(a)$ belongs to
$\Nfi$ and $\la(\si_t(a)) = \nab^{it} \la(a)$.

From this, we get immediately that $\vfi$ is invariant under $\si$.
\item Choose $b \in D(\si_{-\frac{i}{2}}) \cap \Nfi \cap \Nfi^*$
and $y \in \Nfi^* \cap \Nfi$. Then $y^* b$ belongs to $\Nfi \cap \Nfi^*$.

\medskip

The function $S(-\frac{i}{2}) \rightarrow B(H) : u \mapsto \pi(\si_u(b))$ is continuous on
$S(-\frac{i}{2})$, analytic on $S(-\frac{i}{2})^0$ and $\pi(\si_t(b)) = \nab^{it} \pi(b)
\nab^{-it}$ for every $t \in \R$ (by assumption). This implies that $D(\nab^\frac{1}{2}
\pi(b) \nab^{-\frac{1}{2}}) = D(\nab^{-\frac{1}{2}}) $ and $\nab^\frac{1}{2} \pi(b)
\nab^{-\frac{1}{2}} \subseteq \pi(\si_{-\frac{i}{2}}(b))$. \ \ \ (a)

Because $y$ belongs to $\Nfi \cap \Nfi^*$ ,\ $\la(y)$ belongs to $\cu''$ and $\la(y)^* =
\la(y^*)$. Therefore we have that $J \la(y)$ belongs to $D(\nab^{-\frac{1}{2}})$ and
$\nab^{-\frac{1}{2}} J \la(y) = \la(y^*)$. Using (a), we have also that $\pi(b)
\nab^{-\frac{1}{2}} J \la(y)$ belongs to $D(\nab^\frac{1}{2})$ and
$\nab^\frac{1}{2}(\pi(b) \nab^{-\frac{1}{2}} J \la(y)) = \pi(\si_{-\frac{i}{2}}(b)) J
\la(y)$.

These two results imply together that $\la(b y^*)$ belongs to $D(\nab^\frac{1}{2})$ and
$\nab^\frac{1}{2} \la(b y^*) = \pi(\si_{-\frac{i}{2}}(b)) J \la(y)$

Because $b y^*$ belongs to $\Nfi \cap \Nfi^*$, we have that $J \nab^\frac{1}{2} \la(b y^*)
= \la( y b^*)$. Substituting this in the last equality gives us that
$$\la(y b^*) = J \pi(\si_{-\frac{i}{2}}(b)) J \la(y) \text{\ \ \ \ \ \ \ \ \ \ \ \ \ \ (b)}$$

\medskip

Because $\Nfi \cap \Nfi^*$ is dense in $A$, proposition \ref{we2.prop2} implies  that
$D(\si_{-\frac{i}{2}}) \cap \Nfi \cap \Nfi^*$ is a core for $\si_{-\frac{i}{2}}$. We know
also that $\Nfi \cap \Nfi^*$ is a core for $\la$. Using these two facts and equality (b),
it is now rather easy to prove for every $a \in D(\si_{-\frac{i}{2}})$ and $x \in \Nfi$
that $x a^*$ belongs to $\Nfi$ and $\la(x a^*) = J \pi(\si_{-\frac{i}{2}}(a)) J \la(x)$.

Proposition \ref{we3.prop3} implies that $\vfi$ is KMS with respect to $\si$.
\end{enumerate}
\end{demo}

\begin{proposition} \label{we3.prop4}
We have that $\la(\Nfi \cap \Nfi^*)$ is a sub left Hilbert algebra of
$\cu''$ and $\la(\Nfi \cap \Nfi^*)'' = \cu''$.
\end{proposition}
\begin{demo}
Because $\vfi$ is a densely defined lower semi-continuous weight on $A$, we know that
$\la(\Nfi \cap \Nfi^*)$ is dense in $H$. By a result of \cite{Stra},  we see that
$\la(\Nfi \cap \Nfi^*)$ is a sub left Hilbert algebra of $\cu''$.

Choose $n \in \N$ and $a \in \Nfi \cap \Nfi^*$.
Define $$a(n) = \frac{n}{\sqrt{\pi}} \int \exp(-n^2 t^2) \, \si_t(a) \, dt \ . $$
Lemma \ref{we2.lem1} implies that $a(n)$ belongs to $\Nfi \cap \Nfi^*$ and
$$\la(a(n)) = \frac{n}{\sqrt{\pi}} \int \exp(-n^2 t^2) \, \nab^{it} \la(a) \, dt \ . $$

Because $\la(\Nfi \cap \Nfi^*)$ is dense in $H$, this equation implies that the set
$\langle \, \la(a(n)) \mid n \in \N, a \in \Nfi \cap \Nfi^* \,\rangle$ is a core for
$\nab^\frac{1}{2}$. Therefore, we get that $\la(\Nfi \cap \Nfi^*)$ is a core for
$\nab^\frac{1}{2}$.

Hence, left Hilbert algebra theory implies that $\la(\Nfi \cap \Nfi^*)'' = \cu''$.
\end{demo}

\section{Properties of KMS-weights} \label{we4}

In this section we will prove important properties of KMS-weights but we
will start of with some equivalent definitions of KMS-weights.

\begin{proposition}  \label{we4.prop12}
Consider a $C^*$-algebra $A$ and a densely defined weight $\vfi$ on $A$ with
GNS-construction $(H,\la,\pi)$. Let $\si$ be  a norm continuous one-parameter group on
$A$. Then $\vfi$ is a KMS-weight with respect to $\si$ if and only if
\begin{enumerate}
\item The mapping $\la : \Nfi \mapsto H$ is closed.
\item We have that $\vfi  \, \si_t = \vfi$ for every $t \in \R$.
\item There exists a core $K$ for $\la$ such that
\begin{itemize}
\item $K$ is a subset of $D(\si_\frac{i}{2})$, $\si_\frac{i}{2}(K)^*$ is a subset of
$\Nfi$ and $\|\la(x)\| = \|\la(\si_\frac{i}{2}(x)^*)\|$ for every $x \in K$.
\item We have that $\si_t(K) \subseteq K$ for every $t \in \R$.
\end{itemize}
\end{enumerate}
\end{proposition}

This proposition is a direct result of  definition \ref{we3.def1} and proposition
\ref{we3.prop2}. The next proposition is an easy consequence of definition \ref{we3.def2}
and proposition \ref{we3.prop3}.

\begin{proposition}
Consider a $C^*$-algebra $A$ and a densely defined weight $\vfi$ on $A$ with
GNS-construction $(H,\pi,\la)$. Let $\si$ be a norm continuous one-parameter group on A.
Then $\vfi$ is a KMS-weight with respect to $\si$ if and only if
\begin{enumerate}
\item The mapping $\la : \Nfi \mapsto H$ is closed.
\item We have that $\vfi  \, \si_t = \vfi$ for every $t \in \R$.
\item There exists an anti-unitary operator $I$ on $H$ such that the following holds.

There exists a core $K$ for $\la$ and a core $C$ for  $\si_\frac{i}{2}$ such that we have
for every $x \in K$ and $a \in C$ that $x a$ belongs to $\Nfi$ and $\la(x a) = I
\pi(\si_\frac{i}{2}(a))^* I^* \la(x)$.
\end{enumerate}
\end{proposition}

\medskip

Notice that the condition of lower semi-continuity of the weight is not assumed. This is
replaced by the weaker condition of the closedness of the mapping $\la$.

\bigskip\medskip

For the most part of this section, we will fix a \cst-algebra $A$ and a KMS-weight $\vfi$
on $A$ with modular group $\si$. Let $(H,\la,\pi)$ be a GNS-construction for $\vfi$.

\bigskip

Because $\vfi$ is invariant with respect to $\si$, there exists a unique
injective positive operator $\nab$ on $H$ such that
$\nab^{it} \la(a) = \la(\si_t(a))$ for every $a \in \Nfi$.
Then we have that
$\pi(\si_t(a)) = \nab^{it} \pi(a) \nab^{-it}$ for every $t \in \R$ and $a
\in A$.

\medskip

Define the set
$$ C_\vfi =  \{ \, a \in \Nfi \cap \Nfi^* \mid a \text{ is analytic with
respect to } \si \text{ and }  \si_z(a) \text{ belongs to } \Nfi \cap \Nfi^* \text{ for
every } z \in \C \,\} \ .$$
Then  $C_\vfi$ is a sub-$^*$-algebra of $A$. It is also clear
that $\si_z(C)
\subseteq C$ for every $z \in \C$.

Because $\vfi$ is lower semi-continuous and densely defined, we know that $\Nfi \cap
\Nfi^*$ is a core for $\la$. So, proposition \ref{we2.prop2} implies that $C_\vfi$ is a
core for $\Nfi$.

\medskip

Let us also define the set
$$K_\vfi = \{ \, a \in \Nfi \cap D(\si_\frac{i}{2}) \mid \si_\frac{i}{2}(a)^*
\text{ belongs to } \Nfi \,\} \ .$$
Because, $C_\vfi$ is a subset of $K_\vfi$, we have that $K_\vfi$ is also a core for $\la$.
It is also clear that $\|\la(a)\| =  \| \la(\si_\frac{i}{2}(a)^*) \|$ for every $a \in
K_\vfi$. Furthermore, $\si_t(K_\vfi) \subseteq K_\vfi$ for every $t \in \R$.

\begin{definition} \label{we4.def1}
Define the anti-unitary $J$ on $H$ such that $J \la(a) = \la(\si_\frac{i}{2}(a)^*)$ for
every $a \in K_\vfi$.
\end{definition}

We have that $J K_\vfi \subseteq K_\vfi$ and $J(J(v)) = v$ for every $v \in K_\vfi$.
Therefore, $J^2 = 1$.

\medskip

The proof of the following proposition is the same as te proof of proposition
\ref{we3.prop1}.

\begin{proposition} \label{we4.prop1}
Consider $x \in \Nfi$ and $a \in D(\si_\frac{i}{2})$. Then $x a$ belongs to $\Nfi$ and
$$\la(x a) = J \pi(\si_\frac{i}{2}(a))^* J \la(x) \ .$$
\end{proposition}

\medskip

In the same way as in the proof of lemma \ref{we3.lem1}, the previous proposition and
proposition \ref{pa2.prop2} imply the following one.

\begin{proposition}  \label{we4.prop2}
We have that $\vfi$ is regular and has a truncating net consisting of elements of $C_\vfi$.
\end{proposition}

\bigskip

\begin{result}
Let $k$ be a natural number. Then $(C_\vfi)^k$ is a core for $\la$.
\end{result}
\begin{demo}
We know already that $C_\vfi$ is a core for $\la$. So we can suppose that $k \geq 2$.
Choose $x \in C_\vfi$.

Because $(C_\vfi)^{k-1}$ is dense in $A$, there exists a bounded net $(e_k)_{k \in K}$ in
$(C_\vfi)^{k-1}$ such that $(e_k)_{k \in K}$ converges strictly to $1$. Then we have for
every $k \in K$ that $e_k x \in (C_\vfi)^k$. It  is also clear that  $(e_k x)_{k \in K}$
converges to $x$.

Because $\la(e_k x) = \pi(e_k) \la(x)$ for every $k \in K$, we get also that $\bigl(\la(e_k
x)\bigr)_{k \in K}$ converges to $\la(x)$.

Because $C_\vfi$ is a core for $\la$, this implies that $(C_\vfi)^k$ is a core for
$C_\vfi$.
\end{demo}

Combining this result with proposition \ref{we1.prop2}, this implies the following result.

\begin{result}  \label{we4.res1}
Let $k$ be a natural number. Consider $a \in \Nfi \cap \Nfi^*$. Then there exists a
sequence $(a_n)_{n=1}^\infty$ in $(C_\vfi)^k$ such that
\begin{enumerate}
\item $(a_n)_{n=1}^\infty \rightarrow a$
\item $\bigl(\la(a_n)\bigr)_{n=1}^\infty \rightarrow \la(a)$
\item $\bigl(\la(a_n^*)\bigr)_{n=1}^\infty \rightarrow \la(a^*)$.
\end{enumerate}
\end{result}

\bigskip

In the rest of this section, we will prove some familiar results

\medskip

\begin{proposition} \label{we4.prop3}
Consider $a \in D(\si_{-i})$ and $x \in \Mfi$. Then $a x$ and $x \si_{-i}(a)$ belong to
$\Mfi$ and $\vfi(a x) = \vfi(x \si_{-i}(a))$.
\end{proposition}
\begin{demo}
Choose $y,z \in \Nfi$. Because $a^*$ belongs to $D(\si_\frac{i}{2})$, we know that $y a^*$
belongs to $\Nfi$ and
$$ \la(y a^*) = J \pi(\si_\frac{i}{2}(a^*))^* J \la(y) = J \pi(\si_{-\frac{i}{2}}(a)) J \la(y)
\ .$$
Because $a y^* z = (y a^*)^* z$, we get that $a y^* z$ belongs to $\Mfi$ and
\begin{eqnarray*}
& & \vfi(a y^* z)  =   \vfi( (y a^*)^* z ) = \langle \la(z) ,
\la(y a^*) \rangle \\
& & \spat = \langle \la(z) , J \pi(\si_{-\frac{i}{2}}(a)) J \la(y) \rangle
 =  \langle J \pi(\si_{-\frac{i}{2}}(a))^* J \la(z) , \la(y) \rangle
\text{\ \ \ \ \ \ (*)}
\end{eqnarray*}
We have that $\si_{-i}(a)$ belongs to $D(\si_\frac{i}{2})$ and
$\si_\frac{i}{2}(\si_{-i}(a)) = \si_{-\frac{i}{2}}(a)$.
This implies that $z \si_{-i}(a)$ belongs to $\Nfi$ and
$$ \la(z \si_{-i}(a)) = J \pi\bigl(\si_\frac{i}{2}(\si_{-i}(a))\bigr)^* J \la(z) = J
\pi(\si_{-\frac{i}{2}}(a))^* J \la(z) \ . $$
Therefore, $y^* z \si_{-i}(a)$ belongs to $\Mfi$
and
$$ \vfi(y^* z \si_{-i}(a)) = \langle \la( z \si_{-i}(a)) , \la(y) \rangle =
\langle J \pi(\si_{-\frac{i}{2}}(a))^* J \la(z) , \la(y) \rangle \ .$$
Comparing this with equation (*), we see that $\vfi(a y^* z) = \vfi(y^* z \si_{-i}(a))$.
\end{demo}

\begin{corollary}  \label{we4.cor1}
Consider $a \in D(\si_i)$ and $x \in \Mfi$. Then $x a$ and $\si_i(a) x$ belong to $\Mfi$
and $\vfi(x a) = \vfi(\si_i(a) x)$.
\end{corollary}

Now, we prove a generalization of proposition \ref{we4.prop1}.

\begin{proposition} \label{we4.prop5}
Consider $x \in \Nfi$ and $a \in D(\overline{\si}_\frac{i}{2})$. Then $x a$ belongs to
$\Nfi$ and $$\la(x a) = J \pi(\si_\frac{i}{2}(a))^* J \la(x) \ .$$
\end{proposition}
\begin{demo}
Choose $y \in \Nfi$ and $b \in D(\si_\frac{i}{2})$. By  proposition \ref{we4.prop1}, we
know that  $y b$ belongs to $\Nfi$ and $\la(y b) = J \pi(\si_\frac{i}{2}(b))^* J \la(y)$.

We have also that  $b a$ belongs to $D(\si_\frac{i}{2})$ and $\si_\frac{i}{2}(b a) =
\si_\frac{i}{2}(b) \, \si_\frac{i}{2}(a)$. Using  proposition \ref{we4.prop1} once more,
we see that $(y b) a$ belongs to $\Nfi$ and
$$
\la((y b) a) =  J \pi(\si_\frac{i}{2}(b a))^* J \la(y)
= J \pi(\si_\frac{i}{2}(a))^* J \, J \pi(\si_\frac{i}{2}(b))^* J \la(y)
= J \pi(\si_\frac{i}{2}(a))^* J  \la(y b) \ .
$$

Because $(C_\vfi)^2$ is a core for $\la$, we have  that $\langle \, y b  \mid y \in \Nfi,
\, b \in D(\si_\frac{i}{2}) \, \rangle$ is a core for $\la$. Using the results of the
previous discussion, this implies easily that $x a$ belongs to $\Nfi$ and $\la(x a) = J
\pi(\si_\frac{i}{2}(a))^* J \la(x)$.
\end{demo}

Therefore, we can prove the next proposition in the same way as proposition \ref{we4.prop3}.

\begin{proposition}
Consider $a \in D(\overline{\si}_{-i})$ and $x \in \Mfi$. Then $a x$ and $x \si_{-i}(a)$
belong to $\Mfi$ and $\vfi(a x) = \vfi(x \si_{-i}(a))$.
\end{proposition}

\begin{corollary}
Consider $a \in D(\overline{\si}_i)$ and $x \in \Mfi$. Then $x a$ and $\si_i(a) x$ belong
to $\Mfi$ and $\vfi(x a) = \vfi(\si_i(a) x)$.
\end{corollary}

\medskip

Another consequence of proposition \ref{we4.prop1} can be found in the following
proposition.

\begin{result} \label{we4.res2}
Consider $a,b \in A$. Then we have that $J \pi(a) J$
and $\pi(b)$ commute.
\end{result}
\begin{demo}
Choose $c \in D(\si_\frac{i}{2})$.

Take $x \in \Nfi$. Then $x c$ belongs to $\Nfi$ and $\la(x c) = J
\pi(\si_\frac{i}{2}(c))^* J \la(x)$. So $b x c$ belongs to $\Nfi$ and $$\la(b x c) =
\pi(b) \la(x c) = \pi(b) \, J \pi(\si_\frac{i}{2}(c))^* J \la(x) \ .$$

On the other hand, we have that $b x$ belongs to $\Nfi$ and
$\la(b x) = \pi(b) \la(x)$. This implies that $b x c$ belongs
to $\Nfi$ and $$\la(b x c) = J \pi(\si_\frac{i}{2}(c))^* J  \la(b x)
= J \pi(\si_\frac{i}{2}(c))^* J \, \pi(b) \la(x) \ .$$

Comparing these two expressions for $\la(b x c)$, we see that
$$\pi(b) \, J \pi(\si_\frac{i}{2}(c))^* J \la(x) =
J \pi(\si_\frac{i}{2}(c))^* J \, \pi(b) \la(x) \ . $$
This implies that $\pi(b) \, J \pi(\si_\frac{i}{2}(c))^* J
= J \pi(\si_\frac{i}{2}(c))^* J \, \pi(b)$.

Because $\si_\frac{i}{2}$ has dense range in $A$, the result follows.
\end{demo}

\bigskip

\begin{proposition} \label{we4.prop4}
The set $\la(\Nfi \cap \Nfi^*)$ is a core for $\nab^\frac{1}{2}$ and
$\la(a^*) = J \nab^\frac{1}{2} \la(a)$ for every $a \in \Nfi \cap \Nfi^*$.
\end{proposition}
\begin{demo}
Choose $a \in \Nfi \cap \Nfi^*$. We know that there exists a sequence $a_n \in
C_\vfi$ such that $\bigl(\la(a_n)\bigr)_{n=1}^\infty \rightarrow \la(a)$
and $\bigl(\la(a_n^*)\bigr)_{n=1}^\infty \rightarrow \la(a^*)$ (see result \ref{we4.res1}).

Choose $m \in \N$.
Because $a_m$ belongs to $\Nfi \cap D(\si_{-\frac{i}{2}})$ and
$\si_{-\frac{i}{2}}(a_m)$ belongs to $\Nfi$, we have that $\la(a_m)$ belongs to
$D(\nab^\frac{1}{2})$ and $\nab^\frac{1}{2} \la(a_m)
= \la(\si_{-\frac{i}{2}}(a_m))$.

We also know that $a_m^*$ belongs to $\Nfi \cap D(\si_\frac{i}{2})$
and $\si_\frac{i}{2}(a_m)^*$ belongs to $\Nfi^*$. By definition, we have that
$$ J \la(a_m^*) = \la(\si_\frac{i}{2}(a_m^*)^*) =
\la(\si_{-\frac{i}{2}}(a_m)) =  \nab^\frac{1}{2} \la(a_m) \ .$$

This implies that $\bigl(\nab^\frac{1}{2} \la(a_n)\bigr)_{n=1}^\infty $ converges
to $J \la(a^*)$. The closedness of $\nab^\frac{1}{2}$ implies that $\la(a)$
belongs to $D(\nab^\frac{1}{2})$ and $\nab^\frac{1}{2} \la(a) = J \la(a^*)$.

\medskip

Choose $x \in \Nfi \cap \Nfi^*$, \ $n \in \N$ and define
$$ x(n) = \frac{n}{\sqrt{\pi}} \int \exp(- n^2 t^2) \, \si_t(x) \, dt \ ,$$
it is clear that $x(n)$ belongs to $\Nfi \cap \Nfi^*$ and
$$\la(x(n)) = \frac{n}{\sqrt{\pi}} \int \exp(-n^2 t^2) \, \nab^{it} \la(x) \, dt \ .$$
Because $\la(\Nfi \cap \Nfi^*)$ is dense in $H$, this implies that
$\langle \, \la(x(n)) \mid x \in \Nfi \cap \Nfi^* , n \in  \N  \, \rangle$
is a core for $\nab^\frac{1}{2}$. This implies that $\la(\Nfi \cap \Nfi^*)$ is
a core for $\nab^\frac{1}{2}$.
\end{demo}

\medskip

We want to construct a left Hilbert algebra in the usual way. First, we will need a lemma
which guarantees that our next definitions are well defined.

\begin{lemma}
\begin{enumerate}
\item Consider $a_1,a_2,b_1,b_2 \in \Nfi \cap \Nfi^*$ such that $\la(a_1) = \la(a_2)$ and
$\la(b_1) = \la(b_2)$. Then $\la(a_1 b_1) = \la(a_2 b_2)$.
\item Consider $a_1,a_2 \in \Nfi \cap \Nfi^*$ such that $\la(a_1) =
\la(a_2)$. Then $\la(a_1^*) =  \la(a_2^*)$.
\end{enumerate}
\end{lemma}
\begin{demo}
\begin{enumerate}
\item
For any $c \in C_\vfi$, we have that $$ \pi(a_1) \la(c) =
J \pi(\si_\frac{i}{2}(c))^* J \la(a_1)
= J \pi(\si_\frac{i}{2}(c))^* J \la(a_2) = \pi(a_2) \la(c) \ .$$
The density of $\la(C_\vfi)$ in $H$  implies that $\pi(a_1) = \pi(a_2)$.

Therefore,
$$\la(a_1 b_1) = \pi(a_1) \la(b_1) = \pi(a_2) \la(b_2) = \la(a_2
b_2) \ .$$
\item
This follows easily from the previous proposition.
\end{enumerate}
\end{demo}

\begin{definition}
We define ${\cal U} = \la(\Nfi \cap \Nfi^*)$, then ${\cal U}$ is subspace
of $H$. We make ${\cal U}$ into a $^*$-algebra such that :
\begin{enumerate}
\item We have for every $a,b \in \Nfi \cap \Nfi^*$ that $\la(a b) = \la(a)
\la(b)$.
\item We have for every $a \in \Nfi \cap \Nfi^*$ that $\la(a)^* = \la(a^*)$.
\end{enumerate}
Then ${\cal U}$ becomes a left Hilbert algebra on $H$.
\end{definition}

It is not difficult to check that ${\cal U}$ satisfies the conditions
of a left Hilbert algebra. The closability of the mapping ${\cal U} \rightarrow
{\cal U} : v \mapsto v^*$ is a direct consequence of proposition \ref{we4.prop4}.

It is immediately clear that $L_{\la(a)} = \pi(a)$ for every $a \in \Nfi
\cap \Nfi^*$.

Let us define $T$ as the closure of the mapping
${\cal U} \rightarrow {\cal U} : v \mapsto v^*$.
Proposition \ref{we4.prop4} implies that $T =  J \nab^\frac{1}{2}$.

So we see that $J$ is the modular conjugation of ${\cal U}$ and $\nab$ is the
modular operator of ${\cal U}$.

We also see that $J$ and $\nab$ are independent of our choice of $\si$.
Therefore, we will call $\nab$ the modular operator of $\vfi$
and $J$ the modular conjugation of $\vfi$ in the GNS-construction $(H,\la,\pi)$.

\begin{lemma}
Consider $a \in \Nfi$. Then $\la(a)$ is left bounded with respect to
${\cal U}$ and $L_{\la(a)} = \pi(a)$.
\end{lemma}
\begin{demo}
Because $\Nfi \cap \Nfi^*$ is a core for $\la$, there exists a sequence
$(a_n)_{n=1}^\infty$ in $\Nfi \cap \Nfi^*$ such that $(a_n)_{n=1}^\infty$
converges to $a$ and $\bigl(\la(a_n) \bigr)_{n=1}^\infty$ converges to $\la(a)$.

We know already that $L_{\la(a_n)} = \pi(a_n)$ for $n \in \N$. This implies that
$\bigl(L_{\la(a_n)}\bigr)_{n=1}^\infty$ converges to $\pi(a)$. Consequently, we find
that $\la(a)$ is left bounded with respect to ${\cal U}$ and $L_{\la(a)} = \pi(a)$.
\end{demo}

\begin{lemma}
Consider $v \in H$ such that $v$ is left bounded with respect to ${\cal
U}$ and such that there exists an element $a \in A$ such that $L_v = \pi(a)$.
Then $a$ belongs to $\Nfi$ and $\la(a)=v$.
\end{lemma}
\begin{demo}
Choose a truncating net $(u_i)_{i \in I}$  for $\vfi$. For every $i \in I$, we define the
operator $S_i \in B(H)$ such that $S_i \la(x) = \la( x u_i)$ for every $x \in A$.

By \cite{Uffe}, there exist a sequence $(a_n)_{n=1}^\infty$ in $\Nfi \cap \Nfi^*$ such
that $\bigl(\la(a_n)\bigr)_{n=1}^\infty$ converges to $\la(a)$,
$\bigl(L_{\la(a_n)}\bigr)_{n=1}^\infty$
converges strongly to $L_v$ and $\bigl(L_{\la(a_n)}\bigr)_{n=1}^\infty$ is bounded.

Choose $j \in I$.
We have immediately that $\bigl(L_{\la(a_n)} \la(u_j) \bigr)_{n=1}^\infty$
converges to   $L_v \la(u_j)$. We have for every
$n \in \N$ that
$$ L_{\la(a_n)} \la(u_j)  =  \pi(a_n) \la(u_j) = \la(a_n u_j) = S_j
\la(a_n) \ .$$
This implies that $\bigl(L_{\la(a_n)} \la(u_j) \bigr)_{n=1}^\infty$
converges to $S_j v$.

Combining these two results, we get that $L_v \la(u_j) = S_j v$, hence
$$ \la(a u_j) = \pi(a) \la(u_j) = L_v \la(u_j) = S_j v \ .$$

Therefore, $\bigl(\la(a u_i)\bigr)_{i \in I}$ converges to $v$. We also have that $(a u_i)_{i \in
I}$ converges to $a$. The closedness of $\la$ implies that $a$ belongs to $\Nfi$ and
$\la(a) = v$.
\end{demo}

\begin{lemma}
Consider $v \in {\cal U}''$ such that there exist an element $a \in A$ such that
$L_v = \pi(a)$. Then $a$ belongs to $\Nfi \cap \Nfi^*$ and $v= \la(a)$.
\end{lemma}
\begin{demo}
We know that $v$ is left bounded. Therefore, the previous lemma implies that $a$ belongs
to $\Nfi$ and $v=\la(a)$. We have also that $v^*$ is left bounded and $L_{v^*} = (L_v)^*
= \pi(a)^* = \pi(a^*)$. Hence, the previous lemma implies that $a^*$ belongs to $\Nfi$. So
we get that $a$ belongs to $\Nfi \cap \Nfi^*$.
\end{demo}

Using this result, we have the following theorem. This is a well known result for weights
satisfying the other KMS-condition (proved by Combes, see lemma 2.2.3 of \cite{E-V}).

\begin{theorem}
Call $\tilde{\vfi}$ the normal semi-finite faithful weight on $\pi(A)''$ associated to
${\cal U}$.  Then $\tilde{\vfi}\,\pi = \vfi$.
\end{theorem}

\begin{remark}\rm
Consider the case where we get a KMS-weight from a left Hilbert algebra (the third part of
section \ref{we3}, page \pageref{pag1}) :

Consider  a Hilbert space $H$ and a left Hilbert ${\cal W}$ algebra on $H$ with modular
operator $\de$, modular conjugation $I$ and associated faithful semi-finite normal weight
$\psi$.

Let $A$ be a \cst-algebra and $\pi$ be a nondegenerate $^*$-representation of $A$ on $H$
such that $\pi(A) \subseteq {\cal L}({\cal W})$ and suppose there exists a one-parameter
group $\si$ on  $A$ such that $\pi(\si_t(a)) = \de^{it} \pi(a) \de^{-it}$ for every $t \in
\R$ and $a \in A$.

Put $\vfi = \psi \, \pi$, so $\vfi$ is a lower semi-continuous weight on $A$.

Define the following mapping $\la$ from $\Nfi$ into $H$ : \ For every $a \in \Nfi$, there
exists a unique element $v \in H$ such that $v$ is left bounded with respect to ${\cal W}$
and $L_v = \pi(a)$ and we define $\la(a) = v$, so we get that $L_{\la(a)} = \pi(a)$. We
suppose that $\Nfi$ is dense in $A$ and that $\la(\Nfi)$ is dense in $H$.

We saw that $\vfi$ is a KMS-weight on $A$ with modular group $\si$ and GNS-construction
$(H,\la,\pi)$. Now we can apply the construction procedures of this section to get the
left Hilbert algebra $\cu$, the weight $\tilde{\vfi}$ associated to $\cu$, the modular
operator $\nab$ and modular conjugation $J$ of $\vfi$ in the GNS-construction
$(H,\la,\pi)$.

Proposition \ref{we3.prop4} implies that $\cu'' = {\cal W}''$. Therefore we get that
$\tilde{\vfi} = \psi$, $\nab = \de$ and $J = I$.
\end{remark}

\medskip

For KMS-weights, the faithfulness of $\vfi$ and $\pi$ are equivalent :

\begin{proposition} \label{we4.prop7}
We have that $\vfi$ is faithful if and only if $\pi$ is faithful.
\end{proposition}
\begin{demo}
\begin{trivlist}
\item[$\,\,\,\Rightarrow$] This follows easily because $\vfi$ is densely defined.
\item[$\,\,\,\Leftarrow$]
Choose $x \in A$ such that $\vfi(x^* x) = 0$. Then $x$ belongs to $\Nfi$ and $\la(x)=0$.

Take $m \in \N$ and define
$$x(m) = \frac{m}{\sqrt{\pi}} \int \exp(-m^2 t^2) \, \si_t(x) \, dt  \ .$$
We have that $x(m)$ belongs to $\Nfi \cap D(\si_\frac{i}{2})$ and
$$\la(x(m)) = \frac{m}{\sqrt{\pi}} \int \exp(-m^2 t^2) \, \nab^{it} \la(x) \, dt = 0 \  .$$
So we have for every $a \in \Nfi$ that
$$ J \pi\bigl(\si_\frac{i}{2}(x(m))\bigr)^* J \la(a) = \pi(a) \la(x(m)) = 0 \ ,$$
which implies that $\pi\bigl(\si_\frac{i}{2}(x(m))\bigr) = 0$. The faithfulness of $\pi$ implies
that   $\si_\frac{i}{2}(x(m)) = 0$, the injectivity of $\si_\frac{i}{2}$ implies that
$x(m)=0$.

The definition of the elements $x(n) \,\, (n \in \N)$ implies
that $(x(n))_{n=1}^\infty \rightarrow x$. Therefore we have that $x=0$.
\end{trivlist}
\end{demo}

\medskip

We have of course also the following result.

\begin{lemma}
Consider a $^*$-automorphism $\al$ on $A$ such that there exists a strictly positive
number $\lambda$ such that $\vfi \, \al = \lambda \, \vfi$. Define $u$ as the unitary
operator on $H$ such that $u \la(a) = \lambda^{-\frac{1}{2}} \, \la(\al(a))$ for every $a
\in \Nfi$. Then $u \nab = \nab u$ and $u J = J u$. We have moreover that
$\pi(\al(a)) = u \pi(a) u^*$ for every $a \in A$.
\end{lemma}
\begin{demo}
Because $\al$ is a $^*$-homomorphism, it is easy to check for every $a \in \Nfi \cap
\Nfi^*$ that $u \la(a)$ belongs to $D(T)$ and $T u \la(a) = u T \la(a)$. Because $\la(\Nfi
\cap \Nfi^*)$ is a core fore $T$ and $u \la(\Nfi \cap \Nfi^*) = \la(\Nfi \cap \Nfi^*)$ we
get that $u T = T u$. So we get also that $T u^* = u^* T$. The unitarity of $u$ implies
that $$u T^* = (T u^*)^* = (u^* T)^* = T^* u \ . $$ Using the fact that $\nab = T^* T$,
this implies that $u \nab = \nab u$.

So we get also that $u \nab^\frac{1}{2} = \nab^\frac{1}{2} u$. This in turn implies that
$$ u J \nab^\frac{1}{2} = u T = T u = J \nab^\frac{1}{2} u = J u \nab^\frac{1}{2} \ . $$
Because $\nab^\frac{1}{2}$ has dense range, this implies immediately that $u J = J u$.
\end{demo}

\medskip

\begin{corollary} \label{we4.prop10}
Suppose that $\vfi$ is faithful. Consider a $^*$-automorphism $\al$ on $A$ such that there
exists a strictly positive number $\lambda$ such that $\vfi \, \al = \lambda \, \vfi$.
Then $\al \si_t = \si_t \al$ for every $t \in \R$.
\end{corollary}
\begin{demo} Choose $t \in \R$.

Define $u$ as the unitary
operator on $H$ such that $u \la(a) = \lambda^{-\frac{1}{2}} \, \la(\al(a))$ for every $a
\in \Nfi$. Then the previous lemma implies that $u \nab = \nab u$.

This implies immediately that $u \nab^{it} = \nab^{it} u$. So we get for every $a \in
\Nfi$ that
$$
\la\bigl(\al(\si_t(a))\bigr) = \lambda^\frac{1}{2} \, u \la(\si_t(a))
= \lambda^\frac{1}{2} \, u \nab^{it} \la(a)
= \lambda^\frac{1}{2} \, \nab^{it} u \la(a) = \nab^{it} \la(\al(a))
 = \la\bigl(\si_t(\al(a))\bigr)  \ .
$$
so the faithfulness of $\vfi$ implies that $\al(\si_t(a)) = \si_t(\al(a))$.
\end{demo}

\medskip

From the theory of left Hilbert algebras (see e.g. \cite{Stra}), we get the following
result :

\begin{proposition} \label{KMS}
The weight $\vfi$ satisfies the KMS-condition with respect to $\si$ : \ For every $x,y \in
\Nfi \cap \Nfi^*$ there exists a bounded continuous function $f$ from the strip $S(i)$
into $\,\C$ which is analytic on $S(i)^0$ and satisfies :
\begin{itemize}
\item $f(t) = \vfi(\si_t(x)\, y)$ for all $t \in \R$
\item $f(t+i) = \vfi(x \,  \si_t(y))$ for all $t \in \R$.
\end{itemize}
\end{proposition}

\medskip

The next two results follow immediately from proposition \ref{we4.prop5} and its corollary.

\begin{proposition} \label{we4.prop6}
Consider $a \in M(A)$ such that there exists a strictly positive number $\lambda$
such that $\si_t(a) = \lambda^{it} \, a$ for every $t \in A$.
Let $x$ be an element in $\Nfi$. Then $x a$ belongs to $\Nfi$ and
$\la(x a) = \lambda^{-\frac{1}{2}} \, J \pi(a)^* J \la(x)$.
\end{proposition}

\begin{proposition}
Consider $a \in M(A)$ such that there exists a strictly positive number
$\lambda$ such that $\si_t(a) = \lambda^{it} \, a$ for every $t \in \R$.
Let $x$ be an element in $\Mfi$. Then $a x$ and $x a$ belong to $\Mfi$
and $\vfi(a x) = \lambda \, \vfi(x a)$.
\end{proposition}

If $\vfi$ is faithful, it is possible to prove a converse of this corollary.

\medskip

The proof of the next lemma can be copied immediately from the proof of lemma 3.4 of
\cite{Pe-Tak}.

\begin{lemma}
Consider $a \in M(A)$ such that $a \Mfi \subseteq \Mfi$ and $\Mfi a \subseteq \Mfi$. Then
$\si_t(a) \Mfi \subseteq \Mfi$ and $\Mfi \, \si_t(a) \subseteq \Mfi$ for every $t \in
\R$.  Let $x,y$ be elements in $C_\vfi$, then there exists an entire function from $\C$
into $\C$ such that
\begin{enumerate}
\item We have that $f(t) = \vfi(\si_t(a) \, y^* x)$ for every $t \in \R$.
\item We have that $f(t+i) = \vfi(y^* x \, \si_t(a))$ for every $t \in \R$.
\item $f$ is bounded on each horizontal strip.
\end{enumerate}
\end{lemma}

The proof of the following proposition is a slight adaptation of the proof of theorem 3.6
of \cite{Pe-Tak}.

\begin{result} \label{we4.res3}
Suppose that $\vfi$ is faithful. Consider $a \in M(A)$ such that $a \Mfi \subseteq \Mfi$,
$\Mfi \, a \subseteq \Mfi$ and such that there exist a strictly positive number $\lambda$
satisfying $\vfi(a x) = \lambda \, \vfi(x a)$ for every $x \in \Mfi$. Then $\si_t(a) =
\lambda^{it} \, a$ for every $t \in \R$.
\end{result}
\begin{demo}
Choose $x,y \in C_\vfi$.

By the previous lemma, there exists an entire function $f$ from $\C$ into $\C$
such that
\begin{enumerate}
\item We have that $f(t) = \vfi(\si_t(a) \, y^* x)$ for every $t \in \R$.
\item We have that $f(t+i) = \vfi(y^* x \, \si_t(a))$ for every $t \in \R$.
\item $f$ is bounded on each horizontal strip.
\end{enumerate}
First, we prove that $f(z+i) = \lambda^{-1} \, f(z)$ for every $z \in \C$.

\begin{list}{}{\setlength{\leftmargin}{.4 cm}}

\item We have for every $t \in \R$ that
\begin{eqnarray*}
& & f(t+i) - \lambda^{-1} \, f(t) =  \vfi(y^* x \, \si_t(a))
- \lambda^{-1} \,  \vfi(\si_t(a) \, y^* x) \\
& & \spat =  \vfi(y^* x \,\si_t(a)) - \lambda^{-1} \, \vfi(a \, \si_{-t}(y^* x)) \\
& & \spat =  \vfi(y^* x \, \si_t(a)) - \lambda^{-1} \,  \lambda \, \vfi(\si_{-t}(y^* x) \, a) \\
& & \spat = \vfi(y^* x \, \si_t(a)) - \vfi(y^* x \, \si_t(a)) = 0 .
\end{eqnarray*}
By analyticity of the function
$$ \C  \rightarrow \C : z \mapsto f(z+i) - \lambda^{-1} \, f(z) \ , $$
we must have that $f(z+i) - \lambda^{-1} \, f(z) = 0$ for every $z \in \C$.

\end{list}

Now we define the entire function $g$ from $\C$ into $\C$ such that
$g(z) = \lambda^{-iz} \, f(z)$ for every $z \in \C$.
Then we have for every $z \in \C$ that
$$g(z+i) = \lambda^{-i(z+i)} \, f(z+i) = \lambda^{-iz} \, \lambda \,
\lambda^{-1} \,
f(z) = \lambda^{-iz} \, f(z) = g(z)  \ . \text{\ \ \ \ (*)}$$
We know that $f$ is bounded on $S(i)$. We have also that
$$ |g(z)| = |\lambda^{-iz}| \,\, |f(z)| = \lambda^{\text{\scriptsize Im}\, z} \, |f(z)| \leq
\max\{1,\lambda\} \, |f(z)| \ ,$$
for every $z \in S(i)$, so $g$ is also bounded on $S(i)$.

By (*), we have that $g$ is bounded on $\C$. The theorem of Liouville implies
that $g$ is a  constant function.
This implies that
$$\vfi(a y^* x) = f(0) = g(0) = g(t) = \lambda^{-it} \, f(t) = \lambda^{-it} \,
\vfi(\si_t(a) y^* x) $$
for every $t \in \R$.

\medskip

Choose $s \in \R$.
Because $(a-\lambda^{-is} \, \si_s(a)) \Mfi \subseteq \Mfi$, we have also that
$\Nfi (a^* - \lambda^{is} \, \si_s(a^*)) \subseteq \Nfi$.

Hence, we have for every $x,y \in C_\vfi$ that
$$\langle \la(x) , \la\bigl( y (a^* - \lambda^{is} \, \si_s(a^*))\bigr) \rangle
= \vfi( (a - \lambda^{-is} \, \si_s(a)) y^* x ) = 0 \ .$$
Because $\la(C_\vfi)$ is dense in $H$, we see that $\la\bigl( y (a^* - \lambda^{is} \,
\si_s(a^*))\bigr) = 0$ for every $y \in C_\vfi$. The faithfulness of $\vfi$ implies that $y
(a^* - \lambda^{is} \, \si_s(a^*)) = 0$ for every $y \in C_\vfi$. Therefore, the density
of $C_\vfi$ in $A$ implies that $a^* - \lambda^{is} \, \si_s(a^*) = 0$.
\end{demo}

\medskip

The following proposition is an immediate consequence of proposition \ref{we4.prop6}.

\begin{proposition}  \label{we4.prop9}
Consider a strictly positive element $\sde \, \eta \, A$ such that there exist a strictly
positive number $\lambda$ such that $\si_t(\sde) = \lambda^t \, \sde$ for every $t \in \R$.
Then we have the following properties :
\begin{itemize}
\item
For every $x \in \Nfi$ and every $t \in \R$, we have that $x \sde^{it}$ belongs to $\Nfi$
and $\la(x \sde^{it}) =$ \newline $\lambda^{-\frac{t}{2}} \, J \pi(\sde)^{-it} J \la(x)$.
\item For every $x \in \Mfi$ and every $t \in \R$, we have that $x \sde^{it}$
and $\sde^{it} x$ belong to $\Mfi$ and
$\vfi(\sde^{it} x ) = \lambda^t \, \vfi(x \sde^{it})$.
\end{itemize}
\end{proposition}

Using this proposition, we get the following one.

\begin{proposition} \label{we4.prop8}
Consider a strictly positive element $\sde \, \eta \, A$ such that there exist a strictly
positive number $\lambda$ such that $\si_t(\sde) = \lambda^t \, \sde$ for every $t \in
\R$. Let $z$ be a complex number and $a \in \Nfi$ such that $a$ is a left multiplier of
$\sde^z$ and $a \, \sde^z$ belongs to $A$. Then we have the following properties :
\begin{itemize}
\item $a \, \sde^z$ belongs to $\Nfi$ $\Leftrightarrow$ $\la(a)$
belongs to $D(J \pi(\sde)^{\overline{z}} J)$.
\item If $a \, \sde^z$ belongs to $\Nfi$, then
$\la(a \, \sde^z) = \lambda^\frac{iz}{2} \, J \pi(\sde)^{\overline{z}} J
\la(a)$.
\end{itemize}
\end{proposition}
\begin{demo}
We define the strongly continuous one-parameter representation $\al$ on $A$ such that
$\al_t(x) = x \sde^{it}$ for every $x \in A$ and $t \in \R$.

Furthermore, we define the strongly continuous unitary group representation $u$ from $\R$
on $H$ such that $u_t = (J \pi(\sde) J)^{it}$ for every $t \in
\R$.
We know by the previous proposition that
\begin{enumerate}
\item We have that $\al_t(\Nfi) \subseteq \Nfi$ for every $t \in \R$
\item We have that $u_t \la(a) = \lambda^\frac{t}{2} \, \la(\al_t(a))$ for
every $t \in \R$ and $a \in \Nfi$.
\end{enumerate}
Because $a$ belongs to $D(\al_{-iz})$ and $\al_{-iz}(a) = a \, \sde^z$ and $u_{-iz}
= (J \pi(\sde) J)^z$, this proposition follows from propositions \ref{we2.prop4} and \ref{we2.prop3}.
\end{demo}

\medskip

\begin{proposition} \label{we4.prop11}
Consider a strictly positive element $\sde$ affiliated with $A$ such that there exists a
strictly positive number $\lambda$ such that $\si_t(\sde) = \lambda^t \, \sde$ for every
$t \in \R$. Let $a$ be an element in $A$ and $z \in \C$ such that $a$ is a left and right
muliplier of $\sde^z$ and  $a \, \sde^z$,  $\sde^z a$ belong to $\Mfi$. Then $\vfi\bigl( a
\sde^z
\bigr) =
\lambda^{i z} \, \vfi\bigl( \sde^z a\bigr)$. \end{proposition}
\begin{demo}
Define the norm continuous one-parameter group $\al$ on $A$ such that $\al_t(a) =
\sde^{it} a \sde^{-it}$ for every $a \in A$ and $t \in \R$.
Proposition \ref{we4.prop9} implies that
$\vfi  \, \al_t = \lambda^t \, \vfi$ for every $t \in \R$.

Because $a$ is a left multiplier of $\sde^z$, we have that $a \, \sde^z$ is al left multiplier of
$\sde^{-z}$ and $(a \, \sde^z)\,\sde^{-z} = a$.

This implies that $(a \, \sde^z)\,\sde^{-z}$ is a right multiplier of $\sde^z$ and
$\sde^z\, [(a \, \sde^z)\,\sde^{-z}] = \sde^z a $.
From this, we get that $a \, \sde^z$ is a middle multiplier of $\sde^z, \, \sde^{-z}$ and
$$\sde^z\, (a \, \sde^z) \, \sde^{-z} = \sde^z a \, $$

From this all, we get that $a \, \sde^z$ belongs to $D(\al_{-iz})$ and
$\al_{-iz}(a \, \sde^z) = \sde^z a$. The proposition follows from  proposition
\ref{we1.prop3}.
\end{demo}

\medskip

We now prove that a KMS-weight has no proper extensions which are  invariant under its
modular group.

\begin{proposition} \label{prop2}
Consider a lower semi-continuous weight $\eta$ on $A$ which is an extension of $\vfi$ and
such that $\eta$ is invariant under $\si$. Then $\vfi = \eta$.
\end{proposition}
\begin{demo}
Take a GNS-construction $(H_\eta,\la_\eta,\pi_\eta)$ for $\eta$.

Because $\eta$ is invariant under $\si$, we get the existence a positive injective
operator $T$ in $H_\eta$ such that $T^{it}\la_\eta(a) =  \la_\eta(\si_t(a))$ for every $a
\in \cN_\eta$ and $t \in \R$.

\medskip

Choose $y \in \cN_\eta$.

Define for every $n \in \N$ the element
$$y_n = \frac{n}{\sqrt{\pi}} \int \exp(-n^2 t^2) \, \si_t(y) \, dt $$
which is clearly analytic with respect to $\si$. We have also that
$(y_n)_{n=1}^\infty$ converges to $y$.

\medskip

By ??, we have for every $n \in \N$ that $y_n$ belongs to $\cN_\eta$ and
$$\la_\eta(y_n) = \frac{n}{\sqrt{\pi}} \int \exp(-n^2 t^2) \,
T^{it} \la(y) \, dt  \ .$$ This implies immediately that $(\la_\eta(y_n))_{n=1}^\infty$
converges to $\la_\eta(y)$.

\medskip

We can also take an approximate unit $(e_i)_{i \in I}$ for $A$ in $\Nfi$.
Then  $(e_i \, y_n)_{(i,n) \in I \times \N}$ converges to $y$.

We have also for every $i \in I$ and $n \in \N$ that $e_i \, y_n$ belongs to $\cN_\eta$
and $\la_\eta(e_i \, y_n) = \pi_\eta(e_i) \la_\eta(y_n)$. Consequently, the net
$(\la_\eta(e_i \, y_n))_{(i,n) \in I \times \N}$ converges to $\la_\eta(y)$.

\medskip

We have for every $i \in I$ and $n \in \N$ that $e_i$ belongs to $\Nfi$ and $y_n$ belongs
to $D(\si_\frac{i}{2})$ implying that $e_i \, y_n$ belongs to $\Nfi$ by proposition
\ref{we4.prop1}.

\medskip

Because $\vfi \subseteq \eta$, we have moreover for every $i,j \in I$ and $m,n \in \N$ that
$$\|\lafi(e_i \, y_n) - \lafi(e_j \, y_m)\| = \|\la_\eta(e_i \, y_n) - \la_\eta(e_j \, y_m)\|  \ .$$
This last equality  implies that the net $(\lafi(e_i \, y_n))_{(i,n) \in I \times \N}$ is Cauchy and hence convergent in $H_\vfi$. Therefore, the closedness of $\lafi$ implies that $y$ is an element of $\cN_\vfi$.
The proposition follows.
\end{demo}

\bigskip\bigskip

The following proposition will guarantee that the modular group is unique for faithful
KMS-weights.

\begin{proposition}
Consider a \cst-algebra $A$ and a  KMS-weight $\vfi$ on $A$ with GNS-construction
$(H,\la,\pi)$. Let $\si$ and $\tau$ be modular groups for $\vfi$. Then $\pi \, \si_t = \pi
\, \tau_t$ for every $t \in \R$.
\end{proposition}
\begin{demo}
Define the injective positive operators $\nab$ and  $\de$ in $H$ such that $\nab^{it}
\la(a) = \la(\si_t(a))$ and $\de^{it} \la(a) = \la(\tau_t(a))$ for every $a \in \Nfi$
and $t \in \R$. Define also the anti-unitary operators $J$, $I$ on $H$ such that
\begin{itemize}
\item $J \la(a) = \la(\si_\frac{i}{2}(a)^*)$ for every $a \in \Nfi \cap D(\si_\frac{i}{2})$
such that $\si_\frac{i}{2}(a)^*$ belongs to $\Nfi$
\item $I \la(a) = \la(\tau_\frac{i}{2}(a)^*)$ for every $a \in \Nfi \cap D(\tau_\frac{i}{2})$
such that $\tau_\frac{i}{2}(a)^*$ belongs to $\Nfi$
\end{itemize}
Proposition \ref{we4.prop4} implies that $J \nab^\frac{1}{2} =  I \de^\frac{1}{2}$. The
uniqueness of the polar decomposition implies that $\nab = \de$. This implies for every $t
\in \R$ that
$$\pi(\si_t(a)) = \nab^{i t} \pi(a) \nab^{-it} = \de^{it} \pi(a) \de^{-it}
= \pi(\tau_t(a))$$ for every $t \in \R$ and $a \in A$.
\end{demo}

\medskip

\begin{corollary}
Consider a \cst-algebra $A$ and a faithful KMS-weight $\vfi$ on $A$. Then $\vfi$ as a
unique modular group.
\end{corollary}

\bigskip

In the last theorem of this section, we prove that our definition of KMS-weight is
equivalent with the usual one (introduced in \cite{Comb2}).

\begin{theorem}  \label{we4.thm1}
Consider a \cst-algebra $A$, a densely defined lower semi-continuous weight $\vfi$ on $A$
and a one-parameter group $\si$ on $A$. Then $\vfi$ is a KMS-weight with respect to $\si$
if and only if
\begin{enumerate}
\item We have for every $t \in \R$ that $\vfi \, \si_t = \vfi$.
\item For every $x,y \in \Nfi \cap \Nfi^*$, there exists a bounded continuous function
$f$ from $S(i)$ into $\,\C$ which is analytic on $S(i)^0$ and such that :
\begin{itemize}
\item We have for every $t \in \R$ that $f(t) = \vfi(\si_t(x) y)$.
\item We have for every $t \in \R$ that $f(t+i) = \vfi(y  \si_t(x))$.
\end{itemize}
\end{enumerate}
\end{theorem}
\begin{demo}
\begin{trivlist}
\item[$\,\,\,\Rightarrow$] This follows from proposition \ref{KMS}.
\item[$\,\,\,\Leftarrow$] Take a GNS-construction $(H,\la,\pi)$ for $\vfi$.  Define the
injective positive operator $\nab$ in $H$ such that $\nab^{it} \la(a) = \la(\si_t(a))$ for
every $t \in \R$.

\medskip

Define the set
$$ C = \{ \, a \in \Nfi \cap \Nfi^* \mid  a \text{ is
analytic with respect to } \si \text{ and }  \si_z(a) \text{ belongs to } \Nfi \cap \Nfi^*
\text{ for every } z \in \C \, \} \ .
$$
So $C$ is a sub$^*$-algebra of $A$ such that $\si_z(C) \subseteq C$ for every $z \in \C$.

\medskip

Choose $x,y \in C$. Because $x,y$ belong to $\Nfi \cap \Nfi^*$, we have by assumption the
existence of a bounded continuous function $f$ on $S(i)$ which is analytic on $S(i)^0$ and
satisfies :
\begin{itemize}
\item We have for every $t \in \R$ that $f(t) = \vfi(\si_t(x) y)$.
\item We have for every $t \in \R$ that $f(t+i) = \vfi(y  \si_t(x))$.
\end{itemize}
We know that $x^*$ belongs to $\Nfi$ and that $\si_{-i}(x^*)$ belongs to $\Nfi$. By
proposition \ref{we2.prop4}, we get that $\la(x^*)$ belongs to $D(\nab)$ and $\nab
\la(x^*) = \la(\si_{-i}(x^*))$. Define the function $g$ from $S(i)$ into $\C$ such that
$g(z) = \langle \la(y) , \nab^{i \, \overline{z}} \la(x^*) \rangle$ for every $z \in
S(i)$. Then $g$ is continuous on $S(i)$ and analytic on $S(i)^0$. We have for every $t \in
\R$ that
$$g(t) = \langle \la(y) , \nab^{i t} \la(x^*) \rangle = \langle \la(y) , \la(\si_t(x^*))
\rangle = \vfi(\si_t(x) y) = f(t) \ .$$
Therefore we must have that $f=g$ which implies that
$$\vfi(y x) = f(i) = g(i) = \langle \la(y) , \nab \la(x^*) \rangle = \langle \la(y) ,
\la(\si_{-i}(x^*)) \rangle = \vfi( \si_i(x) y ) \ . \text{\ \ \ \ \ (*)}$$

We know already that $C \subseteq \Nfi \cap D(\si_\frac{i}{2})$ and that
$\si_\frac{i}{2}(C)^* \subseteq \Nfi$. Just like in lemma \ref{we3.lem2} , we get that $C$
is a core for $\la$. Choose $a \in C$. If we replace in equation (*)  the element $x$  by
$\si_{-\frac{i}{2}}(a)$ and $y$ by $\si_\frac{i}{2}(a)^*$, we get that
$$ \vfi( \si_\frac{i}{2}(a)^* \, \si_{-\frac{i}{2}}(a) )
= \vfi( \si_\frac{i}{2}(a) \, \si_\frac{i}{2}(a)^*  ) \ . $$
We have that $\si_{-\frac{i}{2}}(a^* a) =  \si_\frac{i}{2}(a)^* \si_{-\frac{i}{2}}(a)$
which belongs to $\Mfi$. On the other hand, we know also that $a^* a$ belongs to $\Mfi$.
Hence, proposition \ref{we1.prop3} implies that $$\vfi(\si_\frac{i}{2}(a)^*
\si_{-\frac{i}{2}}(a)) = \vfi(\si_\frac{i}{2}(a^* a)) = \vfi(a^* a) \ .$$ Combining these
two equalities, we get that $\vfi(a^* a)
=  \vfi( \si_\frac{i}{2}(a) \si_\frac{i}{2}(a)^*  )$.
Therefore, proposition \ref{we4.prop12} implies that  $\vfi$ is KMS with respect to $\si$.
\end{trivlist}
\end{demo}

In order to prove the implication from the right to the left, we do not really need the
lower semi-continuity of $\vfi$. The closedness of the map $\la : \Nfi \rightarrow H$ is
also sufficient.

\bigskip

\section{The tensor product of KMS-weights}

Consider two \cst-algebras $A$ and $B$. Let $\vfi$ be a KMS-weight on $A$ with modular
group $\si$ and $\psi$ a KMS-weight on $B$ with modular group $\tau$.

\medskip

Take also a GNS-construction $(H,\la,\pi)$ for $\vfi$ and a GNS-construction $(K,\ga,\th)$
for $\psi$.
\begin{itemize}
\item With respect to the GNS-construction $(H,\la,\pi)$, we denote the modular conjugation of
$\vfi$ by $J$ and the modular operator of $\vfi$ by $\nab$.
\item With respect to the GNS-construction $(K,\ga,\th)$, we denote the modular conjugation of
$\psi$ by $I$ and the modular operator $\psi$ by $\de$.
\end{itemize}

\medskip

We will define the weight $\vfi \ot \psi$ and show that it is a KMS-weight with respect to
the obvious one-parameter group. The construction of the weight itself is due to Jan
Verding and can be found in \cite{Verd}. He did it in fact in a more general case. Again,
due to the lack of availability of this work, we will include the proofs. We only added
the KMS-characteristics.

\medskip

\begin{definition} \label{tens.def1}
We define the norm continuous one-parameter group $\si \ot \tau$ on $A \ot B$ such that
$(\si \ot \tau)_t = \si_t \ot \tau_t$ for every $t \in \R$.
\end{definition}

Consider $z \in \C$. Then it is easy to see that $\si_z \od \tau_z \subseteq (\si \ot
\tau)_z$. We will prove in \cite{JK3} (and this is not too difficult, using smearing
techniques) that $D(\si_z) \od D(\tau_z)$ is a core for $(\si \ot \tau)_z$. But we will
not need this fact in this paper.

\medskip

\begin{lemma}
The mapping $\la \od \ga : \Nfi \odot \Npsi \mapsto H \ot K$ is closable.
\end{lemma}
\begin{demo}
Define for every $\om \in \cF_\vfi$ the operator $S_\om \in B(H)$ such that $\langle S_\om
\, \la(a) , \la(b) \rangle = \om(b^* a)$ for every $a,b \in \Nfi$.
In a similar way, we define for every $\eta \in \cF_\psi$ the operator $T_\eta \in B(K)$
such that $\langle T_\eta \, \ga(a) , \ga(b) \rangle = \eta(b^* a)$ for every $a,b \in
\Npsi$.

\medskip

Take a sequence $(x_n)_{n=1}^\infty$ in $\Nfi \odot \Npsi$ and $v \in H \ot K$ such that
$(x_n)_{n=1}^\infty$ converges to $0$ and $\bigl((\la \od \ga)(x_n)\bigr)_{n=1}^\infty$
converges to $v$.

Fix $\om \in \cF_\vfi$ and $\eta \in \cF_\psi$ for the moment.

\medskip

Choose  $a \in \Nfi$ and $b \in \Npsi$. It is not difficult to see that
$$\langle (S_\om \ot T_\eta) \, (\la \ot \ga)(x_n) , \la(a) \ot \ga(b) \rangle
= (\om \ot \eta)( (a \ot b)^* \, x_n ) $$
for every $n \in \N$. This implies that the sequence
$$\bigl(\,\langle (S_\om \ot T_\eta) \, (\la \ot \ga)(x_n) , \la(a) \ot \ga(b)
\rangle \, \bigr)_{n=1}^\infty$$
converges to $0$. It is also clear that the sequence
$$\bigl(\,\langle (S_\om \ot T_\eta) \, (\la \ot \ga)(x_n) , \la(a) \ot \ga(b)
\rangle \, \bigr)_{n=1}^\infty$$
converges to $\langle (S_\om \ot T_\eta) \, v , \la(a) \ot \ga(b) \rangle$. So we see that
$\langle (S_\om \ot T_\eta) \, v , \la(a) \ot \ga(b) \rangle = 0$.

Hence we get that $(S_\om \ot T_\eta) \, v = 0$.

\medskip

Therefore, the remarks after result \ref{we1.res2} implies that $v = 0$.
\end{demo}

\medskip

So we can give the following definition.

\begin{definition} \label{tens.def2}
We define the closed linear map $\la \ot \ga$ from within $A \ot B$ into $H \ot K$ such
that $\Nfi \odot \Npsi$ is a core for $\la \ot \ga$ and $(\la \ot \ga)(a \ot b) = \la(a)
\ot \ga(b)$ for every $a \in \Nfi$ and $b \in \Npsi$.
\end{definition}

\begin{remark} \rm
\begin{itemize}
\item We have for every $x \in A \od B$ and $y \in \Nfi \odot \Npsi$ that $x \, y$ belongs
to $\Nfi \odot \Npsi$ and $(\la \ot \ga)(x \, y) = (\pi \ot \th)(x) \, (\la \ot \ga)(y)$.
This implies easily that $D(\la \ot \ga)$ is a left ideal in $A \ot B$ and that $(\la \ot
\ga)(x \, y) = (\pi \ot \th)(x) \, (\la \ot \ga)(y)$ for every $x \in A \ot B$ and $y \in
D(\la \ot \ga)$.
\item Choose $t \in \R$. Then it is easy to check for every $x \in \Nfi \odot
\Npsi$ that $(\si \ot \tau)_t (x)$ belongs to $\Nfi \odot \Npsi$ and $(\la \ot \ga)
\bigl((\si \ot \tau)_t (x)\bigr) = (\nab^{it} \ot \de^{it})(\la \ot \ga)(x)$.

This will imply for every $x \in D(\la \ot \ga)$ that $(\si
\ot \tau)_t (x)$ belongs to $D(\la \ot \ga)$ and $(\la \ot \ga)\bigl((\si
\ot \tau)_t (x)\bigr) = (\nab^{it} \ot \de^{it})(\la \ot \ga)(x)$.
\end{itemize}
\end{remark}

\bigskip

In the beginning of section \ref{we4}, we introduced the sub-$^*$-algebras $C_\vfi$ of
$\Nfi \cap \Nfi^*$ and $C_\psi$ of $\Npsi \cap \Npsi^*$.

We have also that $C_\vfi \subseteq D(\si_z)$ and $\si_z(C_\vfi) \subseteq C_\vfi$ and
that $C_\psi \subseteq D(\tau_z)$ and $\tau_z(C_\psi) \subseteq C_\psi$ for every $z \in
\C$.

We know also that $C_\vfi$ is a core for $\la$ and that $C_\psi$ is a core for $\ga$.

\medskip

Define the $^*$-algebra $C = C_\vfi \od C_\psi$. Then $C$ is a core for $\la \ot \ga$.

By the remarks after definition \ref{tens.def1}, we see that easily that $C$ consists of
analytic elements with respect to $(\si \ot \tau)_z$ and that $(\si \ot \tau)_z(C)
\subseteq C$ for every $z \in \C$.

\medskip

So we get in particular that $C \subseteq D((\si \ot \tau)_{\frac{i}{2}})$ and that $(\si
\ot \tau)_{\frac{i}{2}}(C)^* \subseteq C \subseteq D(\la \ot \ga)$.

\medskip

Using the remarks after definition \ref{tens.def1} and the fact that $\vfi(b^* a)
= \vfi(\si_\frac{i}{2}(a) \, \si_\frac{i}{2}(b)^*)$ for every $a,b \in \Nfi \cap D(\si_\frac{i}{2})$
and $\psi(b^* a) = \psi(\tau_\frac{i}{2}(a) \, \tau_\frac{i}{2}(b)^*)$ for every $a,b \in
\Npsi \cap D(\tau_\frac{i}{2})$, it is not difficult to check that
$$\|(\la \ot \ga)(x)\|^2 = \|(\la \ot \ga)\bigl((\si
\ot \tau)_{\frac{i}{2}}(x)^*\bigr) \|^2 $$ for every $x \in C$.

\medskip

So we see that the ingredients $A \ot B$, $H \ot K$, $\la \ot \ga$, $\pi \ot \th$, $\si
\ot \th$, $\nab \ot \de$ satisfy the conditions of the first construction procedure of
section \ref{we3}. Therefore we can use definition \ref{we3.def1}.

\begin{definition}
We define the weight $\vfi \ot \psi$ on $A \ot B$ such that $(H \ot K, \la \ot \ga, \pi
\ot \th)$ is a GNS-construction for $\vfi \ot \psi$.
\end{definition}

So we have in particular that $\cN_{\vfi \ot \psi} = D(\la \ot \ga)$, so $\Nfi \od \Npsi
\subseteq \cN_{\vfi \ot \psi}$.

\bigskip

Proposition \ref{we3.prop2} implies the following one :

\begin{proposition}
The weight $\vfi \ot \psi$ is a KMS-weight with modular group $\si \ot \tau$.
\end{proposition}

\medskip

Concerning the modular conjugation we have the following obvious result :

\begin{proposition}
With respect to the GNS-construction, the modular conjugation of $\vfi \ot \psi$ is given
by $J \ot I$, the modular operator of $\vfi \ot \psi$ is given by $\nab \ot \de$.
\end{proposition}

Call $J'$ the modular conjugation of $\vfi \ot \psi$. Using definition \ref{we4.def1} and
the remarks after definition \ref{tens.def1}, it is easy to check that $J' (\la(a) \ot
\ga(b)) = J \la(a) \ot I \ga(b)$ for every $a \in C_\vfi$ and $b \in C_\psi$. This implies
immediately that $J' = J \ot I$. The result about the modular operator follow from the
remarks after definition \ref{tens.def2}.

\bigskip

Definition \ref{tens.def2} also implies easily the following result :

\begin{result}
Consider $a \in \Mfi$ and $b \in \Mpsi$. Then $a \ot b$ belongs to $\cM_{\vfi \ot \psi}$
and $(\vfi \ot \psi)(a \ot b) = \vfi(a) \, \psi(b)$.
\end{result}

\medskip

There is also another characterization of $\vfi \ot \psi$.

\medskip

\begin{result}
We have the following properties :
\begin{enumerate}
\item We have for every $\om \in \cF_\vfi$ and $\eta \in \cF_\psi$ that
$\om \ot \eta$ belongs to $\cF_{\vfi \ot \psi}$.
\item We have for every $\om \in \cG_\vfi$ and $\eta \in \cG_\psi$ that
$\om \ot \eta$ belongs to $\cG_{\vfi \ot \psi}$.
\end{enumerate}
\end{result}
\begin{demo}
Consider $\om \in \cF_\vfi$ and $\eta \in \cF_\psi$.

Define the operator $S_\om \in B(H)$ such that $\langle S_\om
\, \la(a) , \la(b) \rangle = \om(b^* a)$ for every $a,b \in \Nfi$.

Define  also the operator $T_\eta \in B(K)$ such that $\langle T_\eta
\, \ga(a) , \ga(b) \rangle = \eta(b^* a)$ for every $a,b \in
\Npsi$.

Then it easy to see that $$\langle (S_\om \ot T_\eta) (\la \ot \ga)(x) , (\la \ot \ga)(y)
\rangle = (\om \ot \eta)(y^* x) $$
for every $x,y \in \Nfi \odot \Npsi$. As usual, this implies that

$$\langle (S_\om \ot T_\eta) (\la \ot \ga)(x) , (\la \ot \ga)(y)
\rangle = (\om \ot \eta)(y^* x) $$
for every $x,y \in \cN_{\vfi \ot \psi}$.

\medskip

Because $0 \leq S_\om \leq 1$ and $0 \leq T_\eta \leq 1$, this implies immediately that
$\om \ot \eta$ belongs to $\cF_{\vfi \ot \psi}$.

The second result follows from the first one.
\end{demo}

\medskip

\begin{proposition}
Consider $x \in (A \ot B)^+$. Then
$$(\vfi \ot \psi)(x) = \sup \, \{ \, (\om \ot \eta)(x) \mid \om \in \cF_\vfi,
\eta \in \cF_\psi \, \} \ .$$
\end{proposition}
\begin{demo}
By the previous proposition and the lower semi-continuity of $\vfi \ot \psi$, we get
immediately the inequality $$(\vfi \ot \psi)(x) \geq \sup \, \{ \, (\om \ot \eta)(x) \mid
\om \in \cF_\vfi, \eta \in \cF_\psi \, \} \ .$$

\medskip

By proposition \ref{we4.prop2}, we have truncating sequences $(u_k)_{k \in K}$ for $\vfi$
and $(v_l)_{l \in L}$ for $\psi$.

We define for every $k \in K$ the operator $S_k \in B(H)$ such that $S_k \la(a) = \la(a
u_k)$ for every $a \in \Nfi$.

In a similar way, we define for every $l \in L$ the operator $T_l \in B(K)$ such that $T_l
\ga(b) = \ga(b u_l)$ for every $b \in \Npsi$.

\medskip

Then we get immediately that
\begin{itemize}
\item We have that $u_k \ot v_l$ belongs to $\cN_{\vfi \ot \psi}$ for every $k \in K, l \in L$.
\item $\|S_k \ot T_l\| \leq 1$ for every $k \in K$ and $l \in L$.
\item The net $(u_k \ot v_l)_{(k,l) \in K \times L}$ converges strictly to 1.
\item The net $(S_k \ot T_l)_{(k,l) \in K \times L}$ converges strongly to 1.
\end{itemize}

Choose $k \in K, l \in L$.

Then it is easy to see that $(S_k \ot T_l)\,(\la \ot \ga)(y) = (\la \ot \ga)(y \, (u_k \ot
v_l))$ for every $y \in \Nfi \od \Npsi$.

As usual, this implies that $(S_k \ot T_l)\,(\la \ot \ga)(y) = (\la \ot \ga)(y \, (u_k \ot
v_l))$ for every $y \in \cN_{\vfi \ot \psi}$.

\medskip

So we see that the net $(u_k \ot v_l)_{(k,l) \in K \ot L}$ satisfies the conditions of
proposition \ref{extr1.prop2}.

\medskip

For every $k \in K$, we define $\om_k \in \cF_\vfi$ such that $\om_k(a) = \vfi(u_k^* a
u_k)$ for every $a \in A$.

For every $l \in L$, we define $\eta_l \in \cF_\psi$ such that $\eta_l(b) = \vfi(v_l^* b
v_l)$ for every $b \in B$.

It is clear that $(\om_k \ot \eta_l)(y) = (\vfi \ot \eta)\bigl((u_k \ot v_l)^* \, y \,
(u_k \ot v_l)\bigr)$ for every  $y \in A \ot B$.

\medskip

Therefore, proposition \ref{extr1.prop2} implies that
$$(\vfi \ot \psi)(x) = \sup \, \{ \, (\om_k \ot \eta_l)(x) \mid k \in K, l \in L \, \}
\leq \sup \, \{ \, (\om \ot \eta)(x) \mid \om \in \cF_\vfi, \eta \in \cF_\psi \, \} \ .$$
\end{demo}

\medskip

\begin{corollary}
We have the following properties :
\begin{enumerate}
\item We have for every $x \in (A \ot B)^+$ that the net
$\bigl(\, (\om \ot \eta)(x) \,\bigr)_{(\om,\eta) \in \cG_\vfi \times \cG_\psi}$ converges
to $(\vfi \ot \psi)(x)$ .
\item We have for every $x \in \cM_{\vfi \ot \psi}$ that the net
$\bigl(\, (\om \ot \eta)(x) \,\bigr)_{(\om,\eta) \in \cG_\vfi \times \cG_\psi}$ converges
to $(\vfi \ot \psi)(x)$ .
\end{enumerate}
\end{corollary}

\section{Absolutely continuous KMS-weights} \label{we5}

Consider a \cst-algebra $A$ and a KMS-weight $\vfi$ on $A$
with modular group $\si$.

Let $\sde$ be a strictly positive element affiliated with $A$ such that
there exists a strictly positive number $\lambda$ such that
$\si_t(\sde) = \lambda^t \, \sde$ for every $t \in \R$.

In this section we are going to use one of the construction procedures to define the
KMS-weight $\vfi(\sde^\frac{1}{2} \,  . \, \sde^\frac{1}{2})$ and prove some useful
properties about this weight.

It is not possible to use the definition of \cite{Pe-Tak}, because it is assumed in this
paper that $\si_t(\sde) = \sde$ for every $t \in \R$. Our construction procedures from
section \ref{we3} allows us to go about it in another way.

\medskip

First, we introduce some notations. Take a GNS-construction $(H,\la,\pi)$ for $\vfi$. We
denote the modular conjugation of $\vfi$ in this GNS-construction by $J$ and the modular
operator by $\nab$.

\medskip

We will extensively have to truncate $\sde$ in order to make things bounded. We will also
need this truncations to behave well with respect to $\si$. For this, we will use the
following elements in $M(A)$.

\begin{notation}
Consider $n \in \N$. Then we define $e_n \in M(A)$ such that
$$ e_n \, a = \frac{n}{\sqrt{\pi}} \int \exp(-n^2 t^2) \, \sde^{it} a \, dt  $$
for every $a \in A$.
\end{notation}

The behaviour of these elements is described in the next proposition.

\begin{proposition} \label{we5.prop2}
We have the following properties.
\begin{enumerate}
\item The sequence $(e_n)_{n=1}^\infty$ is bounded by 1 and converges
strictly to 1.
\item We have for every $n \in \N$ that $e_n$ is strictly analytic with respect tho $\si$.
\item Let $y \in \C$. Then the sequence $(\si_y(e_n))_{n=1}^\infty$ is bounded and
converges strictly to 1.
\item Consider $n \in \N$ and $z \in \C$.
 Then $e_n$  is a left and right multiplier of $\sde^z$ and $\sde^z e_n = e_n \sde^z$.
\item Let $n \in \N$ and $y,z \in \C$. Then $\sde^z  e_n$ is a left and right multiplier of
$\sde^{y}$ and $(\sde^z  e_n) \, \sde^y   = \sde^y \, (\sde^z  e_n) = \sde^{y+z} e_n$.
\item Consider $n \in \N$ and $y,z \in \C$. Then  $\si_y(e_n)$ is a left and right multiplier of $\sde^z$
and $\si_y(e_n) \, \sde^z = \sde^z \, \si_y(e_n)$.
\item Consider $n \in \N$ and $z \in \C$. Then $\sde^z e_n$ is strictly analytic
with respect to $\si$ and $\si_y( \sde^z e_n ) = \lambda^{y z} \, \sde^z \, \si_y(e_n)$
for every $y \in \C$.
\end{enumerate}
\end{proposition}

The results of this proposition will be frequently used in the rest of this section
without reference to this proposition. The proof of it is a standard exercise in the use
of analytic continuations and will be left out.

\bigskip

We define the set
$$N_0 = \{\, a \in A \mid a \text{ is a left multiplier of } \sde^\frac{1}{2}
\text{ and } a \, \sde^\frac{1}{2} \text{ belongs to } \Nfi \,\} \ .$$
Introducing this set is of course a very natural thing to do in this case.

It is clear that this set is a left ideal in $M(A)$.

\begin{lemma}
The mapping $N_0 \rightarrow H : a \mapsto \la( a \, \sde^\frac{1}{2}  )$ is closable
\end{lemma}
\begin{demo}
Choose a sequence $(x_k)_{k=1}^\infty$ in $N_0$ and $v \in H$ such that
$(x_k)_{k=1}^\infty$ converges to $0$ and $\bigl(\la(x_k \,
\sde^\frac{1}{2})\bigr)_{k=1}^\infty$ converges to $v$. Choose $n \in \N$.

\medskip

Take $l \in \N$.

Because $e_n$ is a left multiplier of $\sde^\frac{1}{2}$, we have that
$x_l \, e_n$ is a left multiplier of $\sde^\frac{1}{2}$ and
$$(x_l \, e_n) \, \sde^\frac{1}{2} = x_l \, (e_n  \, \sde^\frac{1}{2})
= x_l \, (\sde^\frac{1}{2} e_n) \ . \hspace{1.5cm} \text{(a)}$$
We know that $x_l$ is a left multiplier of $\sde^\frac{1}{2}$, so this last equality
implies also that $(x_l \, e_n) \, \sde^\frac{1}{2} = (x_l \, \sde^\frac{1}{2}) \,  e_n$.

Because $x_l \, \sde^\frac{1}{2}$ belongs to $\Nfi$ and $e_n$ belongs to
$D(\overline{\si}_\frac{i}{2})$, this equality and proposition \ref{we4.prop1} imply that
$(x_l \, e_n) \, \sde^\frac{1}{2}$ belongs to $\Nfi$ and
$$\la\bigl((x_l \, e_n) \, \sde^\frac{1}{2}\bigr) = J \pi(\si_\frac{i}{2}(e_n))^* J
\la(x_l \, \sde^\frac{1}{2}) = J \pi(\si_\frac{i}{2}(e_n))^* J \ga(x_l) \ .
\hspace{1.5cm} \text{(b)}$$

\medskip

Equality (a) implies that $\bigl( (x_k \, e_n) \, \sde^\frac{1}{2}\bigr)_{k=1}^\infty$
converges to 0, whereas equality (b) implies that \newline $\bigl(\la((x_k \, e_n) \,
\sde^\frac{1}{2})\bigr)_{k=1}^\infty$ converges to $J \pi(\si_\frac{i}{2}(e_n))^* J v$.

Therefore, the closedness of $\la$ implies that $J \pi(\si_\frac{i}{2}(e_n))^* J v = 0$

\medskip

Because $\bigl(\pi(\si_\frac{i}{2}(e_n))\bigr)_{n=1}^\infty$ converges strongly$^*$ to 1,
this last equation implies that $v=0$.
\end{demo}

Therefore, we can give the following definition.

\begin{definition}
We define $\ga$ as the closure of the linear mapping $N_0 \rightarrow H : a \mapsto
\la(a \, \sde^\frac{1}{2})$. The domain of $\ga$ will be denoted by $N$.
\end{definition}

So $\ga$ is a closed linear mapping from $N$ into $H$ such that $N_0$ is a core for $\ga$
and $\ga(a) = \la(a \, \sde^\frac{1}{2})$ for every $a \in N_0$. It is also easy
to check that $N$ is a left ideal in $M(A)$ and that $\ga(x a) = \pi(x) \ga(a)$ for every
$x \in M(A)$ and $a \in N$.

\begin{lemma} \label{we5.lem7}
Consider $a \in \Nfi$, $n \in \N$ and $z \in \C$. Then $a \,(\sde^z e_n)$ belongs to $N_0$
and $\ga(a \,(\sde^z  e_n)) = J \pi(\si_\frac{i}{2}(\sde^{z+\frac{1}{2}}  e_n))^* J \la(a)$.
\end{lemma}
\begin{demo}
We  know that $\sde^z e_n$ is a left multiplier of $\sde^\frac{1}{2}$ and
$(\sde^z e_n) \, \sde^\frac{1}{2} = \sde^{z+\frac{1}{2}} e_n$.

This implies that $a \,(\sde^z e_n)$ is a left multiplier of $\sde^\frac{1}{2}$ and $(a
\, (\sde^z e_n)) \, \sde^\frac{1}{2} = a \, (\sde^{z+\frac{1}{2}}  e_n)$. Because $a$
belongs to $\Nfi$ and $\sde^{z+\frac{1}{2}} e_n$ belongs to
$D(\overline{\si}_\frac{i}{2})$, proposition \ref{we4.prop5} implies that $(a
\, (\sde^z  e_n)) \, \sde^\frac{1}{2}$ belongs to $\Nfi$ and   $\la\bigl(\, (a
\,(\sde^z e_n)) \, \sde^\frac{1}{2}\,\bigr) = J
\pi(\si_\frac{i}{2}(\sde^{z+\frac{1}{2}}  e_n))^* J \la(a) \ . $
This gives by definition that $a \,(\sde^z  e_n)$ belongs to $N_0$ and $\ga(a \,(\sde^z
e_n)) = J \pi(\si_\frac{i}{2}(\sde^{z+\frac{1}{2}}  e_n))^* J \la(a)$.
\end{demo}

\begin{result}
The set $N$ is dense in $A$  and the set $\ga(N)$ is dense in $H$.
\end{result}
\begin{demo}
Choose $a \in \Nfi$. The previous lemma implies for every $n \in \N$ that $a \,e_n$ belongs
to $N_0$. Because $(e_n)_{n=1}^\infty$ converges strictly to 1, we see that $(a \,
e_n)_{n=1}^\infty$ converges to $a$. Hence, the density of $\Nfi$ in $A$ implies that
$N_0$ is dense in $A$.

\medskip

Choose $a \in \Nfi$. The previous lemma implies for every $n \in \N$ that $a \,
(\sde^{-\frac{1}{2}}  e_n)$ belongs to $N_0$ and $\ga(a \,(\sde^{-\frac{1}{2}}  e_n)) = J
\pi(\si_\frac{i}{2}(e_n))^* J \la(a)$. Because $\bigl(\pi(\si_\frac{i}{2}(e_n))
\bigr)_{n=1}^\infty$
converges strongly$^*$ to 1. This implies that $\bigl(\ga(a \,(\sde^{-\frac{1}{2}}
e_n))\bigr)_{n=1}^\infty$ converges to $\la(a)$. Hence, the density of $\la(\Nfi)$ in $H$
implies that $\ga(N_0)$ is dense in $H$.
\end{demo}

So we have objects $A$,$H$,$N$,$\la$,$\pi$ on which we want to apply the second procedure
of section \ref{we3}. In order to do so, we have to introduce the KMS-characteristics.

\begin{definition}
We define the norm continuous one-parameter group $\si'$ on $A$ such that
\newline $\si_t'(a) = \sde^{it} \si_t(a) \sde^{-it}$ for every $t \in A$ and $a \in A$.
\end{definition}

It is not difficult to see that $\si_t'(\sde) = \lambda^t \, \sde$ for every $t \in \R$.

\medskip

Result \ref{we4.res2} implies the following lemma.

\begin{lemma} \label{we5.lem1}
We have for every $s,t \in \R$ that $J \pi(\sde)^{it} J$ and $\pi(\sde)^{is}$ commute.
\end{lemma}

The following result is true because $\pi(\si_s(a)) = \nab^{is} \pi(a) \nab^{-i s}$ for
every $a \in A$, $s \in \R$ and the fact that $\si_s(\sde^{i t}) = \lambda^{i s t} \,
\sde^{i t}$ for $s,t \in \R$ (which is true by assumption).

\begin{lemma}
Consider $s,t \in \R$. Then $\nab^{is} \, \pi(\sde)^{it} =  \lambda^{i s t} \,
\pi(\sde)^{it} \, \nab^{is}$
\end{lemma}

This will imply the following lemma.

\begin{lemma} \label{we5.lem2}
We have for every $s,t \in \R$ that $J \pi(\sde)^{it} J \, \pi(\sde)^{it}$ and $\nab^{is}$
commute.
\end{lemma}
\begin{demo}
We have that
\begin{eqnarray*}
J \pi(\sde)^{it} J \, \pi(\sde)^{it} \, \, \nab^{is} & = & \lambda^{-i s t} \, J
\pi(\sde)^{it} J \, \nab^{is} \, \pi(\sde)^{it} = \lambda^{-i s t} \, J \pi(\sde)^{it} \,
\nab^{is} \, J \, \pi(\sde)^{it} \\ & \stackrel{(*)}{=} & \lambda^{-i s t} \, \lambda^{i s
t} \, J \, \nab^{is} \, \pi(\sde)^{it} J \, \pi(\sde)^{it}
= \nab^{is} \, \, J \pi(\sde)^{it} J \, \pi(\sde)^{it} \ .
\end{eqnarray*}
In (*), we used that $J$ is antilinear.
\end{demo}

Lemmas \ref{we5.lem1} and \ref{we5.lem2} imply that the mapping $\R \rightarrow B(H): t
\mapsto J \pi(\sde)^{it} J $   $\pi(\sde)^{it} \, \nab^{it}$ is a strongly continuous
unitary group representation on $H$. By the Stone theorem, this justifies the following
definition.

\begin{definition}
We define the positive injective operator $\nabp$ (pronounced 'nabla prime') in $H$ such
that $\nabp^{it} = J \pi(\sde)^{it} J \, \pi(\sde)^{it} \, \nab^{it}$ for every $t
\in \R$.
\end{definition}

\begin{proposition} \label{we5.prop1}
Consider $a \in \Nfi$ and $t \in \R$. Then $\si_t'(a)$ belongs to $\Nfi$
and $\la(\si_t'(a)) = \lambda^\frac{t}{2} \, \nabp^{it} \la(a)$.
\end{proposition}
\begin{demo}
We know by definition that $\si_t(a)$ belongs to $\Nfi$ and $\la(\si_t(a)) = \nab^{it}
\la(a)$. This implies that $\sde^{it} \si_t(a)$ belongs to $\Nfi$ and that $\la(\sde^{it}
\si_t(a)) = \pi(\sde)^{it} \la(\si_t(a)) = \pi(\sde)^{it} \, \nab^{it} \la(a)$.

Therefore, proposition \ref{we4.prop9} implies that $\sde^{it} \si_t(a) \sde^{-it}$
belongs to $\Nfi$ and
$$\la(\sde^{it} \si_t(a) \sde^{-it}) = \lambda^\frac{t}{2} \, J \pi(\sde)^{it} J
\la( \sde^{it} \si_t(a))  = \lambda^\frac{t}{2} \, J \pi(\sde)^{it} J \, \pi(\sde)^{it}
\, \nab^{it} \la(a)  \  .$$
Looking at the definitions of $\si_t'$ and $\nabp$, the lemma follows.
\end{demo}

\begin{proposition}
Consider $a \in N$ and $t \in \R$. Then $\si_t(a)$ belongs to $N$
and $\ga(\si_t(a)) = \lambda^{-\frac{t}{2}} \, \nab^{it} \ga(a)$.
\end{proposition}
\begin{demo}
Choose $b \in N_0$.

Because $b$ is a left multiplier of $\sde^\frac{1}{2}$
and $\si_t(\sde) = \lambda^t \, \sde$, we get that $\si_t(b)$ is a left multiplier of
$\sde^\frac{1}{2}$
and $$\si_t(b) \, \sde^\frac{1}{2} = \lambda^{-\frac{t}{2}} \,
\si_t( b \, \sde^\frac{1}{2} ) \ . $$
Because $b \, \sde^\frac{1}{2}$ belongs to $\Nfi$, this implies that
$\si_t(b) \, \sde^\frac{1}{2}$ belongs to $\Nfi$ and
$$\la(\si_t(b) \, \sde^\frac{1}{2}) = \lambda^{-\frac{t}{2}} \,
\la(\si_t(b \, \sde^\frac{1}{2}))
= \lambda^{-\frac{t}{2}} \, \nab^{it} \la(b \, \sde^\frac{1}{2})
= \lambda^{-\frac{t}{2}} \, \nab^{it} \ga(b) \ . $$

So we see that $\si_t(b)$ belongs to $N_0$ and
$$\ga(\si_t(b)) = \la(\si_t(b) \, \sde^\frac{1}{2})
= \lambda^{-\frac{t}{2}} \, \nab^{it} \ga(b) \ . $$

The result follows easily because $N_0$ is a core for $\ga$.
\end{demo}

\begin{proposition}
Consider $a \in N$ and $t \in \R$. Then $\si_t'(a)$ belongs to $N$
and $\ga(\si_t'(a)) = \nabp^{it} \ga(a)$.
\end{proposition}
\begin{demo}
Choose $b \in N_0$.

Because $b$ is a left multiplier of $\sde^\frac{1}{2}$
and $\si_t'(\sde) = \lambda^t \, \sde$, we get that $\si_t'(b)$ is a left multiplier of
$\sde^\frac{1}{2}$ and
$$\si_t'(b) \, \sde^\frac{1}{2} = \lambda^{-\frac{t}{2}} \, \si_t'(b \,
\sde^\frac{1}{2}) \ . $$
Because $b \, \sde^\frac{1}{2}$ belongs to $\Nfi$, this implies that
$\si_t'(b) \, \sde^\frac{1}{2}$ belongs to $\Nfi$ and
$$\la(\si_t'(b) \, \sde^\frac{1}{2}) = \lambda^{-\frac{t}{2}} \,
\la(\si_t'(b\,\sde^\frac{1}{2})) \stackrel{(*)}{=} \lambda^{-\frac{t}{2}} \,
 \lambda^\frac{t}{2} \, \nabp^{it} \la(b \, \sde^\frac{1}{2}) =  \nabp^{it} \ga(b) \ , $$
where we used proposition \ref{we5.prop1}  in equality (*).
The results follows easily because $N_0$ is a core for $\ga$.
\end{demo}

\bigskip

\begin{proposition}
Consider $x \in N$ and $a \in D(\si_\frac{i}{2}')$. Then $x a$ belongs to $N$ and
$$\ga(x a) = J \pi(\si_\frac{i}{2}'(a))^* J \ga(x) \ .$$
\end{proposition}
\begin{demo}
Define the norm continuous one-parameter group $\tau$ on $A$ such that
$\tau_t(a) = \sde^{i t} a \sde^{-i t}$ for every $a \in A$ and $t \in \R$.
Then we have that  $\si_s \tau_t = \tau_t \si_s$ for every $s,t \in \R$
and $\si_t'(a) = \tau_t( \si_t(a))$ for every $t \in \R$.

This implies that $\si_\frac{i}{2} \, \tau_\frac{i}{2}$ is closable and that the closure
is equal to $\si_\frac{i}{2}'$ \ (see \cite{JK3}).

\medskip

Take $y \in N_0$. So $y$ is a left multiplier of $\sde^\frac{1}{2}$ and
$y \, \sde^\frac{1}{2}$ belongs to $\Nfi$.

\begin{list}{}{\setlength{\leftmargin}{.4 cm}}

\item Choose $b \in D(\si_\frac{i}{2} \, \tau_\frac{i}{2})$.

Because $b$ belongs to $D(\tau_\frac{i}{2})$, the element $b$ is a middle multiplier of
$\sde^{-\frac{1}{2}},\,\sde^\frac{1}{2}$ and $\sde^{-\frac{1}{2}} b \, \sde^\frac{1}{2} =
\tau_\frac{i}{2}(b)$.

We have by assumption that $y$ is a left multiplier of $\sde^\frac{1}{2}$. This implies
that $y \, \sde^\frac{1}{2}$ is a left multiplier of $\sde^{-\frac{1}{2}}$ and $(y \,
\sde^\frac{1}{2}) \, \sde^{-\frac{1}{2}} = y$.

This in turn implies that $((y \, \sde^\frac{1}{2})\,\sde^{-\frac{1}{2}}) \, b$ is a left
multiplier of $\sde^\frac{1}{2}$ and
$$[((y \, \sde^\frac{1}{2})\,\sde^{-\frac{1}{2}}) \, b] \, \sde^\frac{1}{2}
= (y \, \sde^\frac{1}{2}) \, (\sde^{-\frac{1}{2}} b \, \sde^\frac{1}{2})
= (y \, \sde^\frac{1}{2}) \, \tau_\frac{i}{2}(b) \ .$$

Combining this with the fact that $((y \, \sde^\frac{1}{2})\,\sde^{-\frac{1}{2}}) b  = y
b$, we see that $y b$ is a left multiplier of $\sde^\frac{1}{2}$ and $(y b) \,
\sde^\frac{1}{2} = (y \, \sde^\frac{1}{2}) \, \tau_\frac{i}{2}(b)$.

Because $y \, \sde^\frac{1}{2}$ belongs to $\Nfi$ and $\tau_\frac{i}{2}(b)$
belongs to $D(\si_\frac{i}{2})$, this last equation implies that $(y b) \,
\sde^\frac{1}{2}$ belongs to $\Nfi$ and
$$\la\bigl((y b) \, \sde^\frac{1}{2}\bigr) = J \pi\bigl(\si_\frac{i}{2}
(\tau_\frac{i}{2}(b))\bigr)^* J  \la( y \, \sde^\frac{1}{2} )
 = J \pi\bigl(\si_\frac{i}{2}(\tau_\frac{i}{2}(b))\bigr)^* J \ga(y) \ .$$

Therefore, $y b$ belongs to $N_0$ and $\ga(y b) = J
\pi(\si_\frac{i}{2}(\tau_\frac{i}{2}(b)))^* J  \ga(y)$

So we see that $\ga(y b) = J \pi(\si_\frac{i}{2}'(b))^* J
\ga(y)$. \ \ \ \ \ \ (*)

\end{list}

Because $D(\si_\frac{i}{2} \, \tau_\frac{i}{2})$ is a core for $\si_\frac{i}{2}'$, there exists a sequence $(a_n)_{n=1}^\infty$ in
$D(\si_\frac{i}{2} \, \tau_\frac{i}{2})$ such that $(a_n)_{n=1}^\infty$ converges to $a$
and $(\si_\frac{i}{2}'(a_n))_{n=1}^\infty$ converges to $\si_\frac{i}{2}'(a)$.

\medskip

By result (*) we have for every $n \in \N$ that
$y a_n$ belongs to $N_0$ and $$\ga(y a_n) =   J \pi(\si_\frac{i}{2}'(a_n))^* J \ga(y)$$

From this, we conclude that the sequence $\bigl(\ga(y a_n)\bigr)_{n=1}^\infty$ converges
to $J \pi(\si_\frac{i}{2}'(a))^* J \ga(y)$. It is also clear that $(y a_n)_{n=1}^\infty$
converges to $y a$ so the closedness of $\ga$ implies that $y a$ belongs to $N$ and $\ga(y
a) = J^* \pi(\si_\frac{i}{2}'(a))^* J \ga(y)$.

\medskip

The proposition follows because $N_0$ is a core for $\ga$.
\end{demo}

\bigskip

At last we are in a position to define $\vfi(\sde^\frac{1}{2} \, . \, \sde^\frac{1}{2})$.

\begin{definition}
We have proven in the previous part of the section that the ingredients
$A$,$H$,$N$,$\la$,$\pi$,$\si'$,$\nabp$ and $J$ satisfy the conditions of the second
construction procedure of section \ref{we3}. The weight which is obtained from this
ingredients is denoted by $\vfi_\sde$ (see definition \ref{we3.def2}).
\end{definition}

We want to forget about this definition, but will instead repeat the determining properties.
Previously, we have defined a norm-continuous one-parameter group $\si'$
on $A$ such that $\si_t'(a) = \sde^{it} \si_t(a) \sde^{-it}$ for every
$t \in \R$ and $a \in A$. Then

\begin{theorem}
We have that $\vfi_\sde$ is a KMS-weight with modular group $\si'$.
\end{theorem}

The weight $\vfi_\sde$ is completely determined by the following proposition.

\begin{theorem}
We have that $(H,\ga,\pi)$ is a GNS-construction for $\vfi_\sde$
\end{theorem}

Remember that the set
$$N_0 = \{\, a \in A \mid a \text{ is a left multiplier of } \sde^\frac{1}{2}
\text{ and }  a \, \sde^\frac{1}{2} \text{ belongs to } \Nfi \,\} \ .$$
is a core for $\ga$ and that $\ga(a) = \la(a \, \sde^\frac{1}{2})$ for every $a
\in N_0$.

\bigskip

Repeating the results of lemmas \ref{we5.lem1} and \ref{we5.lem2} gives
\begin{itemize}
\item We have for every $s,t \in \R$ that $J \pi(\sde)^{it} J$ and $\pi(\sde)^{is}$ commute.
\item We have for every $s,t \in \R$ that $J \pi(\sde)^{it} J \, \pi(\sde)^{it}$ and
$\nab^{is}$ commute.
\end{itemize}
Also, the injective positive operator $\nabp$ in $H$ was defined in such a way that
$\nabp^{it} = J \pi(\sde)^{it} J \pi(\sde)^{it} \nab^{it}$ for every $t \in \R$. Then :

\begin{proposition}
We have that $\nabp$ is the modular operator for $\vfi_\sde$ in the GNS-construction
$(H,\ga,\pi)$.
\end{proposition}

Later, we will prove that the modular conjugation for $\vfi_\sde$
is a scalar multiple of $J$.

\bigskip

The following relative invariance properties are also valid.
\begin{itemize}
\item Consider $a \in \Nfi$ and $t \in \R$. Then $\si_t'(a)$ belongs to $\Nfi$
and $\la(\si_t'(a)) = \lambda^\frac{t}{2} \, \nabp^{it} \la(a)$.
\item Consider $a \in \cN_{\vfi_\sde}$ and $t \in \R$. Then $\si_t(a)$ belongs to
$\cN_{\vfi_\sde}$ and $\ga(\si_t(a)) = \lambda^{-\frac{t}{2}} \, \nab^{it} \ga(a)$.
\end{itemize} These imply the following proposition.

\begin{proposition} Consider $t \in \R$. Then $\vfi \, \si_t' = \lambda^t \, \vfi$
and $\vfi_\sde \, \si_t = \lambda^{-t} \, \vfi_\sde$.
\end{proposition}

Because $(H,\ga,\pi)$ is a GNS-construction for $\vfi_\sde$, proposition \ref{we4.prop7}
implies the following result.

\begin{proposition}
The weight $\vfi_\sde$ is faithful $\Leftrightarrow$ $\vfi$ is faithful
\end{proposition}

\bigskip

Because $\si_t'(\sde) = \lambda^t \, \sde$ for every $t \in \R$, also the following
proposition will be true (see proposition \ref{we5.prop2}).

\begin{proposition}
We have the following properties.
\begin{enumerate}
\item We have for every $n \in \N$ that $e_n$ is strictly analytic with respect tho $\si'$.
\item Let $y \in \C$. Then the sequence $(\si_y'(e_n))_{n=1}^\infty$ is bounded and converges
strictly to 1.
\item Consider $n \in \N$ and $y,z \in \C$. Then $\si_y'(e_n)$ is a left and right
multiplier of $\sde^z$ and $\si_y'(e_n) \, \sde^z = \sde^z \, \si_y'(e_n)$.
\item Consider $n \in \N$ and $z \in \C$. Then $\sde^z e_n$ is strictly analytic
with respect to $\si'$ and $\si_y'( \sde^z e_n ) = \lambda^{y z} \, \sde^z \, \si_y'(e_n)$
for every $y \in \C$.
\end{enumerate}
\end{proposition}

\bigskip

In the next part of this section, we will prove a formula which formally says that
$\vfi_\sde(x) = \vfi(\sde^\frac{1}{2}\, x \, \sde^\frac{1}{2})$ (and something more
general).

\begin{lemma} \label{we5.lem9}
We have the following properties :
\begin{itemize}
\item We have for every $a \in \Nfi$, $n \in \N$ and $z \in \C$ that $a \, (\sde^z e_n)$
belongs to $\Nfi \cap {\cal N}_{\vfi_\sde}$.
\item We have for every $x \in \Mfi$, $m,n \in \N$ and $y,z \in \C$ that $(\sde^y e_m)\, x
\, (\sde^z e_n)$ belongs to $\Mfi \cap {\cal M}_{\vfi_\sde}$.
\end{itemize}
\end{lemma}
\begin{demo}
\begin{itemize}
\item From lemma \ref{we5.lem7}, we know already that $a \, (\sde^z e_n)$
belongs to ${\cal N}_{\vfi_\sde}$. Because $a$ belongs to $\Nfi$ and $\sde^z e_n$ belongs
to $D(\si_\frac{i}{2})$, proposition \ref{we4.prop5} implies also that $a \, (\sde^z e_n)$
belongs to $\Nfi$.
\item This follows immediately from the first property.
\end{itemize}
\end{demo}

\begin{lemma}   \label{we5.lem5}
Consider $x \in \Mfi$. Then $\bigl(\vfi(e_n \, x \, e_n)\bigr)_{n=1}^\infty$ converges to
$\vfi(x)$.
\end{lemma}
\begin{demo}
Choose $a,b \in \Nfi$. We have for every $n \in \N$ that $a \, e_n$, $b \, e_n$ belong to
$\Nfi$ and $\la(a \, e_n) = J \pi(\si_\frac{i}{2}(e_n))^* J \la(a)$, \ $\la(b \, e_n) = J
\pi(\si_\frac{i}{2}(e_n))^* J \la(b)$ which implies that
$$\vfi(e_n \, b^* a \, e_n) = \langle J \pi(\si_\frac{i}{2}(e_n))^* J \la(a) ,
J \pi(\si_\frac{i}{2}(e_n))^* J \la(b) \rangle \ .$$ Because
$\bigl(\pi(\si_\frac{i}{2}(e_n))\bigr)_{n=1}^\infty$ converges strongly$^*$ to 1, we see
that $\bigl(\vfi(e_n \, b^* a \, e_n)\bigr)_{n=1}^\infty$ converges to $\langle \la(a) ,
\la(b) \rangle$ which is equal to $\vfi(b^* a)$.
\end{demo}

\bigskip

\begin{lemma} \label{we5.lem3}
Consider $x \in \cN_{\vfi_\sde}$, $n \in \N$ and $z \in \C$.
Then $x \, (\sde^z e_n)$  belongs to $\Nfi$.
\end{lemma}
\begin{demo}
Because $\sde^z e_n$ belongs to $D(\overline{\si'}_\frac{i}{2})$, we have immediately that
$x \, (\sde^z e_n)$ belongs to  $\cN_{\vfi_\sde}$.

By definition, there exists a sequence $(x_k)_{k=1}^\infty$ in $N_0$ such that
$(x_k)_{k=1}^\infty$ converges to $x$ and  \newline $(\ga(x_k))_{k=1}^\infty$ converges to
$\ga(x)$.

\medskip

Choose $l \in \N$. Because $x_l$ is a left multiplier of $\sde^\frac{1}{2}$ and
$\sde^{z-\frac{1}{2}} e_n$ is a right multiplier of $\sde^\frac{1}{2}$, we have that
$$(x_l \, \sde^\frac{1}{2}) (\sde^{z-\frac{1}{2}} e_n)
= x_l \, (\sde^\frac{1}{2} \, (\sde^{z-\frac{1}{2}} e_n)) = x_l \, (\sde^z e_n) \ . $$

Because $x_l \, \sde^\frac{1}{2}$ belongs to $\Nfi$ and $\sde^{z-\frac{1}{2}} e_n$ belongs to
$D(\overline{\si}_\frac{i}{2})$, this implies that $x_l \, (\sde^z e_n)$ belongs to $\Nfi$ and
$$\la(x_l \, (\sde^z e_n)) = J \pi(\si_\frac{i}{2}(\sde^{z-\frac{1}{2}} e_n))^* J
\la(x_l \, \sde^\frac{1}{2}) = J \pi(\si_\frac{i}{2}(\sde^{z-\frac{1}{2}} e_n))^*
J \ga(x_l) \ . $$

\medskip

From this last equation, we infer that $\bigl(\la(x_k \,(\sde^z e_n))\bigr)_{k=1}^\infty$
converges to $J \pi(\si_\frac{i}{2}(\sde^{z-\frac{1}{2}} e_n))^* J \ga(x)$. It is also
clear that $\bigl(x_k \,(\sde^z e_n)\bigr)_{k=1}^\infty$ converges to $x \, (\sde^z e_n)$. Therefore, the
closedness of $\la$ implies that $x \, (\sde^z e_n)$ belongs to $\Nfi$.
\end{demo}

Using this lemma, we can even do better.

\begin{lemma} \label{we5.lem4}
Consider $x \in \cN_{\vfi_\sde}$, $n \in \N$ and $z \in \C$. Then $x \, (\sde^z e_n)$ belongs
to $N_0$.
\end{lemma}
\begin{demo}
Because $\sde^z e_n$ is a left multiplier of $\sde^\frac{1}{2}$, we get that
$x \, (\sde^z e_n)$ is a left multiplier of $\sde^\frac{1}{2}$ and
$$(x \, (\sde^z e_n)) \, \sde^\frac{1}{2}
= x \, ((\sde^z e_n) \, \sde^\frac{1}{2}) = x \, (\sde^{z+\frac{1}{2}} e_n)  $$
which belongs
to $\Nfi$ by the previous lemma. The result follows by the definition of $N_0$.
\end{demo}

\medskip

The following results can be proven in a similar way as for elements in
$\Mfi$ and $\Nfi$.

\begin{lemma} \label{we5.lem8}
We have the following properties.
\begin{itemize}
\item We have for every $a \in {\cal N}_{\vfi_\sde}$, $n \in \N$ and $z \in \C$ that
$a \, (\sde^z e_n)$ belongs to $\Nfi \cap {\cal N}_{\vfi_\sde}$.
\item We have for every $x \in {\cal M}_{\vfi_\sde}$, $m,n \in \N$ and $y,z \in \C$ that
$(\sde^y e_m) \, x \, (\sde^z e_n)$ belongs to $\Mfi \cap {\cal M}_{\vfi_\sde}$.
\end{itemize}
\end{lemma}

\begin{lemma} \label{we5.lem6}
Consider $x \in {\cal M}_{\vfi_\sde}$. Then $\bigl(\vfi_\sde(e_n\,x\,e_n)\bigr)_{n=1}^\infty$
converges to $\vfi_\sde(x)$.
\end{lemma}

\medskip

We have also the next result.

\begin{corollary} \label{we5.cor1}
The set $\Nfi \cap N_0$ is a core for $\ga$ and a core for $\la$.
\end{corollary}
\begin{demo}
Choose $x \in \cN_{\vfi_\sde}$. We have by the previous lemmas that $x \, e_n$ belongs to
$\Nfi \cap N_0$ for every $n \in \N$. By proposition \ref{we4.prop5}, we know that $\ga(x
\, e_n) = J' \pi(\si_\frac{i}{2}'(e_n))^* J' \ga(x)$ for every $n \in \N$ (here $J'$ denotes
the modular conjugation of $\vfi_\sde$). This implies that $(\ga(x \, e_n))_{n=1}^\infty$
converges to $\ga(x)$. It is also clear that $(x \, e_n)_{n=1}^\infty$ converges to $x$.

The statement about $\la$ is proven in a similar way.
\end{demo}

\medskip

\begin{lemma}
Consider $x \in \cM_{\vfi_\sde}$, $z \in \C$ and $n \in N$. Then $e_n \, x \, e_n$ belongs to
$\cM_{\vfi_\sde}$, $(\sde^{1-z} e_n) \, x \, (\sde^z e_n)$ belongs to $\cM_\vfi$ and
$$\vfi_\sde(e_n \, x \, e_n) = \lambda^{i(\frac{1}{2}-z)} \, \vfi\bigl( \, (\sde^{1-z} e_n)\, x
\, (\sde^z e_n)\, \bigr) \ . $$
\end{lemma}
\begin{demo}
Choose $a,b \in \cN_{\vfi_\sde}$. By lemma \ref{we5.lem4}, we know that $a \, e_n, b \, e_n$
belong to $N_0$. This implies that
\begin{eqnarray*}
\vfi_\sde( e_n \, b^* a \, e_n) & = & \langle \ga(a \, e_n) , \ga(b \, e_n) \rangle
= \langle \la\bigl( (a \, e_n) \, \sde^\frac{1}{2}\bigr) ,
 \la\bigl((b \, e_n) \, \sde^\frac{1}{2}\bigr) \rangle \\
& = & \langle \la\bigl( a \, (\sde^\frac{1}{2} e_n) \bigr) , \la\bigl( b \, (\sde^\frac{1}{2} e_n) \bigr)
\rangle = \vfi\bigl(\,(\sde^\frac{1}{2} e_n) \, b^* a \, (\sde^\frac{1}{2} e_n)\,\bigr) \ .
\end{eqnarray*}
We know that $\sde^{1 - z} e_n$ is a right multiplier of $\sde^{z - \frac{1}{2}}$ and
that $\sde^{z - \frac{1}{2}} \, (\sde^{1 - z} e_n) = \sde^\frac{1}{2} e_n$.

This implies that $(\sde^{1 - z} e_n)\, b^* a \,(\sde^\frac{1}{2} e_n)$ is a right
multiplier of $\sde^{z - \frac{1}{2}}$ and
$$ \sde^{z - \frac{1}{2}} \,  [(\sde^{1 - z} e_n)\, b^* a \,(\sde^\frac{1}{2} e_n)]
= (\sde^\frac{1}{2} e_n) \, b^* a \, (\sde^\frac{1}{2} e_n) \ . $$

So   $\sde^{z - \frac{1}{2}} \, [(\sde^{1 - z} e_n)\,b^* a\,(\sde^\frac{1}{2} e_n)]$ belongs
to $\Mfi$ by lemma \ref{we5.lem8}.

\medskip

We know on the other hand that $\sde^\frac{1}{2} e_n$ is a left multiplier of $\sde^{z -
\frac{1}{2}}$ and that $(\sde^\frac{1}{2} e_n) \, \sde^{z - \frac{1}{2}} = \sde^z e_n$.
This implies that $(\sde^{1-z} e_n) \, b^* a \, (\sde^\frac{1}{2} e_n)$ is a left
multiplier of $\sde^{z - \frac{1}{2}}$ and
$$ [(\sde^{1-z} e_n) \, b^* a \, (\sde^\frac{1}{2} e_n)] \, \sde^{z - \frac{1}{2}}
= (\sde^{1-z} e_n) \, b^* a \, (\sde^z e_n)  \ .$$

So $[(\sde^{1-z} e_n) \, b^* a \, (\sde^\frac{1}{2} e_n)] \, \sde^{z - \frac{1}{2}}$
belongs also to $\Mfi$ by \ref{we5.lem8}.

Therefore proposition \ref{we4.prop11} implies that
\begin{eqnarray*}
\vfi_\sde(e_n\,b^* a\,e_n) & = &  \vfi\bigl(\,(\sde^\frac{1}{2} e_n)\, b^* a\,(\sde^\frac{1}{2}
e_n)\,\bigr)
=  \vfi\bigl(\,\sde^{z - \frac{1}{2}} \, [(\sde^{1 - z} e_n)\, b^* a \,(\sde^\frac{1}{2} e_n)]
\,\bigr) \\
& = & \lambda^{i(\frac{1}{2}-z)} \, \vfi\bigl(\,[(\sde^{1-z} e_n)\, b^* a\, (\sde^\frac{1}{2}
e_n)] \, \sde^{z - \frac{1}{2}}\, \bigr) \\
& = & \lambda^{i(\frac{1}{2}-z)} \, \vfi\bigl(\,(\sde^{1-z} e_n)\, b^* a \,(\sde^z e_n)\,\bigr)
 \ .
\end{eqnarray*}
\end{demo}

\begin{proposition}
Consider $x \in \cM_{\vfi_\sde}$ and $z \in \C$ such that $x$ is a middle multiplier of
$\sde^{1-z},\,\sde^z$ and $\sde^{1-z} x \sde^z$ belongs to $\Mfi$. Then
$$\vfi_\sde(x) = \lambda^{i(\frac{1}{2}-z)} \, \vfi\bigl(\sde^{1-z} x \sde^z\bigr) \ .$$
\end{proposition}
\begin{demo}
Choose $n \in \N$. By the previous lemma, we know that $e_n \, x \, e_n$ belongs to
$\cM_{\vfi_\sde}$, $(\sde^{1-z} e_n) \,  x \, (\sde^ze_n)$ belongs to $\Mfi$ and
$$\vfi_\sde(e_n x e_n) = \lambda^{i(\frac{1}{2}-z)} \,
\vfi\bigl(\,(\sde^{1-z} e_n) x (\sde^z e_n)\,\bigr) \ .$$

Because $e_n$ is a left multiplier of $\sde^{1-z}$ and $e_n$ is a right multiplier of
$\sde^z$, we get that
$$e_n \, (\sde^{1-z} x \sde^z) \, e_n = (e_n \sde^{1-z})\, x (\sde^z \, e_n)
= (\sde^{1-z} e_n)\, x (\sde^z \, e_n) \ .$$
So $e_n \, (\sde^{1-z} x \sde^z) \, e_n$ belongs to $\cM_\vfi$ and
$$\vfi_\sde(e_n x e_n) = \lambda^{i(\frac{1}{2}-z)} \,
\vfi(e_n \, (\sde^{1-z} x \sde^z) \, e_n) \ . $$
The proposition follows by
lemmas \ref{we5.lem5} and \ref{we5.lem6}.
\end{demo}

\medskip

\begin{remark}\rm
We would like to mention the following special cases.
\begin{enumerate}
\item Consider $x \in \cM_{\vfi_\sde}$ such that $x$ is a middle multiplier of
$\sde^\frac{1}{2}, \, \sde^\frac{1}{2}$ and $\sde^\frac{1}{2} x \sde^\frac{1}{2}$ belongs
to $\Mfi$. Then $$\vfi_\sde(x) = \vfi\bigl(\sde^\frac{1}{2} x \sde^\frac{1}{2} \bigr) \
.$$
\item Consider $x \in \cM_{\vfi_\sde}$ such that $x$ is a left multiplier of $\sde$
and $x \sde$ belongs to $\Mfi$. Then $$\vfi_\sde(x) = \lambda^{-\frac{i}{2}} \,
\vfi(x \sde) \ . $$
\item Consider $x \in \cM_{\vfi_\sde}$ such that $x$ is a right multiplier of $\sde$
and $\sde x$ belongs to $\Mfi$. Then $$\vfi_\sde(x) = \lambda^\frac{i}{2} \, \vfi(\sde x)
\ .$$
\end{enumerate}
\end{remark}

\medskip

\begin{proposition}
Consider an  element $a \in A^+$ and $z \in \C$ such that $a$ is a middle multiplier of
$\sde^{1-z},\,\sde^z$ and  $\sde^{1-z} a \sde^z$ belongs to $\Mfi$. Then $a$ belongs to
${\cal M}_{\vfi_\sde}^+$.
\end{proposition}
\begin{demo}
Take $b \in A$ such that $b^* b = a$.

Choose $m \in \N$. Because $\sde^{z-1} e_m$ is a left multiplier of $\sde^{1-z}$ and
$\sde^{-z} e_m$ is a right multiplier of $\sde^z$, we have that
$$(\sde^{z-1} e_m) \, (\sde^{1-z} a \sde^z ) \, (\sde^{-z} e_m)
= ((\sde^{z-1} e_m)\,\sde^{1-z}) \, a \, (\sde^z  \, (\sde^{-z} e_m))
= e_m \, a \, e_m = (b e_m)^* (b e_m)  \ , $$
Because $\sde^{1-z} a \sde^z$ belongs to $\Mfi$, this equality  and lemma \ref{we5.lem9}
imply that $(b e_m)^* (b e_m)$ belongs to ${\cal M}_{\vfi_\sde}$. Hence $b e_m$ belongs to
$\cN_{\vfi_\sde}$.

\medskip

Take $k,l \in \N$. By the previous part, we already know that $e_l b^* b e_k$ belongs to
${\cal M}_{\vfi_\sde}$.

Because $e_l$ is a right multiplier of $\sde^{1-z}$ and $e_k$ is a left multiplier of
$\sde^z$, we have that $e_l b^* b e_k$ is a middle multiplier of $\sde^{1-z},\,\sde^z$ and
$$ \sde^{1-z} \, (e_l b^* b e_k) \, \sde^z = (\sde^{1-z} e_l) \, b^* b \, (e_k \, \sde^z)
= (\sde^{1-z} e_l) \, a \, (e_k \, \sde^z)$$
Knowing that $e_l$ is a left multiplier of $\sde^{1-z}$ and $e_n$ is a right multiplier of
$\sde^z$, we get moreover that
$$ e_l \, (\sde^{1-z} a \sde^z) \, e_k  = (e_l \, \sde^{1-z}) \, a \, (\sde^z e_k)
= (\sde^{1-z} e_l) \, a \, (e_k \sde^z )= \sde^{1-z} \, (e_l b^* b e_k) \, \sde^z \ .$$
Hence $\sde^{1-z} \, (e_l b^* b e_k) \, \sde^z$ belongs to $\Mfi$ by lemma \ref{we5.lem9}.

Therefore the previous proposition implies that
$$ \vfi_\sde( e_l  b^* b  e_k ) = \lambda^{i(\frac{1}{2}-z)}
\, \vfi\bigl(\sde^{1-z} \, (e_l b^* b e_k) \, \sde^z \bigl)
= \lambda^{i(\frac{1}{2}-z)} \, \vfi( e_l \, (\sde^{1-z} a \sde^z) \, e_k ) \ .$$

\medskip

Therefore, we have for all $m,n \in \N$ that
\begin{eqnarray*}
& & \hspace{-0.3cm} \| \ga(b e_m) - \ga(b e_n) \|^2 = \vfi_\sde(e_m \, b^* b \, e_m)
- \vfi_\sde(e_m \, b^* b \, e_n)  - \vfi_\sde(e_n \, b^* b \, e_m) +
\vfi_\sde(e_n \, b^* b \, e_n) \\
& & \hspace{-0.3cm} \spat = \lambda^{i(\frac{1}{2}-z)}
\,\, [\,\vfi(e_m \, (\sde^{1-z} a \sde^z) \, e_m)
- \vfi(e_m \, (\sde^{1-z} a \sde^z) \, e_n)
- \vfi(e_n \, (\sde^{1-z} a \sde^z) \, e_m)
 + \vfi(e_n \, (\sde^{1-z} a \sde^z) \, e_n)\,]
 \end{eqnarray*}

Using this equality and the obvious generalization of lemma \ref{we5.lem5}, we see that
$(\ga(b e_n))_{n=1}^\infty$ is Cauchy and hence convergent. Because we also have  that $(b
e_n)_{n=1}^\infty$ converges to $b$, the closedness of $\ga$ imply that $b$ belongs to
${\cal N}_{\vfi_\sde}$. Hence, $a$ is an element in ${\cal M}_{\vfi_\sde}$.
\end{demo}

\medskip

The next proposition is proven in a similar way.

\begin{proposition}
Consider an element $a \in {\cal M}_{\vfi_\sde}$ and $z \in \C$ such that $a$ is a middle
multiplier of $\sde^{1-z},\,\sde^z$ and $\sde^{1-z} a \sde^z$ belongs to $A^+$. Then
$\sde^{1-z} a \sde^z$ belongs to $\Mfi^+$.
\end{proposition}

\medskip

\begin{corollary}
Consider an element $a \in A^+$ such that $a$ is a middle multiplier of
$\sde^\frac{1}{2},\,\sde^\frac{1}{2}$ and $\sde^\frac{1}{2} a \sde^\frac{1}{2}$ belongs to
$A$. Then $a$ belongs to ${\cal M}_{\vfi_\sde}^+$ $\Leftrightarrow$ $\sde^\frac{1}{2} a
\sde^\frac{1}{2}$ belongs to $\Mfi^+$.
\end{corollary}

\bigskip\bigskip

In the  next part of this section, we will prove that the modular conjugation of
$\vfi_\sde$ with respect to $(H,\ga,\pi)$ is a scalar multiple of $J$.

\begin{notation}
We define  the set $D = \{\, x \in \Nfi \cap \Nfi^* \mid \text{We have for every } z \in
\C$

\vspace{2mm}
$\text{ that } $ \begin{minipage}[t]{14cm}
\begin{enumerate}
\item $x$ is a right multiplier of $\sde^z$ and $\sde^z x$ belongs to $\Nfi \cap \Nfi^*$
\item $x$ is a left multiplier of $\sde^z$ and $x \sde^z$ belongs to $\Nfi \cap \Nfi^*$
\hspace{3cm} $\}$ \ .
\end{enumerate} \end{minipage}
\end{notation}

\medskip

It is not difficult to check that $D$ is a sub-$^*$-algebra of
$N_0 \cap N_0^* \subseteq \cN_{\vfi_\sde} \cap \cN_{\vfi_\sde}^*$.

\medskip

\begin{lemma}
Consider $x \in \cN_{\vfi_\sde} \cap \cN_{\vfi_\sde}^*$. Then there exists a sequence
$(x_n)_{n=1}^\infty$ in  $D$ such that
\begin{enumerate}
\item $(x_n)_{n=1}^\infty$ converges to $x$
\item $(\ga(x_n))_{n=1}^\infty$ converges to $\ga(x)$
\item $(\ga(x_n^*))_{n=1}^\infty$ converges to $\ga(x^*)$
\end{enumerate}
\end{lemma}
\begin{demo}
We define for every $n \in \N$ the element $x_n = e_n x e_n$. Using lemma \ref{we5.lem3},
it is not difficult to check that $x_n$ belongs to $D$ for every $n \in \N$. We have also
for every $n \in \N$ that \begin{itemize}
\item $\ga(x_n) = \ga(e_n x e_n) =J \pi(\si_\frac{i}{2}'(e_n))^* J \pi(e_n) \ga(x)$
\item $\ga(x_n^*) = \ga(e_n x^* e_n) = J \pi(\si_\frac{i}{2}'(e_n))^* J \pi(e_n) \ga(x^*) $
\end{itemize}

This implies that $(x_n)_{n=1}^\infty$ converges to $x$,
$(\ga(x_n))_{n=1}^\infty$ converges to $\ga(x)$ and
$(\ga(x_n^*))_{n=1}^\infty$ converges to $\ga(x^*)$
\end{demo}

It follows immediately that $\ga(D)$ is dense in $H$.

\medskip

\begin{proposition}
We have that $\lambda^\frac{i}{4} \, J$ is the modular operator for
$\vfi_\sde$ in the GNS-construction $(H,\ga,\pi)$.
\end{proposition}
\begin{demo}
Denote the modular conjugation of $\vfi_\sde$ in the GNS-construction $(H,\ga,\pi)$ by $J'$.

\vspace{1mm}

Choose $a \in D$.
\begin{itemize}
\item We have that $a^*$ is a left multiplier of $\sde^\frac{1}{2}$  and  that
$a^* \sde^\frac{1}{2}$ belongs to $\Nfi$. This implies that $\ga(a^*) =
\la(a^* \sde^\frac{1}{2})$ by the definition of $\ga$. We know moreover
that $a \in D(\sde^\frac{1}{2})$ and that $\sde^\frac{1}{2} a$ belongs to $\Nfi \cap
\Nfi^*$. This implies that $\la(\sde^\frac{1}{2} a)$ belongs to $D(\nab^\frac{1}{2})$ and
$$J \nab^\frac{1}{2} \la(\sde^\frac{1}{2} a)
= \la( (\sde^\frac{1}{2} a)^*) = \la(a^* \sde^\frac{1}{2})
= \ga(a^*) \ .$$

\item Because $a$ belongs to $\Nfi$, $a$ belongs to $D(\sde^\frac{1}{2})$ and
$\sde^\frac{1}{2} a$ belongs to $\Nfi$, proposition \ref{we1.prop1}  implies  that
$\la(a)$ belongs to $D(\pi(\sde)^\frac{1}{2})$ and $\pi(\sde)^\frac{1}{2} \la(a) =
\la(\sde^\frac{1}{2} a)$. So, using the previous part, we get that $\la(a)$ belongs to
$D(\nab^\frac{1}{2} \, \pi(\sde)^\frac{1}{2} )$ and
$$ J \, \nab^\frac{1}{2} \, \pi(\sde)^\frac{1}{2}  \la(a)
= J \, \nab^\frac{1}{2} \la( \sde^\frac{1}{2} a ) = \ga(a^*) \ . $$

\item We know that $a$ is a left multiplier of $\sde^\frac{1}{2}$. This implies that
$a \sde^\frac{1}{2}$ is a left multiplier of $\sde^{-\frac{1}{2}}$ and that
$(a \sde^\frac{1}{2}) \, \sde^{-\frac{1}{2}} = a$ which belongs to $\Nfi$. Therefore,
proposition \ref{we4.prop8} implies that $\la(a \sde^\frac{1}{2})$ belongs to
$D((J \pi(\sde) J)^{-\frac{1}{2}})$ and
$$ (J \pi(\sde) J)^{-\frac{1}{2}} \la(a \sde^\frac{1}{2})
= \lambda^\frac{i}{4} \, \la((a \sde^\frac{1}{2}) \, \sde^{-\frac{1}{2}})
= \lambda^\frac{i}{4} \, \la(a) \ .$$
By definition, we have that $\ga(a) = \la(a \sde^\frac{1}{2})$. Consequently, we have that
$\ga(a)$ belongs to \newline $D((J \pi(\sde) J)^{-\frac{1}{2}})$ and
$$(J \pi(\sde) J)^{-\frac{1}{2}} \ga(a) = \lambda^\frac{i}{4} \, \la(a) \ .
$$
Using the previous part, we see that $\ga(a)$ belongs to   $D(\nab^\frac{1}{2} \,
\pi(\sde)^\frac{1}{2} \, (J \pi(\sde) J)^{-\frac{1}{2}})$ and
$$ J \nab^\frac{1}{2} \, \pi(\sde)^\frac{1}{2} \, (J  \pi(\sde) J)^{-\frac{1}{2}}
\ga(a) = \lambda^{-\frac{i}{4}} \, J \nab^\frac{1}{2} \, \pi(\sde)^\frac{1}{2} \la(a)
 = \lambda^{-\frac{i}{4}} \, \ga(a^*) \ .$$

\end{itemize}

The only thing we need from the previous discussion is that
$\ga(a)$ belongs to $D(\nab^\frac{1}{2} \, \pi(\sde)^\frac{1}{2} \,
(J  \pi(\sde) J)^{-\frac{1}{2}})$ and
$$ J \nab^\frac{1}{2} \, \pi(\sde)^\frac{1}{2} \, (J \pi(\sde) J)^{-\frac{1}{2}}
\ga(a) =  \lambda^{-\frac{i}{4}} \, \ga(a^*) \ . \text{\ \ \ \ \ \ (*)}$$

We have that $\pi(\sde)$ and $J \pi(\sde) J$ are strongly commuting (lemma
\ref{we5.lem1}). Define the positive injective operator $S$ in $H$ such that $S^{it} =
\pi(\sde)^{it} \, (J \pi(\sde) J)^{-it}$ for every $t \in \R$. Therefore,
$\pi(\sde)^\frac{1}{2} (J \pi(\sde) J)^{-\frac{1}{2}} \subseteq S^\frac{1}{2}$.

We also have that $\nab$ and $S$ are strongly commuting (lemma \ref{we5.lem2}). By
definition, we have that $\nabp^{it} = \nab^{it} S^{it}$ for $t \in \R$. This implies that
$\nab^\frac{1}{2} S^\frac{1}{2} \subseteq \nabp^\frac{1}{2}$.

Combining these two facts, we arrive at the conclusion that
$\nab^\frac{1}{2} \pi(\sde)^\frac{1}{2} (J \pi(\sde) J)^{-\frac{1}{2}}$
$\subseteq \nabp^\frac{1}{2}$.

Using (*), this gives that $\ga(a)$ belongs to $D(\nabp^\frac{1}{2})$
and $J \nabp^\frac{1}{2} \ga(a) = \lambda^{-\frac{i}{4}} \, \ga(a^*)$.

On the other hand, the fact that $a \in \cN_{\vfi_\sde} \cap \cN_{\vfi_\sde}^*$
implies that $J' \nabp^\frac{1}{2} \ga(a) = \ga(a^*)$.

Combining these two equalities results in the equality
$\lambda^\frac{i}{4} \, J \ga(a^*) = J' \ga(a^*)$.

Because $\ga(D)$ is dense in $H$, we get that $J' = \lambda^\frac{i}{4} \, J$.
\end{demo}

\bigskip

The following two lemmas will be used to prove that $\sde$ is uniquely determined in the
GNS-representation of $\vfi$ by the weight $\vfi_\sde$.

\begin{lemma}  \label{new.lem1}
Consider an element $a \in \Nfi \cap \cN_{\vfi_\sde}$. Then  $\la(a)$ belongs to $D(J
\pi(\sde^\frac{1}{2}) J)$ and $$\lambda^\frac{i}{4} \, J \pi(\sde^\frac{1}{2})
J \la(a) = \ga(a) \ .$$ Consequently,  $\|J \pi(\sde^\frac{1}{2})J \la(a)\|^2 =
\vfi_\sde(a^* a)$.
\end{lemma}
\begin{demo}
Choose $m \in \N$.

Because $a \in \Nfi$, we know by lemma \ref{we5.lem9} that $a e_m$ belongs to $\Nfi$.

We also know that $e_m$ is a left multiplier of $\sde^\frac{1}{2}$, which implies that $(a
\, e_m) \, \sde^\frac{1}{2}$ is a left multiplier of $\sde^\frac{1}{2}$ and $(a \, e_m)
\, \sde^\frac{1}{2} = a \,(e_m \, \sde^\frac{1}{2})$. This last equality
ensures that also $(a e_m) \sde^\frac{1}{2}$ belongs to $\Nfi$.

Hence, proposition  \ref{we4.prop8} implies that $\la(a e_m)$ belongs to $D(J
\pi(\sde)^\frac{1}{2} J)$ and $$\lambda^\frac{i}{4} \, J \pi(\sde^\frac{1}{2})
J \la(a e_m) = \la( (a e_m) \sde^\frac{1}{2})$$

So by the definition of $\ga$, we see that $a e_m$ belongs to $\cN_{\vfi_\sde}$ and
$$\lambda^\frac{i}{4} \, J \pi(\sde^\frac{1}{2}) J \la(a e_m) = \ga(a e_m)$$
Because $a \in \cN_{\vfi_\sde}$ and $e_m \in D(\overline{\si}_\frac{i}{2}')$, this implies
that
$$\lambda^\frac{i}{4} \, J \pi(\sde^\frac{1}{2}) J \la(a e_m)
= J \pi(\si_\frac{i}{2}'(e_m))^* J \ga(a)$$

Because $a \in \Nfi$ and $e_m \in D(\si_\frac{i}{2})$, we have also that
$$ \la(a e_m) = J \pi(\si_\frac{i}{2}(e_m))^* J \la(a) $$

\medskip

So we see that $(\la(a e_n))_{n=1}^\infty$ converges to $\la(a)$ and that
$\bigl(\,\lambda^\frac{i}{4} \, J \pi(\sde^\frac{1}{2}) J \la(a e_m)\,\bigr)_{n=1}^\infty$
converges to $\ga(a)$.

Therefore the closedness of $J \pi(\sde^\frac{1}{2}) J$ implies that $\la(a)$ belongs to
$D( J \pi(\sde^\frac{1}{2}) J)$ and
$$\lambda^\frac{i}{4} \, J \pi(\sde^\frac{1}{2}) J \la(a) = \ga(a) \ . $$
\end{demo}

\medskip

\begin{lemma}  \label{new.lem2}
The set $\la(\Nfi \cap \cN_{\vfi_\sde})$ is a core for $J \pi(\sde^\frac{1}{2}) J$.
\end{lemma}
\begin{demo}
Choose $m \in \N$.

We know already that $\si_\frac{i}{2}(e_m)^*$ is a left and right multiplier of
$\sde^\frac{1}{2}$ and that $\si_\frac{i}{2}(e_m)^* \, \sde^\frac{1}{2} =
\sde^\frac{1}{2} \, \si_\frac{i}{2}(e_m)^*$.

This implies easily that $\pi(\si_\frac{i}{2}(e_m))^* \, \pi(\sde^\frac{1}{2})
\subseteq \pi(\sde^\frac{1}{2}) \, \pi(\si_\frac{i}{2}(e_m))^* \in B(H)$.

So we get that $$(J \pi(\si_\frac{i}{2}(e_m))^* J) \, (J \pi(\sde^\frac{1}{2}) J)
\subseteq (J \pi(\sde^\frac{1}{2}) J) \, (J \pi(\si_\frac{i}{2}(e_m))^* J) \in B(H) \ .$$

\medskip

Because  $\bigl(\,J \pi(\si_\frac{i}{2}(e_n))^* J\,\bigr)_{n=1}^\infty$ is a bounded
sequence which converges strongly to 1 and  $\la(A)$ is dense in $H$, it is not too
difficult to infer from this last equality that the set
$$\langle \, J \pi(\si_\frac{i}{2}(e_m))^* J \la(a) \mid n \in \N, a \in \Nfi \, \rangle $$
is a core for $J \pi(\sde^\frac{1}{2}) J$.

\medskip

By lemma \ref{we5.lem9} and proposition \ref{we4.prop5}, we have for every $a \in \Nfi$
and $n \in \N$ that $a\,e_n$ belongs to $\Nfi \cap \cN_{\vfi_\sde}$ and $\la(a \, e_m) = J
\pi(\si_\frac{i}{2}(e_m))^* J \la(a)$. So we see that $\la(\Nfi \cap \cN_{\vfi_\sde})$ is
a core for $J \pi(\sde^\frac{1}{2}) J$.
\end{demo}

\bigskip

In a last proposition, we  will prove that the strictly positive element $\sde$ is
uniquely determined in the GNS-representation of $\vfi$ by the weight $\vfi_\sde$.

\begin{proposition}
Consider a \cst-algebra $A$ and a KMS-weight $\vfi$ on $A$ with modular group $\si$ and
GNS-construction $(H,\la,\pi)$. Let $\sde_1$, $\sde_2$ be strictly positive elements
affiliated with $A$ such that there exists strictly positive numbers $\lambda_1$,
$\lambda_2$ satisfying  $\si_t(\sde_1) = \lambda_1^t \, \sde_1$ and $\si_t(\sde_2) =
\lambda_2^t \, \sde_2$ for $t \in \R$  .

If $\vfi_{\sde_1} = \vfi_{\sde_2}$, then $\pi(\sde_1) = \pi(\sde_2)$.
\end{proposition}
\begin{demo}
Because $\vfi_{\sde_1} = \vfi_{\sde_2}$, we have that $\cN_{\vfi_{\sde_1}} =
\cN{\vfi_{\sde_2}}$ so $\la(\Nfi \cap \cN_{\vfi_{\sde_1}}) = \la(\Nfi \cap
\cN{\vfi_{\sde_2}})$.

By lemma \ref{new.lem2}, this gives a common core $\la(\Nfi \cap \cN_{\vfi_{\sde_1}})$ for
$J \pi(\sde_1^\frac{1}{2}) J$ and $J \pi(\sde_2^\frac{1}{2}) J$.

Using lemma \ref{new.lem1}, we have moreover for every $a \in \Nfi \cap
\cN_{\vfi_{\sde_1}}$ that
$$ \| J \pi(\sde_1^\frac{1}{2}) J \la(a) \|^2 = \vfi_{\sde_1}(a^* a)
= \vfi_{\sde_2}(a^* a) = \| J \pi(\sde_1^\frac{1}{2}) J \la(a) \|^2 \ . $$

This implies easily that $D(J \pi(\sde_1^\frac{1}{2}) J) = D(J \pi(\sde_2^\frac{1}{2}) J)$
and $\| J \pi(\sde_1^\frac{1}{2}) J v \| = \| J \pi(\sde_2^\frac{1}{2}) J v \|$ for every
$v \in D(J \pi(\sde_1^\frac{1}{2} J)$.

As a consequence, we get that $J \pi(\sde_1) J = J \pi(\sde_2) J$ which gives us that
$\pi(\sde_1) = \pi(\sde_2)$.
\end{demo}

\medskip

\begin{corollary}
Consider a \cst-algebra $A$ and a  faithful KMS-weight $\vfi$ on $A$ with modular group
$\si$. Let $\sde_1$, $\sde_2$ be strictly positive elements affiliated with $A$ such that
there exists strictly positive numbers $\lambda_1$, $\lambda_2$ satisfying  $\si_t(\sde_1)
= \lambda_1^t \, \sde_1$ and $\si_t(\sde_2) = \lambda_2^t \, \sde_2$ for $t \in \R$.
If $\vfi_{\sde_1} = \vfi_{\sde_2}$, then $\sde_1 = \sde_2$.
\end{corollary}

\end{document}